\documentclass[a4paper,12pt]{amsart}

\usepackage{latexsym}
\usepackage{amsfonts}
\usepackage{amsmath}
 \newtheorem{introtheo}{Th\'eor\`eme}

\newtheorem{lem}{Lemme}
\newtheorem{theo}[lem]{Th\'eor\`eme}
\newtheorem{prop}[lem]{Proposition}
\newtheorem{cor}[lem]{Corollaire}

\newtheorem{defin}[lem]{Definition}

\newcommand\preuve{\noindent {\bf Preuve.\ }}
\newcommand\supp{{\rm Supp}\,}

\def\SN{\mathbb N}    
\def\SP{\mathbb P}    
\def\SR{\mathbb R}   
\def\SZ{\mathbb Z}    

\def\Aa{{\mathcal A}}           
\def\Cc{{\mathcal C}}           
\def\Ee{{\mathcal E}}           
\def\Ff{{\mathcal F}}           
\def\Gg{{\mathcal G}}           
\def\Ss{{\mathcal S}}           

\renewcommand{\and}{et\ }
\def\datefrancais{%
\def\today{\ifnum\day=1\relax 1\/$^{\rm er}$\else
  \number\day\fi \space\ifcase\month\or
  janvier\or f\'evrier\or mars\or avril\or mai\or juin\or
  juillet\or ao\^ut\or septembre\or octobre\or novembre\or
  d\'ecembre\fi
  \space\number\year}}

\renewcommand{\and}{et\ }
\newcommand{\Mean}{-\hspace{-.45cm}\int}

\newcommand{\mean}{{\scriptstyle -}\hspace{-.33cm}\int}

\begin{document}


\protect\pagenumbering{roman}\protect

\title[Espaces critiques pour
Navier-Stokes]{Espaces critiques pour le syst\`eme
des equations de Navier-Stokes incompressibles}

\author{P. Auscher et  Ph. Tchamitchian} 
\email{auscher@u-picardie.fr \and
tchamphi@math.u-3mrs.fr}

\address{Universit\'e d'Amiens, Facult\'e de
math\'ematiques et d'Informatique,  
 33, rue
Saint Leu, F-80039 Amiens Cedex 1, et LAMFA, CNRS,
UPRES-A 6119\\ \and \\
Universit\'e d'Aix-Marseille III,
Facult\'e des
Sciences et Techniques de Saint-J\'er\^ome,
 Avenue Escadrille
Normandie-Niemen,
F-13397 Marseille Cedex 20,
et LATP, CNRS, UMR 6632. 
 }

\maketitle

\centerline{7 { mai} 1999}

\begin{abstract}
Nous d\'egageons 
dans ce travail des conditions abstraites portant
sur un espace fonctionnel E qui assurent
l'existence globale pour toute donn\'ee initiale
suffisamment petite dans E, ou l'exis\-ten\-ce
locale sans condition de taille, de solutions
pour une classe d'\'equa\-tions paraboliques
semi-lin\'eaires, dont le syst\`eme de
Navier-Stokes incompressibles dans l'espace
constitue un exemple fondamental. Nous donnons
\'egalement un crit\`ere abstrait de
r\'egularit\'e des solutions obtenues. Ces
conditions sont simples \`a v\'erifier dans tous
les cas connus: espaces de Lebesgue, de Lorentz,
de Besov, de Morrey, {\it et caetera}. Elles
s'adaptent au cas d'espaces E non invariants par
translation : nous d\'etaillons l'exemple de
certains espaces 2-microlocaux.
\end{abstract}

\maketitle

\bigskip

{\bf AMS Classification numbers:} 35K55, 35Q30,
35R05, 35S50, 42B25. 

\bigskip

{\bf Mots-clefs:} Navier-Stokes systems; mild solutions; Littlewood-Paley decomposition; maximal spaces



\clearpage
\tableofcontents


\section*{Introduction}

Les \'equations \'etudi\'ees ici sont des \'equations paraboliques
semi-lin\-\'e\-ai\-res, \'ecrites sous la forme int\'egrale
\begin{equation}\label{eq 1}
  u = Su_{0} + B(u,u),
 \end{equation}
dont l'inconnue, not\'ee
$u$, est une distribution temp\'er\'ee d\'efinie sur $]0,\infty[\times
\SR^3$. On d\'esigne par $u_{0}$ une distribution donn\'ee de
$\Ss'(\SR^3)$, et par $Su_{0}$ l'image de $u_{0}$ sous l'action du
semi-groupe de la chaleur~: $Su_{0} (t) = e^{t\Delta}u_{0}, t > 0$.
Enfin, $B$ est une application bilin\'eaire sym\'etrique, formellement
d\'efinie par
\begin{equation}\label{eq  2}
 B(u,v) (t) = \int^t_{0} e^{(t - \tau) \Delta} P(D) \{u(\tau)
v(\tau)\} d\tau,
\end{equation}
o\`u $P(D)$ est un op\'erateur  pseudo-diff\'erentiel homog\`ene de degr\'e 1,
dont le symbole  $P(\xi)$ est suppos\'e non nul et $C^\infty$ en-dehors
de 0.\medskip

Bien que nous restreignant au cas scalaire, nous pourrions sans
difficult\'e consid\'erer le cas vectoriel~: le syst\`eme de
Navier-Stokes incompressible est alors un exemple fondamental,
d'ailleurs \`a l'origine de ce travail, dans lequel $u$ et $u_{0}$ ont
trois composantes scalaires et sont \`a divergence nulle, le terme
bilin\'eaire s'\'ecrivant
$$B(u,v) (t) = -\frac{1}{2}\int^t_{0}\ e^{(t -
\tau)\Delta}\ \SP
\nabla\ .\ (u(\tau)\otimes v(\tau) + v(\tau)\otimes u(\tau))\ d\tau,$$
o\`u $\SP$ est le projecteur de Leray dans $\SR^3$.\medskip

Reprenant dans [K] un sch\'ema formalis\'e par Weissler ([W]),
Kato r\'esoud (\ref{eq 1})\footnote{Plus
exactement, Kato r\'esoud le syst\`eme de
Navier-Stokes incompressible. Mais, une fois le
formalisme mis en place, la condition
d'incompressibilit\'e ne joue plus aucun r\^ole,
et il est plus simple de ne consid\'erer que le
cas scalaire g\'en\'eral.} pour toute donn\'ee $u_{0}\in
L^3(\SR^3)$ suffisament petite, et obtient des
solutions dans
${\mathcal C} ([0,\infty[ ; L^3 (\SR^3))$
\footnote{Dans tout l'article,  l'espace  $ {\mathcal C} (I ; E)$ des
fonctions continues de $I$ \`a valeurs dans $E$ est muni
implicitement de la norme $\sup_{t\in I} \|u(t)\|_E$.}, alors
que l'application $B$ n'est pas continue sur cet espace (voir Oru
\cite{bib O}). Sa m\'ethode repose sur la d\'efinition d'un espace de
Banach $\Ff$, inclus dans ${\mathcal C} ([0,\infty[ ; L^3 (\SR^3))$,
tel que\par
\begin{enumerate}
\item[i)] $Su_{0}\in \Ff$ si $u_{0} \in
L^3(\SR^3),$
\item[ii)]  $B$ soit continue de $\Ff\times \Ff$ 
dans $\Ff$.
\end{enumerate}
 Il ne lui reste plus qu'\`a utiliser les
it\'erations successives de Picard pour
r\'esoudre (\ref{eq 1}), et c'est de cette
derni\`ere \'etape que provient la condition de
taille sur
$u_{0}$.\medskip

La m\^eme m\'ethode lui a \'egalement permis d'obtenir des solutions
locales pour toute donn\'ee $u_{0}$, \`a partir d'un espace $\Ff_{T}$
inclus dans $\Cc ([0,T[ ; L^3 (\SR^3))$ v\'erifiant les analogues de
i) et ii). Dans ce cas, la condition de taille porte sur $T =
T(u_{0})$.\medskip

Cette m\'ethode, que nous appellerons dor\'enavant m\'ethode KW, a \'et\'e
relue et d\'evelopp\'ee par Giga et Miyakawa ([G,M]), Taylor ([T]),
Kozono et Yamazaki ([Ko,Y]), Cannone (\cite{bib C}),
Meyer ([M]), Planchon ([P]), Barraza ([B]), et d'autres auteurs. De
nouvelles solutions de (\ref{eq 1})
ont \'et\'e obtenues (notamment des solutions autosimilaires ou
asymptotiquement autosimilaires) en rempla\c cant $L^3 (\SR^3)$ par
d'autres espaces $E$ bien choisis : l'espace de Lorentz
$L^{3,\infty}$ (Barraza, Meyer), les espaces de Besov $\dot B^{-1 +
\frac{3}{p},\infty}_{p}$, $p<\infty$ (Cannone, Planchon),
les espaces de Morrey (Giga et Miyakawa, Taylor, Cannone,
Lemari\'e), certains espaces de Besov plac\'es au-dessus
d'espaces de Morrey (Kozono et Yamazaki). Dans chaque cas, la
construction suit la d\'emarche de Kato, et devient sp\'ecifique dans
le choix de $\Ff$ et la preuve de la continuit\'e de $B$.\medskip

La question principalement \'etudi\'ee ici est la suivante : quels
sont les espaces $E$ pour lesquels la m\'ethode KW fonctionne?\medskip

Autrement dit : peut-on caract\'eriser les espaces $E$ pour lesquels
on peut trouver un espace $\Ff$ de fonctions continues du temps \`a valeurs banachiques  satisfaisant aux condition i) et ii),
conduisant ainsi \`a l'existence d'une solution de (1) pour toute
donn\'ee $u_{0}\in E$ assez petite?\medskip

Cela n\'ecessite de pr\'eciser la formalisation de la m\'ethode KW : tel
est l'objet
de la premi\`ere partie.\medskip

La seconde partie aborde le coeur du probl\`eme. On commence par
d\'elimiter la classe des espaces $E$ consid\'er\'es, d'abord en se
restreignant aux espaces, dits invariants, sur lesquels le groupe
affine $ax + b$ agit en accord avec les propri\'et\'es d'invariance
de l'ensemble des solutions de (1). On d\'efinit ensuite la notion de
compatibilit\'e avec la non-lin\'earit\'e contenue dans $B$. Celle-ci
exprime que chaque bloc de la d\'ecomposition de Littlewood-Paley d'un
produit $fg$ appartient \`a $E$, lorsque $f,g\in E$ et sont soumis \`a
des conditions spectrales. Plus pr\'ecis\'ement, on suppose
l'existence d'une suite $(\eta_{n})_{n\in \SZ}$ telle que, si
$\supp
\widehat f\subset \Gamma_{k} = \{\xi ; 2^{k - 1}\le \vert \xi \vert
\le 2^{k + 1}\}$ et $\supp \widehat g\subset \Gamma_{l}, k,l\in
\SZ$, alors on a pour tout $j\in  \SZ$
$$\Vert \Delta_{j} (fg)\Vert_{E}\le
\eta_{\max(k-j, l - j)}\ 2^{k + l-j}\Vert
f\Vert_{E} \Vert g\Vert_{E}.$$
(On renvoie \`a la section 2.1 pour les notations
et les
\'enonc\'es pr\'ecis.) Cette hypoth\`ese n'est pas contraignante,
dans la mesure o\`u tous les espaces invariants connus la v\'erifient.
Il faut d'ailleurs souligner que cette condition ne suppose pas que
l'appartenance \`a $E$ soit caract\'eris\'ee par une propri\'et\'e de
la d\'ecomposition de Littlewood-Paley.

C'est le comportement de la suite $(\eta_{n})$ qui discrimine les
espaces, ainsi que le montrent les trois principaux r\'esultats de ce
travail, qui peuvent \^etre r\'esum\'es de la fa\c con suivante
(rapide, mais impr\'ecise).\medskip

\begin{introtheo} La  m\'ethode  KW  fonctionne  pour  l'espace\ \
$E$\ \ d\`es que $\displaystyle{\sum_{n\ge 0} \eta_{n} <
+\infty}$.\end{introtheo} Cela signifie qu'on peut trouver $\Ff$
(respectivement $\Ff_{T}$) v\'erifiant i) et ii). En
revanche, l'inclusion de $\Ff$ dans $\Cc([0,\infty[;E)$
(respectivement $\Ff_{T}$ dans $\Cc ([0,T[ ; E)$) n'est pas a
priori satisfaite.\medskip

\begin{introtheo} Si de plus $\displaystyle \sum_{n\ge 0}
n \eta_{n} < +\infty$, on peut faire en sorte que $\Ff$ (resp. $\Ff_{T}$)
soit inclus
dans $\Cc ([0,\infty[;E)$ (resp. $\Ff_{T}$ dans $\Cc ([0,T[ ; E)$) si $E$
est s\'eparable.\end{introtheo}
Il existe aussi un r\'esultat de r\'egularit\'e analogue si $E$ est
le dual non s\'eparable d'un espace de Banach s\'eparable.\medskip

La d\'emonstration de ces deux r\'esultats est constructive.\medskip

\begin{introtheo} Si $(\eta_n)$ est d\'ecroissante, $(2^{3n/2}\eta_n)$ croissante \`a partir d'un certain rang et $\displaystyle \sum_{n\ge 0} \eta^{2}_{n} = +\infty$, alors il
existe un espace $E$ invariant et compatible, associ\'e \`a
$(\eta_{n})$, pour lequel il est impossible de faire fonctionner la
m\'ethode KW : cela signifie qu'il n'existe pas d'espace de fonctions
continues du temps $\Ff$ ou $\Ff_{T}$ v\'erifiant i) et
ii). \end{introtheo}
Les contre-exemples de ce dernier th\'eor\`eme comprennent
l'espace de Besov $\dot B^{-1,\infty}_{\infty}$. Voir la section 2.2 pour une condition plus g\'en\'erale impos\'ee \`a $(\eta_n)$. \medskip

Les th\'eor\`emes A et B sont d\'emontr\'es dans cette m\^eme partie,
le th\'eor\`eme C dans la quatri\`eme partie. Auparavant, on montre
dans la troisi\`eme partie que tous les exemples connus rel\`event des
th\'eor\`emes A et B.\medskip

Enfin, dans la cinqui\`eme et derni\`ere partie, on d\'ecrit en
d\'etail de nouveaux exemples, inspir\'es directement des espaces
2-microlocaux de Bony. La motivation est ici de construire des
solutions de (\ref{eq 1}) pour des donn\'ees initiales $u_{0}$ les plus
singuli\`eres possibles. On est amen\'e \`a d\'efinir des espaces de
distributions singuli\`eres sur un ferm\'e de $\SR^3$, pour lequel on
suppose qu'une sorte de densit\'e locale, appel\'ee fonction de
densit\'e, ob\'eit \`a la condition de Dini. On montre alors que de
tels espaces entrent dans le cadre du th\'eor\`eme A (convenablement
g\'en\'eralis\'e au cas d'espaces non invariants) et permettent donc
de faire fonctionner la m\'ethode KW. En particulier, on obtient
ainsi de nouvelles solutions autosimilaires.\medskip

Afin de simplifier l'exposition, nous ne consid\'erons que
l'existence et la r\'egularit\'e globales de solutions de (1), sauf
dans la section 2.6, qui d\'ecrit les modifications \`a apporter pour
obtenir les r\'esultats locaux\footnote{Ces r\'esultats r\'epondent \`a une question qui nous a \'et\'e pos\'ee par J.-Y. Chemin.}.\medskip

Les m\'ethodes d\'ecrites dans ce travail ne sont pas sp\'ecifiques \`a
la dimension 3, et s'adaptent en toute dimension.

\clearpage
\pagenumbering{arabic}

\section{Formalisation abstraite}\bigskip

On se donne une fois pour toutes un op\'erateur bilin\'eaire $B$ de la
forme (\ref{eq 2}). Le but de cette partie est de d\'ecrire de fa\c con
abstraite la m\'ethode KW. Celle-ci permet d'obtenir des
r\'esultats du type ``il existe $\alpha > 0$ tel que, pour toute
donn\'ee $u_{0}\in E$ v\'erifiant $\Vert u_{0}\Vert_{E} < \alpha$,
il existe $u\in \Ff$ solution de (\ref{eq 1}) et telle que $u(t)$ tende
vers $u_{0}$ en un certain sens, lorsque $t$ tend vers 0''.\bigskip

\subsection{Propri\'et\'e d'invariance}$\phantom{x}$\medskip

L'homog\'en\'e\"{\i}t\'e de l'op\'erateur $P(D)$
implique que, si
$u(t,x)$ est une solution de (\ref{eq 1}) pour la donn\'ee $u_{0} (x)$,
alors quels que soient $x_{0}\in \SR^3$ et $\lambda > 0,\ \lambda
u(\lambda^2t, \lambda x - x_{0})$ est \'egalement solution de \ref{eq
1}, pour la donn\'ee $\lambda u_{0} (\lambda x - x_{0})$. Par
cons\'equent, il est raisonnable de se restreindre aux espaces $E$
et $\Ff$ tels que
\begin{equation}\label{eq  3}
\forall x_{0}\in \SR^3\ \ \forall \lambda > 0\ \ \ \ \ \Vert \lambda
u_{0} (\lambda \cdot - x_{0})\Vert_{E} = \Vert u_{0}\Vert_{E},
\end{equation}
\begin{equation}\label{eq  4}
\forall x_{0}\in \SR^3\ \ \forall \lambda > 0\ \ \ \ \ \Vert \lambda
u (\lambda^{2}\cdot, \lambda \cdot- x_{0})\Vert_{\Ff} = \Vert
u\Vert_{\Ff}.
\end{equation}

Par exemple, si on cherche $E$ parmi les espaces de Lebesgue
$L^p(\SR^3)$, alors la valeur $p = 3$ est naturelle. En effet, si la
m\'ethode de Kato fonctionnait pour $L^p(\SR^3)$ avec $p\not= 3$, la
condition de taille sur $u_{0}$ pourrait \^etre supprim\'ee : en
choisissant bien $\lambda$, on aurait $\Vert v_{0}\Vert_{L^p} <
\alpha$, o\`u $v_{0} (x) = \lambda u_{0} (\lambda x)$, ce qui
donnerait une solution globale $u(t,x) = \frac{1}{\lambda}\  v\ 
\big(\frac{t}{ \lambda^2}, \frac{x}{  \lambda}\big)$ associ\'ee \`a
$u_{0}$, quelle que soit $u_{0}$. On ne sait pas si un tel r\'esultat
est vrai ou faux.

\begin{defin} On dit que $E$ est invariant s'il v\'erifie (\ref{eq
3}).\end{defin}\bigskip

\subsection{Couples admissibles}$\phantom{x}$\medskip

La m\'ethode KW fonctionne avec deux espaces de Banach $E$ et
$F$ ayant les propri\'et\'es $(P1), (P2)$ et $(P3)$ suivantes.\medskip

\noindent {\it Propri\'et\'e (P1)}\par
\begin{itemize}
\item[a)] $E$ s'injecte contin\^ument dans $\Ss'$,
\item[b)] $E$ est invariant.
\end{itemize}\medskip

\noindent {\it Propri\'et\'e (P2)}
\begin{itemize}
\item[a)] $F$ s'injecte contin\^ument dans $\Ss'$,
\item[b)] la norme de $F$ est invariante par translation, et
$f(\lambda \cdot)\in F$ si et seulement si $f\in F$, pour tout $\lambda >
0$,
\item[c)] sur tout compact de $]0,\infty[$, les normes
$$\Vert f\Vert_{t,F} = \sqrt t \Vert f(\sqrt {t} \, 
\cdot)\Vert_{F}$$ sont uniform\'ement \'equivalentes entre elles,
\footnote{Cette propri\'et\'e est v\'erifi\'ee d\`es
que $\displaystyle \lim_{\lambda\rightarrow 1} f(\lambda \, \cdot) = f$
pour tout $f\in F$.}
\item[d)] l'op\'erateur $e^\Delta$ est continu de $E$ dans
$F$.
\end{itemize}\medskip

\noindent {\it Propri\'et\'e (P3)}\par
\noindent Si $\Ff$ d\'esigne l'ensemble des fonctions continues de
$]0,\infty[$ \`a valeurs dans $F$, not\'ees $u$ ou $u(t)$, $t > 0$, telles que
$$\Vert u\Vert_{\Ff} = \sup_{t>0} \Vert u(t)\Vert_{t,F} < +\infty,$$
alors\par
\begin{itemize}
\item[a)] $B$ est continue de $\Ff\times \Ff$ dans $\Ff$,
\item[b)] $\displaystyle \lim_{t\rightarrow 0} B(u,v) (t) = 0$
dans $\Ss'$, pour tous $u,v\in \Ff$.
\end{itemize}

\noindent On note $\Vert B\Vert$ la plus petite des constantes $C$
telles que
\begin{equation}\label{eq bicontinuite}
\forall u,v\in \Ff\ \ \ \ \ \Vert B(u,v)\Vert_{\Ff}\le C \Vert
u\Vert_{\Ff}\Vert v\Vert_{\Ff}.
\end{equation}

\begin{defin} Lorsque les propri\'et\'es ci-dessus sont satisfaites,
le couple $(E,F)$ est dit admissible.\end{defin}

De tels espaces fournissent des solutions \`a (\ref{eq 1}) en vertu
du r\'esultat abstrait suivant.

\begin{theo} Soit $\Ff$ un espace de Banach et $B$ un op\'erateur
bilin\'e\-ai\-re continu de $\Ff\times \Ff$ dans $\Ff$. Alors, pour tout
$a\in \Ff$ tel que $\Vert a\Vert_{\Ff}\le \frac{1}{ 4\Vert B\Vert}$,
il existe $u\in \Ff$ solution de l'\'equation
\begin{equation}\label{eq 5}
u = a + B(u,u).
\end{equation}
De plus, il existe des op\'erateurs $T_{k}, k\ge 1$, tels que :

\begin{itemize}
\item[i)] chaque $T_{k}$ est la restriction \`a la diagonale de $\Ff^k$
d'un op\'erateur $k$-lin\'eaire continu de $\Ff^k$ dans $\Ff$;\par
\item[ii)] il existe une constante absolue $C$ telle que, pour tous $k\ge
1$ et $a\in \Ff$
$$\Vert T_{k} (a)\Vert_{\Ff}\le \frac{C}{ \Vert B\Vert} k^{-3/2} (4\
\Vert B\Vert\ \ \Vert a\Vert_{\Ff})^k ;$$
\item[iii)] si $\Vert a\Vert_{\Ff}\le \frac{1}{ 4\Vert B\Vert}$, alors
\begin{equation}\label{eq 6}
u = \displaystyle \sum^\infty_{k = 1} T_{k} (a).
\end{equation}
\end{itemize}
Enfin, on a toujours $\Vert u\Vert_{\Ff}\le \frac{1}{ 2\Vert B\Vert}$,
et $u$ est l'unique solution de (\ref{eq 5}) dans la boule ferm\'ee
$\overline{B_{\Ff}} (0,\frac{1}{ 2\Vert B\Vert})$.\end{theo}
\medskip

L'existence de $u$ est bien connue sous la condition $\Vert
a\Vert_{\Ff} < \frac{1}{ 4\Vert B\Vert}$~: voir par exemple Cannone
(\cite{bib C}, p.37), qui utilise le th\'eor\`eme des contractions de
Picard pour l'obtenir. L'approche par d\'eveloppement multilin\'eaire
de la solution qui est choisie ici est un peu plus pr\'ecise.\medskip

\preuve On d\'efinit les op\'erateurs $T_{k}$ de
proche en proche par les relations
\begin{equation}\label{eq 7}
T_{1} (a) = a,
\end{equation}
\begin{equation}\label{eq 8}
T_{k} (a) = \displaystyle \sum^{k - 1}_{l = 1} B(T_{l}(a), T_{k -
l}(a)),\ k\ge 2.
\end{equation}
Par construction et d'apr\`es $(P3)$, les $T_{k}$ sont la restriction
\`a la diagonale de $\Ff^k$ d'op\'erateurs $k$-lin\'eaires, et il
existe des constantes $a_{k}$ telles que
$$\Vert T_{k}(a)\Vert_{\Ff}\le a_{k} \Vert a\Vert^k_{\Ff}$$
pour tout $a\in \Ff$.\medskip

Pour estimer les $a_{k}$ on part de l'in\'egalit\'e de r\'ecurrence
$$a_{k}\le \Vert B\Vert \sum^{k - 1}_{l = 1} a_{l}\ a_{k - l},$$
avec la condition initiale $a_{1} = 1$. On en d\'eduit que, pour
tout $k\ge 1$
\begin{equation}\label{eq 9}
a_{k}\le \Vert B\Vert^{k - 1}\ c_{k},
\end{equation}
o\`u les coefficients $c_{k}$ sont tels que
\[\left\{\begin{array}{llll}
c_{1} &= &1\\
c_{k} &= &\displaystyle\sum^{k - 1}_{l = 1} c_{l} c_{k - l}, &k\ge 2.
\end{array}
\right. \]
Ce sont les nombres de Catalan (voir par exemple Comtet, tome
1 [Co]), donn\'es par la formule
$$c_{k} = \frac{(2k - 2) !}{ k ! (k - 1) !}.$$
Leur s\'erie g\'en\'eratrice se calcule ais\'ement :
\begin{equation}\label{eq 10}
\displaystyle \sum^\infty_{k = 1} c_{k} z^k = \frac{1 - \sqrt{1 -
4z}}{2},
\end{equation}
et lorsque $k$ tend vers $\infty$, on a
$$c_{k} \sim\frac{1}{ 4\sqrt \pi} k^{-3/2} 4^k.$$
Cela d\'emontre les points i) et ii) pour les op\'erateurs $T_{k}$
d\'efinis en (\ref{eq 7}) et (\ref{eq 8}).\medskip

Si $\Vert a\Vert_{\Ff}\le \frac{1}{ 4\Vert B\Vert}$, on pose $u =
\displaystyle \sum^\infty_{k = 1} T_{k} (a)$ : cette s\'erie est
normalement convergente dans $\Ff$. On calcule
\begin{align*}
B(u,u)&= \sum^\infty_{l = 1}\sum^\infty_{m = 1} B
(T_{l}(a), T_{m}(a))\\
&= \sum^\infty_{k = 2} \sum^{k - 1}_{l = 1} B(T_{l}(a),
T_{k - l}(a))\\
&= u - a
\end{align*}
d'apr\`es (\ref{eq 7} - \ref{eq 8}), ce qui montre que $u$ est une
solution de (\ref{eq 5}). De plus, il r\'esulte de (\ref{eq 9} - \ref{eq
10}) que
\begin{equation}\label{eq 11}
\Vert u\Vert_{\Ff}\le \frac{1 - \sqrt {1 - 4\Vert B\Vert \Vert
a\Vert_{\Ff}}}{2\Vert B\Vert}.
\end{equation}
En particulier, $\Vert u\Vert_{\Ff}\le \frac{1}{ 2\Vert
B\Vert}.$\medskip

Il reste \`a obtenir l'unicit\'e de $u$ dans la boule ferm\'ee
$\overline{B_{\Ff}}\big(0,\frac{1}{ 2\Vert B\Vert}\big)$. Elle est facile
\`a d\'emontrer lorsque $\Vert u\Vert_{\Ff} < \frac{1}{ 2\Vert B\Vert}$ :
si $v$ est une solution de (\ref{eq 5}), $\Vert v\Vert_{\Ff} \le
\frac{1}{ 2\Vert B\Vert}$, on a
$$u - v = B(u + v, u - v),$$
d'o\`u $u = v$, puisque $\Vert B\Vert\ \Vert u + v\Vert_{\Ff} < 1$.
Dans le cas g\'en\'eral, o\`u il est possible d'avoir $\Vert
u\Vert_{\Ff} = \frac{1}{2\Vert B\Vert}$, on d\'efinit pour tout $N\ge
1$ l'\'el\'ement $v_{N}$ de $\Ff$ par
\begin{equation}\label{eq 12}
v = T_{1}(a) + \dots + T_{N}(a) + v_{N},
\end{equation}
et on prouve que
\begin{equation}\label{eq 13}
\Vert v_{N}\Vert_{\Ff}\le \frac{1}{ \Vert B\Vert}
\sum^\infty_{k = N + 1}\ c_{k}\ 4^{-k}.
\end{equation}
Faisant tendre $N$ vers $\infty$, il vient $v = u$.\medskip

L'in\'egalit\'e (\ref{eq 13}) est vraie si $N = 1$, car d'une part
$v_{1} = B(v,v)$, donc $\Vert v_{1}\Vert_{\Ff}\le \Vert B\Vert\ \Vert
v\Vert^2_{\Ff}\le \frac {1}{ 4\Vert B\Vert}$, et d'autre part
$\displaystyle \sum^\infty_{k = 2}\ c_{k}\ 4^{-k} = \dfrac{1}{4}$
d'apr\`es (\ref{eq 10}).\medskip

Supposant (\ref{eq 13}) prouv\'ee au rang $N$, on injecte (\ref{eq 12})
dans l'\'equation (\ref{eq 5}), et on trouve
\begin{multline*}
v_{N + 1} =  \sum_{\substack{ l + m\ge N + 2\\
 1\le l\, , \, m\le N}}\ \ B(T_{l}(a), T_{m}(a)) \\ 
+ 2\ B
(v_{N}, T_{1}(a) + \dots + T_{N}(a)) + B(v_{N},v_{N}).
\end{multline*}
Puisqu'on sait que $\Vert T_{k}(a)\Vert_{\Ff}\le \frac{1}{ \Vert B\Vert}
c_{k} 4^{-k}$ pour tout $k\ge 1$, on obtient en utilisant
l'hypoth\`ese de r\'ecurrence
\begin{align*} \Vert v_{N + 1}\Vert_{\Ff}&\le \frac{1}{ \Vert
B\Vert}
\Biggl\{  \sum_{\substack{ l + m\ge N + 2\\
 1\le l\, , \, m\le N}}\  
c_{l}\ c_{m}\ 4^{-l-m}\\ &\ \ \ \ \ \ \
\ +2\bigg(
\sum^N_{k = 1}\ c_{k}\ 4^{-k}\bigg) \bigg(
\sum^\infty_{k = N + 1}\ c_{k}\ 4^{-k}\bigg)\\
&\ \ \ \ \ \ \ \ + \bigg( \sum^\infty_{k = N + 1}\ c_{k}\
4^{-k}\bigg)^2\Biggr\}\\
&\le\frac{1}{ \Vert
B\Vert} \sum_{k\ge N + 2} \bigg(
\sum^{k - 1}_{l = 1}\ c_{l}\
c_{k - l}\bigg)\ 4^{-k}\\
&\le\frac{1}{ \Vert
B\Vert} \sum_{k\ge N + 2}\ c_{k}\ 4^{-k}.
\end{align*}
Ainsi, (\ref{eq 13}) est prouv\'ee pour tout $N$, et la d\'emonstration
est achev\'ee.\qed\medskip

Ce th\'eor\`eme s'interpr\`ete comme un r\'esultat d'analyticit\'e au
voisinage de $0\ :\ u = 0$ est \'evidemment une solution de (\ref {eq
5}) lorsque $a = 0$, et les solutions $u$ construites lorsque $a$ est
petit sont obtenues par perturbation et d\'eveloppement en s\'erie
autour de 0. Elles d\'ependent analytiquement de $a$, pour la
topologie forte de $\Ff$. Enfin, par un ph\'enom\`ene analogue au
prolongement continu jusqu'au bord des s\'eries enti\`eres \`a
coefficients positifs, on peut r\'esoudre (\ref{eq 5}) sous la
condition limite $\Vert a\Vert_{\Ff} = \frac{1}{ 4\Vert B\Vert}$. Le
r\'esultat est optimal, comme le montre l'exemple \'el\'ementaire
$\Ff = \SR$ et $B(u,u) = u^2$.\medskip

Revenant \`a l'\'equation (\ref{eq 1}), on obtient le corollaire suivant.

\begin{prop} Si $(E,F)$ est un couple admissible, il existe $\alpha >
0$ tel que, pour tout $u_{0}\in E$ avec $\Vert u_{0}\Vert_{E}\le
\alpha$, l'\'equation (\ref {eq 1}) admet une solution $u\in \Ff$
telle que
\begin{equation}\label{eq 14}
\displaystyle \lim_{t\rightarrow 0}u(t) = u_{0}
\end{equation}
dans $\Ss'$. De plus, $u$ s'\'ecrit
$$u = \sum^\infty_{k = 1} T_{k} (Su_{0}),$$
o\`u les op\'erateurs $T_{k}$ sont donn\'es par le th\'eor\`eme 3.
Enfin, $\Vert u\Vert_{\Ff}\le \frac{1}{ 2\Vert B\Vert}$, et $u$ est
l'unique solution de (\ref{eq 1}) dans la boule ferm\'ee
$\overline{B_{\Ff}}(0,\frac{1}{ 2\Vert B\Vert})$.
\end{prop}

\preuve Soit $(E,F)$ un couple admissible et
$u_{0}\in E$. La continuit\'e de $e^\Delta$ de $E$ vers $F$ et la
d\'efinition de $\Ff$ impliquent $e^{t\Delta} u_{0}\in F$ pour tout
$t > 0$, avec
$$\Vert e^{t\Delta} u_{0}\Vert_{t,F}\le \Vert e^\Delta \Vert_{F,E}
\Vert u_{0}\Vert_{E}.$$
Ecrivant $\Delta e^\Delta = \Delta e^{\Delta/2} e^{\Delta/2}$, on
obtient de m\^eme $\Delta e^{t\Delta} u_{0}\in F$, avec
$$\Vert \Delta e^{t\Delta} u_{0}\Vert_{t,F}\le \frac{c}{ t} \Vert
u_{0}\Vert_{E}.$$
Puisque les normes $\Vert .\Vert_{t,F}$ sont uniform\'ement
\'equivalentes sur tout compact de $]0,\infty[$, ceci implique la
d\'erivabilit\'e de $t\longmapsto e^{t\Delta} u_{0}$, de $]0,\infty[$
dans $F$, et a fortiori la continuit\'e. On a donc $Su_{0}\in \Ff$, et
$$\Vert Su_{0}\Vert_{\Ff}\le \Vert e^\Delta\Vert_{F,E} \Vert
u_{0}\Vert_{E}.$$
Appliquant le th\'eor\`eme 3 avec\ \ $a = Su_{0}$,\ on peut r\'esoudre
(\ref{eq 1}) d\`es que\ \ $4\Vert B\Vert\ \Vert e^\Delta\Vert_{F,E}
\Vert u_{0}\Vert_{E}\le 1$. La relation (\ref{eq 14}) provient de
(P3b) et du fait que $E$ s'injecte contin\^ument dans $\Ss'$. Le reste
d\'ecoule directement du th\'eor\`eme 3.\qed\vskip 1 cm

On peut maintenant formuler pr\'ecis\'ement la question centrale
de ce travail : quels sont les espaces $E$ pour lesquels
on peut construire un espace $F$ formant avec $E$ un couple $(E,F)$
admissible ?\medskip

On y r\'epond dans la partie suivante, avant de donner divers exemples.



\section{Les bons espaces de Banach}\bigskip

\subsection{Espaces fonctionnels compatibles avec la
non-lin\'earit\'e}$\phantom{x}$\medskip

On d\'ecrit dans ce paragraphe les hypoth\`eses faites a priori sur
les espaces consid\'er\'es. Celle-ci ont pour but de permettre
l'utilisation de la d\'ecomposition de Littlewood-Paley, et de
l'algorithme associ\'e pour calculer un produit, qu'on rappelle
maintenant bri\`evement.\medskip

On se donne une fonction $\varphi^0$ de classe $C^\infty$ sur
$\SR^+$, avec $\varphi^0 = 1$ sur $[0,\frac{1}{ 4}]$ et $\supp
\varphi^0\subset
[0,1]$, et on pose $\psi^0 = \varphi^0(\frac{\cdot}{ 4}) - \varphi^0$. Si
$j\in \SZ$, on note $B_{j}$ la boule ferm\'ee  $ \overline
B (0,2^j)$ dans $\SR^3 $, $ \Gamma_{j}$ la couronne
$$\Gamma_{j} = \{\xi\in \SR^3\ ;\ 2^{j - 1}\le \vert \xi \vert\le
2^{j + 1}\},$$
et $S_{j}, \Delta_{j}$ les op\'erateurs $\varphi^0(-4^{-j}\Delta),
\psi^0(-4^{-j}\Delta)$ o\`u $\Delta$ est le Laplacien sur $\SR^3$. On
emploiera aussi les op\'erateurs
$\widetilde \Delta_{j} = \Delta_{j - 2} + \Delta_{j - 1} + \Delta_{j}
+
\Delta_{j + 1} + \Delta_{j + 2}$, qui v\'erifient l'identit\'e utile
\begin{equation}\label{eq 15}
\widetilde \Delta_{j}\ \Delta_{j} = \Delta_{j},
\end{equation}
et on \'ecrira $\widetilde \Gamma_{j}$ pour $\Gamma_{j - 2}\cup
\Gamma_{j - 1}\cup \Gamma_{j}\cup \Gamma_{j + 1}\cup \Gamma_{j +
2}$.\medskip

On a alors $\Delta_{j} = S_{j + 1} - S_{j}$ et $\displaystyle
\lim_{j\rightarrow +\infty} S_{j} = I$ dans $\Ss'$, d'o\`u il
r\'esulte que
\begin{equation}\label{eq 16}
f = \displaystyle \sum_{j\in \SZ} \Delta_{j} f
\end{equation}
pour toute distribution temp\'er\'ee $f$ telle que $\displaystyle
\lim_{j\rightarrow -\infty} S_{j} f = 0$.\medskip

En s'inspirant de Meyer (\cite{bib M}), on adopte la terminologie
suivante.

\begin{defin}\label{def 5}
On appelle espace fonctionnel (sous-entendu : adapt\'e \`a la
d\'ecomposition de Littlewood-Paley) tout espace de Banach $E$ tel
que
\begin{itemize}
\item[a)] $\Ss\subset E\subset\Ss'$, les injections \'etant continues;
\item[b)] deux cas sont possibles : dans le cas 1, $\Ss$ est dense dans
$E$, et dans le cas 2, $\Ss$ est dense dans un espace de Banach dont $E$
est le dual;
\item[c)] si $f\in E,\ \displaystyle \lim_{j\rightarrow -\infty} S_{j} f =
0$
dans $\Ss'$.
\end{itemize}

\end{defin}\medskip

 Si $f\in E$, on a
$\displaystyle \lim_{j\rightarrow +\infty} S_{j} f = f$  et de m\^eme $\displaystyle \lim_{t\rightarrow 0}
e^{t\Delta} f = f$ pour la topologie forte de $E$ dans le cas 1 et pour
la topologie faible $*$ dans le cas 2. Ceci induira la propri\'et\'e
analogue sur les solutions de (\ref{eq 1}) dans les espaces
favorables.\medskip

Remarquer qu'un espace fonctionnel n'est pas n\'ecessairement
caract\'eris\'e par la d\'ecomposition de Littlewood-Paley. Par
exemple, tous les espaces $L^p,\ 1\le p\le \infty$, sont des espaces
fonctionnels au sens de la d\'efinition pr\'ec\'edente.\medskip

Si $f$ et $g$ sont deux distributions temp\'er\'ees v\'erifiant
(\ref{eq 16}) et telles que le produit $fg$ ait un sens, on a formellement
\begin{equation}\label{eq 17}
fg = 
\sum_{\substack{\scriptstyle k,l\in \SZ\\ \scriptstyle
\vert k - l\vert \le 2}} \Delta_{k} f\ \Delta_{l} g
+ \displaystyle \sum_{\substack{\scriptstyle k,l\in \SZ\\ \scriptstyle
\vert k - l\vert \ge 3}} \Delta_{k} f\ \Delta_{l} g .
\end{equation}
Si $\vert k - l\vert\ge 3$, la transform\'ee de Fourier de
$\Delta_{k}
f \Delta_{l} g$ est support\'ee dans la couronne $\Gamma_{j - 1}\cup
\Gamma_{j}\cup \Gamma_{j + 1}$, o\`u $j = \max(k,l)$, et si
$\vert k - l\vert \le 2$, elle est support\'ee dans la boule $B_{j +
2}$.
La formule (\ref{eq 17}) isole ainsi les produits dont le spectre
contient 0 (chevauchement spectral) des autres dont le spectre ne
contient pas 0 (s\'eparation spectrale).

\begin{defin}\label{def 6}
Un espace fonctionnel $E$ est dit compatible (sous-en\-ten\-du : avec la
non-lin\'earit\'e de l'\'equation (\ref{eq 1})) lorsque :
\begin{itemize}
\item[a)] $E$ est invariant,
\item[b)] il existe une suite $\eta = (\eta_{n})_{n\in \SZ}$ telle que, si
$f,g\in E$ avec $\supp \widehat f\subset \Gamma_{k}, \supp
\widehat g\subset \Gamma_{l}, k,l\in \SZ$, et si $j\in \SZ$, alors
$\Delta_{j}(fg)\in E$, et
\begin{equation}\label{eq 18}
\Vert \Delta_{j} (fg)\Vert_{E}\le \eta_{\max(k - j, l - j)}
\ 2^{k + l - j}\ \Vert f \Vert_{E} \Vert g\Vert_{E}.
\end{equation}
\end{itemize}

\end{defin}\medskip

\noindent La propri\'et\'e b) et l'in\'egalit\'e (\ref{eq 18})
appellent
plusieurs remarques.

D'abord le fait que le produit $fg$ est bien d\'efini sous les
hypoth\`eses \'enonc\'ees, qui entra\^{\i}nent $f = S_{k + 1} f$ et $g =
S_{l + 1} g$, ce qui montre que $f$ et $g$ sont des fonctions de
classe $C^\infty$. Elles sont de plus born\'ees, en vertu de
l'injection de $E$ dans $\Ss'$ et de l'invariance par translation de
$E$ (voir la Proposition 13 pour un \'enonc\'e plus pr\'ecis). Par
cons\'equent, $\Delta_{j} (fg)$ est \'egalement bien d\'efini.\medskip

Ensuite, il convient de distinguer les deux m\^emes cas que dans la
formule (\ref{eq 17}). Si $\vert k - l\vert \ge 3$ (s\'eparation
spectrale), on a $\Delta_{j} (fg) = 0$ d\`es que $\vert\mbox{max}
(k,l) -
j\vert \ge 3$. Compte tenu de la forme particuli\`erement simple que
prend alors (\ref{eq 17}), on voit que (\ref{eq 18}) se ram\`ene \`a
\begin{equation}\label{eq 19}
\Vert fg\Vert_{E}\le C\ 2^{\min(k,l)}\ \Vert f \Vert_{E} \Vert
g\Vert_{E},
\end{equation}
pour une constante $C$ ind\'ependante de $f,g,k,l$.\medskip

Le cas sensible est celui o\`u $\vert k - l\vert\le 2$ (chevauchement
spectral) : (\ref{eq 18}) se simplifie en
\begin{equation}\label{eq 20}
\Vert \Delta_{j} (fg)\Vert_{E}\le C\ {\eta_{k - j}}\ 4^k\ 2^{-j}\ \Vert
f \Vert_{E}
\Vert g\Vert_{E},
\end{equation}
quitte \`a modifier la suite $\eta$, avec $C$ constante ne
d\'ependant
pas de $f,g,j,k$. Il n'y a en revanche aucune raison de supposer
$fg\in E$ (voir section 3.3 pour des exemples). Noter toutefois que
$\Delta_{j} (fg) = 0$ d\`es que $j\ge k + 5$ : on posera presque
toujours
$\eta_{n} = 0$ pour $n\le -5$.\medskip

La constante $C$ de (\ref{eq 20}), et surtout la suite $\eta$,
d\'ependent du choix de la fonction $\varphi^0$ d\'efinissant
les op\'erateurs $\Delta_{j}$. Mais leur existence ne d\'epend que de
l'espace $E$, et si $\widetilde\varphi^0$ est une autre fonction
d\'efinissant d'autres op\'erateurs $\widetilde\Delta_{j}$,
l'in\'egalit\'e correspondant \`a (\ref{eq 20}) reste vraie, avec
une
suite $\widetilde\eta = (\widetilde\eta_{n})_{n\in \SZ}$ qui
v\'erifie
$\widetilde \eta_{n}\le C (\eta_{n - 2} + \dots + \eta_{n + 2})$ pour
une certaine constante $C$.

Enfin, il faut souligner que la forme de l'in\'egalit\'e (\ref{eq
18}) est presque enti\`erement dict\'ee par l'invariance de $E$. Si,
en effet, on suppose seulement l'existence, pour tous $j,k,l\in \SZ$,
d'une constante $C(j,k,l)$ telle que
$$\Vert \Delta_{j} (fg)\Vert_{E}\le C(j,k,l) \Vert f \Vert_{E}
\Vert g\Vert_{E}$$
quand $\supp \widehat f\subset \Gamma_{k}$ et $\supp
\widehat g\subset \Gamma_{l}$, alors on d\'eduit de l'invariance de
$E$ l'existence de constantes $D(m,n),\ m,n\in \SZ$, v\'erifiant
l'\'egalit\'e
$$C(j,k,l) = D(k - j, l - j) 2^{k + l - j}.$$
L'hypoth\`ese suppl\'ementaire implicite dans (\ref{eq 18}) est donc
seulement que les constantes $D(m,n)$ ne d\'ependent que de
$\mbox{max} (m,n)$. Cette hypoth\`ese n'intervient d'ailleurs que
dans le cas ``facile'' o\`u $\vert k - l\vert\ge 3$, c'est-\`a-dire
dans
l'in\'egalit\'e (\ref{eq 19}).\bigskip

\subsection{Th\'eor\`emes d'existence, de r\'egularit\'e, et
contre-exemples}$\phantom{x}$\medskip

\begin{defin}\label{def 7}
On appelle bon espace tout espace fonctionnel compatible tel que
$\eta\in l^1(\SZ)$.
\end{defin}
\noindent Il r\'esulte de ce qui pr\'ec\`ede que la propri\'et\'e
``$\eta\in  l^1(\SZ)$'' ne d\'epend pas d'un choix particulier des
op\'erateurs
$\Delta_{j}$.

\begin{theo}\label{th 8}
Si $E$ est un bon espace, il existe un espace de Banach $F$ formant
avec lui un couple admissible.
\end{theo}

On en d\'eduit, en appliquant la Proposition 4, l'existence  d'une solution de
l'\'equation (\ref{eq 1}) dans
l'es\-pa\-ce $\Ff$ construit \`a partir de $F$, pour toute donn\'ee $u_{0}\in E$ assez
petite.

On peut obtenir une version plus pr\'ecise de ce th\'eor\`eme, en
renfor\c cant un peu l'hypoth\`ese.

\begin{theo}\label{th 9}
Soit $E$ un bon espace tel que $\displaystyle \sum_{n\ge 0} n
\eta_{n} <
+\infty$.\par
\noindent Alors il existe un espace de Banach $G\subset E$ formant
avec $E$ un couple admissible. De plus, si $u\in \Gg$ est une
solution
de (\ref{eq 1}) pour une donn\'ee initiale $u_{0}\in E$ quelconque,
alors\par
\noindent - $u\in \Cc ([0,\infty[;E)$ si $\Ss$ est dense dans
$E$ (cas 1, D\'efinition \ref{def 5}),\par
\noindent - $u\in \Cc (]0,\infty[;E)$ et $\displaystyle
\lim_{t\rightarrow 0} u(t) = u_{0}$ pour la topologie $\ast$ faible
si $\Ss$ est dense dans le pr\'edual de $E$ (cas 2).
\end{theo}

Ce r\'esultat est celui qui g\'en\'eralise le plus directement la
construction de Kato, dans la mesure o\`u il donne des solutions de
(\ref{eq 1}) r\'eguli\`eres. Il est apparent\'e \`a une conjecture
\'enonc\'ee par Meyer dans [M], qu'on peut r\'esumer en ``si
$E$ est un espace de Banach invariant qui s'injecte contin\^ument
dans
l'espace de Morrey $M^3_{2}$ (voir section 3.4), alors on peut r\'esoudre
(\ref{eq 1}) dans $\Cc ([0,\infty[;E)$ (cas 1) pour toute donn\'ee
$u_{0}\in E$ assez petite (modifier comme ci-dessus dans le cas 2)''.
Cette conjecture pose un probl\`eme de r\'egularit\'e, puisque
l'existence d'une solution de (\ref{eq 1}) pour toute donn\'ee assez
petite dans $E$ est garantie par l'inclusion de $E$ dans $M^3_{2}$,
qui est un bon espace (voir Proposition 23). Le th\'eor\`eme 9 donne
un
r\'esultat positif dans la direction de cette conjecture.\medskip

La question de l'optimalit\'e du th\'eor\`eme 8 motive le dernier
th\'eor\`eme de cette section. On a besoin de la notion suivante.

\begin{defin}\label{def 10}
Une suite positive $(\eta_{n})_{n\in \SZ}$ est r\'eguli\`ere s'il
existe une constante $C > 0$ telle que
$$\frac{1}{ C} \eta_{n + 1}\le \eta_{n}\le C\ \eta_{n + 1}$$
pour tout $n\in  \SZ$ v\'erifiant $\eta_{n} > 0$.
\end{defin}

\begin{theo}\label{th 11}
Soit $\eta = (\eta_{n})_{n\in \SZ}$ une suite d\'ecroissante
r\'eguli\`ere
telle que $$\displaystyle \sum_{n\ge 0} 2^{-3n}\
\displaystyle\inf_{k\ge n}\ 2^{3k}\,\eta^2_{k} = +\infty.$$
Alors, il existe un espace fonctionnel $E$, compatible et associ\'e \`a
$\eta$ via (\ref{eq 18}), tel que, \'etant donn\'e   un espace de 
Banach
$F$ v\'erifiant (P2), l'in\'egalit\'e (\ref{eq bicontinuite}) est en
d\'efaut. 
\end{theo}

En  d'autres termes, il n'y a pas d'espace $F$ formant avec $E$ un couple
admissible; cela ne signifie pas pour autant qu'il n'y ait pas de
solutions de (\ref{eq 1}) dans l'un de nos espaces $\Ff$.
 
 L'hypoth\`ese pr\'ec\'edente implique $\displaystyle \sum_{n\ge 0}\
\eta_{n}^2 =
+\infty$, et, r\'eciproquement, est v\'erifi\'ee d\`es que
$\eta\notin l^2(\SN)$ et que la suite $(2^{3n/2}\ \eta_{n})_{n\in \SN}$
est croissante \`a partir d'un certain rang.\medskip

On voit qu'il subsiste une lacune
\`a combler entre l'hypoth\`ese $l^1$ donnant le th\'eor\`eme 8 et
l'hypoth\`ese du type non-$l^2$ donnant les contre-exemples.\medskip

Il faut \'egalement prendre garde au fait que cette\ discussion de
l'opti-\par
\noindent malit\'e porte sur des classes d'espaces de Banach, et non
pas
sur les espaces consid\'er\'es s\'epar\'ement. En effet, il existe
aussi,
pour toute suite $\eta\notin l^2$,
un espace de Banach $\widetilde E$, compatible et associ\'e \`a
$\eta$,
pour lequel on peut trouver un espace $F$ formant avec $\widetilde E$
un couple admissible. Il suffit de prendre $\widetilde E = E\cap L^3
(\SR^3)$, o\`u $E$ est n'importe quel espace v\'erifiant (\ref{eq
20})
avec $\eta$ : la construction originelle de Kato convient \`a
$\widetilde E$.

Les contre-exemples du th\'eor\`eme 11 sont toutefois
suffisamment naturels pour inclure l'espace de Besov $\dot
B^{-1,\infty}_{\infty}$, comme on le verra dans la quatri\`eme partie
(voir \'egalement la Proposition 22).

\begin{cor}\label{co 12}
Il est impossible de faire fonctionner la m\'ethode KW \`a
partir de $E = \dot B^{-1,\infty}_{\infty} (\SR^3)$.
\end{cor}

\noindent L'importance de ce r\'esultat provient de la maximalit\'e
de $\dot B^{-1,\infty}_{\infty}$, observ\'ee par Meyer [M].

\begin{prop}\label{pro 13}
Tout espace de Banach invariant et qui s'injecte contin\^ument dans
$\Ss'$ est inclus dans $\dot B^{-1,\infty}_{\infty}$, avec injection
continue \'egale\-ment. Plus pr\'ecis\'ement, il existe une constante
$C\ge0$ telle que pour toute $f\in E$ et tout $j\in \SZ$,
$$
\|\Delta_j f\|_\infty \le \, C\, 2^j \, \|f\|_E.
$$
\end{prop}

\preuve  Soit $k (x) = e^{-\frac{1}{ 4}\vert
x\vert^2}$. Si
$E$ est un espace de Banach inclus dans $\Ss'$, il existe une constante
$C$ telle que 
$$
|(f,k)| \le C \|f\|_E
$$
pour tout $f \in E$. Puisque $E$ est invariant par translation
cela signifie
$$\Vert e^\Delta f\Vert_{L^\infty}\le C\Vert f\Vert_{E}.$$
 Comme $E$
v\'erifie (\ref{eq 3}), on en d\'eduit par changement d'\'echelle
$$\Vert e^{t\Delta} f\Vert_{L^\infty}\le C t^{-1/2}\ \Vert f
\Vert_{E}$$
pour tout $t > 0$.

Or, $\displaystyle \sup_{t > 0}\ \sqrt t\ \Vert e^{t \Delta}
f\Vert_{L^\infty}$ est une norme \'equivalente \`a la norme de $\dot
B^{-1,\infty}_{\infty}$ (voir par exemple Cannone \cite{bib C} pour
une d\'emonstration), d'o\`u la proposition.
\qed\medskip

La preuve des deux r\'esultats positifs (th\'eor\`emes \ref{th 8} et
\ref{th 9}) est donn\'ee \`a la suite de ces lignes : l'espace $F$
est
un espace de type Besov construit au-dessus de $E$. Le th\'eor\`eme
\ref{th 11} est de d\'emonstration plus d\'elicate : la quatri\`eme
partie lui est consacr\'ee, apr\`es que dans la troisi\`eme on montre
comment on peut retrouver l'ensemble des r\'esultats pr\'ealablement
connus.\bigskip

\subsection{Existence :\ preuve du th\'eor\`eme
8}$\phantom{x}$\medskip

Soit $E$ un bon espace, et $N$ un r\'eel $> 0$, pour le moment
quelconque. L'espace $F$ sera l'un des espaces not\'es
$C^{N,\infty}_{E}$ : par d\'efinition, $f\in C^{N,\infty}_{E}$
signifie que
$$f = \displaystyle \sum_{j\in \SZ} \Delta_{j} f$$
dans $\Ss'$, que $\Delta_{j} f\in E$ pour tout $j$, et que
$$\Vert f\Vert_{C^{N,\infty}_{E}} = \displaystyle \sup_{j\in \SZ}
(1 + 2^j)^N\ \Vert \Delta_{j} f\Vert_{E} < +\infty.$$
Cet espace est bien complet.\medskip

Puisque, d'apr\`es la Proposition 13, on a
$$\Vert \Delta_{j} f\Vert_{L^\infty}\le C\ 2^j\ \Vert \Delta_{j}
f\Vert_{E},$$
l'espace $C^{N,\infty}_{E}$ est inclus dans $L^\infty$ quand $N > 1$.
En particulier, le produit de deux \'el\'ements de $C^{N,\infty}_{E}$
est d\'efini dans ce cas. Il n'y a en revanche aucune raison pour que
$C^{N,\infty}_{E}$ soit inclus dans $E$, quel que soit $N$.\medskip

On choisit maintenant et pour toute la suite le param\`etre $N$ parmi
les entiers pairs $\ge 4$, et on note $F = C^{N,\infty}_{E}$. Il
s'agit de
d\'emontrer que le couple $(E,F)$ est admissible.\medskip

La propri\'et\'e (P1) est vraie par hypoth\`ese sur $E$.\medskip

La propri\'et\'e (P2) r\'esulte d'une s\'erie de remarques
simples.\medskip

Tout d'abord, l'injection de $F$ dans $L^\infty$ implique celle de
$F$
dans $\Ss'$ (contin\^ument). Pour v\'erifier que la classe $\Ss$
s'injecte dans $F$, on utilise l'existence d'une fonction $\psi^N\in
C^\infty (\SR^+),\ \supp \psi^N = [\frac{1}{ 4}, 4]$, telle que,
pour tout $j\in \SZ$
\begin{equation}\label{eq 21}
2^{jN} \Delta_{j} = \Delta_{j,N} (-\Delta)^{N/2},
\end{equation}
o\`u $\Delta_{j,N} = \psi^N (-4^{-j}\Delta)$. Si $f\in \Ss\subset E$, alors
$(-\Delta)^{N/2} f\in
\Ss$. En vertu de l'invariance par translation de $E$,  $L^1$ est un module
de convolution sur $E$, donc il existe une constante $C$ telle que
$$\Vert 2^{jN} \Delta_{j} f\Vert_{E}\le C$$
pour tout $j$. Ceci implique $f\in F$. On v\'erifie sans peine la
continuit\'e de l'injection.\medskip

Ensuite la norme $F$ est invariante par translation parce que $E$ est
invariant, et pour la m\^eme raison, $F$ est stable par homoth\'etie
de rapport $\lambda,\ \lambda > 0$. Si $t > 0,\ f\in F$, alors $\Vert
f\Vert_{t,F}$ est \'equivalente \`a
$$\sup_{j\in \SZ} \big(1 + 2^j \sqrt t\big)^N\, \Vert \Delta_{j}
f\Vert_{E},$$ et ce uniform\'ement par rapport \`a $t$. On identifie dans la
suite ces deux normes. Il devient \'evident que les $\Vert
.\Vert_{t,F}$, o\`u $t$ appartient \`a un compact de $]0,\infty[$,
sont uniform\'ement \'equivalentes entre elles.\medskip

Enfin, si $f\in E$ alors $e^\Delta f\in E$  ce qui implique
$$\sup_{j\in \SZ} \Vert \Delta_{j} e^\Delta f\Vert_{E} < +\infty.$$
On a de m\^eme $(-\Delta)^{N/2}\, e^\Delta f\in E$, donc d'apr\`es
(\ref{eq 21})
$$\sup_{j\in \SZ} 2^{jN}\, \Vert \Delta_{j} e^\Delta f\Vert_{E} <
+\infty.$$
Ceci montre que $e^\Delta$ est continu de $E$ dans $F$, et ach\`eve
de
prouver que (P2) est vraie.\medskip

On note la formule utile
\begin{equation}\label{eq 22}
\Delta_{j}\ e^{t\Delta} = (2^j \sqrt t)^{-N}\ (-t\Delta)^{N/2}\
e^{t\Delta}\ \Delta_{j,N}.
\end{equation}
\medskip

Il reste \`a v\'erifier (P3). Soient $u,v\in \Ff$, qu'on suppose pour
simplifier un peu de norme 1 tous les deux : $\Vert u\Vert_{\Ff}
= \Vert v\Vert_{\Ff} = 1$. On doit d'abord montrer l'existence d'une
constante $C$ telle que
\begin{equation}\label{eq 23}
\Vert \Delta_{j} B(u,v) (t)\Vert_{E}\le C \big(1 + 2^j\sqrt t\big)^{-N}
\end{equation}
pour tous $j\in \SZ, t > 0$, o\`u $B(u,v)(t)$ est la distribution
temp\'er\'ee d\'efinie par
$$B(u,v) (t) = \int^t_{0} e^{(t - \tau)\Delta}\ P(D) u(\tau) v(\tau)
d\tau.$$
Cette formule est \`a interpr\'eter comme la s\'erie
$$
B(u,v) (t) = \sum_{j \in \SZ}\int^t_{0} \Delta_j\{e^{(t -
\tau)\Delta}\ P(D) u(\tau) v(\tau)\} d\tau,
$$ dont la convergence dans $\Ss'$ est une
cons\'equence de (\ref{eq 23}). En effet,   la Proposition 13  entra\^{\i}ne
$$\Vert \Delta_{j} B(u,v) (t)\Vert_{L^\infty}\le C 2^j,$$
ce qui implique  que $B(u,v)(t) \in \dot B^{-1,\infty}_{\infty}$ et que,  $
\Delta_{j} B(u,v) (t)$ \'etant le terme g\'en\'eral, la
s\'erie converge pour la topologie faible  $\ast$ de $\dot
B^{-1,\infty}_{\infty}$.
\medskip

Le point de d\'epart est l'in\'egalit\'e \'el\'ementaire
$$\Vert \Delta_{j} B(u,v) (t)\Vert_{E}\le \int^t_{0} \Vert \Delta_{j}
e^{(t -
\tau) \Delta}\ P(D) u(\tau) v(\tau)\Vert_{E}\ d\tau.$$\medskip

\begin{lem}\label{le 14}
Il existe pour tout $p > 0$ une constante $C$, ne d\'ependant que de
$p$, telle que
$$\Vert \Delta_{j} e^{(t - \tau) \Delta}\ P(D) u(\tau)
v(\tau)\Vert_{E}\le C \big(1 + 2^j \sqrt{t - \tau}\big)^{-p}2^j\ \Vert
\Delta_{j} u(\tau) v(\tau)\Vert_{E}.$$
\end{lem}
\medskip

\preuve
On part de (\ref{eq 15}) pour \'ecrire $\widetilde \Delta_{j}^2
\Delta_{j} = \Delta_{j}.$
\noindent L'op\'erateur $P(D)$ \'etant ``\`a coefficients constants'',
on en d\'eduit
$$\Vert \Delta_{j} e^{(t - \tau) \Delta}\ P(D) u(\tau)
v(\tau)\Vert_{E}\le \Vert \widetilde \Delta_{j} e^{(t - \tau)
\Delta}\Vert\, \Vert \widetilde \Delta_{j} P(D)\Vert\,\Vert
\Delta_{j} u(\tau)  v(\tau)\Vert_{E}.$$
 L'action de $L^1$ sur $E$ par convolution implique
$\Vert
\widetilde \Delta_{j} e^{(t - \tau) \Delta}\Vert\le C$, et de
plus, si $2^j\sqrt{t - \tau} > 1$, la formule (\ref{eq 22})
s'applique (quitte \`a  choisir \`a la place de $N$ un autre entier
pair sup\'erieur \`a $p$) et donne
$$\Vert \widetilde \Delta_{j} e^{(t - \tau)\Delta}\Vert\le
C \big(2^j\sqrt{t - \tau}\big)^{-p}.$$
D'autre part, il existe $\psi^\sharp\in C^\infty(\SR^+)$ \`a support
dans
$]0,+\infty[$ telle que
$$\widetilde \Delta_{j} P(D) = 2^j \psi^\sharp (-4^{-j}\Delta).$$
Ceci implique
$$\Vert\widetilde \Delta_{j} P(D)\Vert\le C 2^j.$$
Le lemme en r\'esulte directement.
\qed
\medskip

On est ainsi ramen\'e, pour prouver (\ref{eq 23}), \`a estimer la
quantit\'e
$$B_{j}(t) = 2^j\int^t_{0} \big(1 + 2^j\sqrt{t - \tau}\big)^{-p}\ \Vert
\Delta_{j}( u(\tau) v(\tau))\Vert_{E}\ d\tau,$$
o\`u $p$ est un r\'eel $> 0$ qu'on choisira plus tard.

On utilise la formule (\ref{eq 17}), sous une forme modifi\'ee
obtenue
apr\`es regroupement de termes s'\'ecrivant
\begin{multline}\label{eq 24}
\Delta_{j}(fg)=\Delta_{j} (\Delta_{j} f\, S_{j - 2} g)
+\Delta_{j} (S_{j - 2} f\, \Delta_{j} g)\\
\\+\Delta_{j} \bigg( \sum_{k\ge j - 4} \Delta_{k} f
\, \widetilde \Delta_{k} g\bigg).
\end{multline}

\noindent Il suffit d'indiquer comment s'estiment les deux termes
suivants :
$$R_{j}(t) = 2^j \int^t_{0} \big(1 + 2^j\sqrt{t - \tau}\big)^{-p}\ \Vert
\Delta_{j} u(\tau) S_{j - 2} v(\tau)\Vert_{E}\ d\tau,$$
$$C_{j}(t) = 2^j \int^t_{0} \big(1 + 2^j\sqrt{t - \tau}\big)^{-p}\bigg\Vert
\Delta_{j} \bigg(\sum_{k\ge j - 4} \Delta_{k} u(\tau) \widetilde
\Delta_{k}  v(\tau)\bigg)\bigg\Vert_{E}\ d\tau.$$\medskip

Si $v(\tau)\in F$, alors $S_{j - 2} v(\tau) = \displaystyle
\sum_{j'\le j - 3}\ \Delta_{j'} v(\tau)$ dans $\Ss'$.
L'in\'egalit\'e (\ref{eq 19}) donne alors
\begin{align*}
\Vert \Delta_{j} u(\tau) S_{j - 2}
v(\tau)\Vert_{E}&\le \sum_{j'\le j -
3} \Vert \Delta_{j} u(\tau)\, \Delta_{j'} v(\tau)\Vert_{E}\\
&\le C\,  2^j\, \Vert \Delta_{j} u(\tau)\Vert\, \sup_{j'}
\Vert\Delta_{j'} v(\tau)\Vert_{E}\\
&\le C\, 2^j\, \big(1 + 2^j\sqrt \tau\big)^{-N}.
\end{align*}
On obtient pour le terme rectangle
\begin{align*}
R_{j}(t)&\le  C \int^t_{0} \big(1 + 2^j \sqrt{t -
\tau}\big)^{-p} \big(1 + 2^j
\sqrt \tau\big)^{-N} 4^j d\tau\\
&\le C \int^{4^jt}_{0} (1 + 4^j t - s)^{-p/2} (1 +
s)^{-N/2} ds.
\end{align*}
En choisissant $p\ge N$, il vient
\begin{equation}\label{eq 25}
R_{j}(t)\le C\ \min (1,4^j t) \big(1 + 2^j \sqrt t\big)^{-N}.
\end{equation}
\medskip

Pour le terme carr\'e, on utilise l'in\'egalit\'e (\ref{eq 20}), ce
qui donne
$$C_{j}(t)\le C\int^t_{0} \big(1 + 2^j \sqrt {t - \tau}\big)^{-p}\
\sum_{k\ge j - 4} \eta_{k - j} \big(1 + 2^k \sqrt \tau\big)^{-2N} 4^k
d\tau.$$\medskip
\noindent Si $2^j \sqrt t\ge 1$, il vient
$$C_{j}(t)\le C\sum_{k\ge j - 4} \eta_{k - j} \int^{4^kt}_{0} (1 +
4^j
t - 4^{j - k} s)^{-p/2} (1 + s)^{-N} ds.$$\medskip

\noindent En choisissant $p\ge 2N$ on obtient
$$C_{j}(t)\le C\sum_{k\ge j - 4}\eta_{k - j} (1 + 4^j t)^{-N},$$
c'est-\`a-dire
\begin{equation}\label{eq 26}
C_{j}(t)\le C \big(1 + 2^j \sqrt t\big)^{-2N}.
\end{equation}
\medskip

Lorsque $2^j \sqrt t < 1$, un calcul analogue donne
$$C_{j}(t)\le C\sum_{k\ge j - 4}\eta_{k - j} \ \min (1,4^k t),$$
soit
\begin{equation}\label{eq 27}
C_{j}(t)\le C \varepsilon (4^j t),
\end{equation}
o\`u $\varepsilon$ est la fonction d\'efinie sur $\SR^+$ par la
formule
$$\varepsilon (s) = \sum_{n\ge  - 4}\eta_{n} \ \mbox{min} (1,4^n s).$$
C'est une fonction positive, croissante, born\'ee, et surtout
v\'erifiant $\displaystyle \lim_{s\rightarrow 0} \varepsilon (s) =
0$.\medskip

Les in\'egalit\'es (\ref{eq 25} - \ref{eq 26} - \ref{eq 27}) donnent
finalement une version plus forte de (\ref{eq 23}), qui s'\'ecrit
\begin{equation}\label{eq 28}
\Vert \Delta_{j} B(u,v) (t)\Vert_{E}\le C\ \varepsilon (4^j t) (1 +
2^j
\sqrt t)^{-N}.
\end{equation}
\medskip

La fonction $\varepsilon$ permet de prouver que $B(u,v) (t)$ tend
vers 0 dans $\Ss'$ quand $t$ tend vers 0. En effet, si $f\in \Ss$, on
a
$$(B(u,v) (t), f) = \sum_{j\in \SZ} (\Delta_{j} B(u,v) (t),
\widetilde
\Delta_{j} f)$$
 La Proposition 13 et l'in\'egalit\'e
(\ref{eq 28}) entra\^{\i}nent
$$\Vert \Delta_{j} B(u,v) (t)\Vert_{L^\infty}\le C\ \varepsilon (2^j
t) 2^j,$$
ce qui implique
$$\vert (B(u,v) (t), f)\vert \le C\ \sum_{j\in \SZ}\ \varepsilon (2^j
t)\ 2^j\ \Vert \widetilde \Delta_{j} f\Vert_{L^1}.$$
L'espace des distributions temp\'er\'ees $f$ telles que
$$\sum_{j\in \SZ}\ 2^j \Vert\widetilde \Delta_{j} f\Vert_{L^1} <
+\infty$$ et $\displaystyle \lim_{j\rightarrow -\infty}
S_{j} f = 0$
est l'espace de Besov homog\`ene $\dot B^{1,1}_{1} (\SR^3)$ ; c'est le
pr\'edual de $\dot B^{-1,\infty}_{\infty}$ (voir Triebel \cite{bib Tr}). 
Par convergence domin\'ee, on obtient donc
$$\lim_{t\rightarrow 0} B(u,v) (t) = 0$$
pour la topologie faible  $\ast$  et de m\^eme, pour tout $t > 0$ fix\'e
$$\lim_{j\rightarrow -\infty} S_{j} B(u,v) (t) = 0$$
pour cette topologie, et par cons\'equent dans $\Ss'$. Ainsi $B(u,v)
(t)\in F$ pour tout $t$, et $\displaystyle \sup_{t > 0} \Vert B(u,v)
(t)\Vert_{t,F} < +\infty$.\medskip

Il reste \`a montrer que $B(u,v) (t)$ d\'epend contin\^ument de $t >
0$ pour la topologie de $F$.\medskip

Soient $t > 0$ et $h > 0,\ h\le \frac{t}{ 4}$. Si $\alpha\in ]0,\frac{t}{
4}]$, on a
\begin{multline*}
B(u,v) (t + h) - B(u,v) (t)
= \int^{\alpha + h}_{0} e^{(t + h - \tau)\Delta}\
P(D) u(\tau) v(\tau) d\tau\\
-\int^{\alpha}_{0} e^{(t - \tau)\Delta}\ P(D)
u(\tau) v(\tau) d\tau\\
+\int^t_{\alpha} e^{(t - \tau)\Delta}\ P(D) [u(\tau
+ h)
v(\tau + h) - u(\tau) v(\tau)] d\tau.
\end{multline*}
\medskip

En reprenant les calculs qui ont conduit \`a (\ref{eq 28}), on prouve
\begin{multline}\label{eq 29}
\Vert B(u,v) (t + h) - B(u,v) (t)\Vert_{t,F}\\
\le C\ \sup_{j\in \SZ} \varepsilon (4^j\alpha +
4^j h) \big(1 + 2^j\sqrt t\big)^{-1}
+C\ \sup_{j\in \SZ}\ \varepsilon (4^j\alpha) \big(1 +
2^j
\sqrt t\big)^{-1}\\
+C\  \sup_{\tau\in [\alpha,t]}\big\{ \Vert u(\tau + h) -
u(\tau)\Vert_{\tau,F} + \Vert v(\tau + h) - v(\tau)\Vert_{\tau,F}
\big\}.
\end{multline}
\noindent Par d\'efinition de $\Ff$, on a pour tout $\alpha > 0$
fix\'e
$$\lim_{h\rightarrow 0}\ \sup_{\tau\in [\alpha,t]}\big\{ \Vert u(\tau + h) -
u(\tau)\Vert_{\tau,F} + \Vert v(\tau + h) - v(\tau)\Vert_{\tau,F}
\big\} =
0.$$
On a aussi pour tout $t > 0$
$$\lim_{\alpha\rightarrow 0}\ \sup_{j\in \SZ}\ \varepsilon
(4^j\alpha)
\big(1 + 2^j\sqrt t\big)^{-1} = 0,$$
d'o\`u il r\'esulte que $\Vert B(u,v) (t + h) - B(u,v)
(t)\Vert_{t,F}$ tend vers 0 avec $h > 0$.\medskip

Le cas $h < 0$, tout \`a fait analogue, est laiss\'e au lecteur. Le
th\'eor\`eme 8 est ainsi compl\`etement d\'emontr\'e.\bigskip

\subsection{R\'egularit\'e :\ preuve du th\'eor\`eme
9}$\phantom{x}$\medskip

L'espace associ\'e \`a $E$ est, lorsque $\displaystyle \sum_{n\ge 0}
n \eta_{n} < +\infty$, diff\'erent de l'espace $F$ pr\'ec\'edent ; on
le note $G$. C'est l'espace de Besov inhomog\`ene $B^{N,\infty}_{E}$
construit au-dessus de $E$, d\'efini par les conditions $S_{0} f\in
E$ et $\Delta_{j} f\in E$ pour tout $j\ge 0$, avec
$$\Vert f\Vert_{G} = \Vert S_{0} f\Vert_{E} + \sup_{j\ge 0}\ 2^{jN}
\Vert\Delta_{j} f\Vert_{E} < +\infty.$$
L'entier $N$ est choisi comme \`a la section pr\'ec\'edente.\medskip

Si $t > 0$, la norme $\Vert .\Vert_{t,G}$ est \'equivalente \`a
$$\Vert S_{j(t)} f\Vert_{E} + \sup_{j\ge j(t)} \big(2^j \sqrt t\big)^N\
\Vert \Delta_{j} f\Vert_{E},$$
o\`u $j(t)$ est d\'efini par les in\'egalit\'es $2^{-j(t)}\le \sqrt
t < 2^{-j(t) + 1}$, et ce uniform\'ement par rapport \`a $t$. On
identifie donc les deux normes. Enfin, on note $\Gg$ au lieu de $\Ff$
l'espace des $u$ telles que
$$\Vert u\Vert_{\Gg} = \sup_{t > 0} \Vert u(t)\Vert_{t,G} < +\infty$$
qui d\'ependent contin\^ument de $t$.\medskip

La preuve du th\'eor\`eme 9 est parall\`ele \`a celle du th\'eor\`eme
8,
et presque compl\`etement laiss\'ee au lecteur. Deux points
m\'eritent
d'\^etre d\'etaill\'es, qui expliquent l'introduction de
l'hypoth\`ese $\displaystyle \sum_{n\ge 0}\ n \eta_{n} <
+\infty$.\medskip

Le premier est la d\'emonstration de l'in\'egalit\'e
$$\Vert S_{j(t)} B(u,v)(t)\Vert_{E}\le C\ \Vert u\Vert_{\Gg}\ \Vert
v\Vert_{\Gg}.$$
On part de
$$\Vert S_{j(t)} B(u,v)(t)\Vert_{E}\le \sum_{j < j(t)}\ \Vert
\Delta_{j} B(u,v) (t)\Vert_{E},$$
et, d\'emontrant que (\ref{eq 28}) est encore vrai, on obtient si
$\|u\|_{\Gg}= \|v\|_{\Gg}=1$, 
\begin{align*}
\Vert S_{j(t)} B(u,v)(t)\Vert_{E}&\le C\ \sum_{j < j(t)}\
\sum_{n\ge -4}\ \eta_{n} \min (1,4^{n + j}t)\\
&\le C\ \sum_{n\ge -4} (5 + n)\ \eta_{n}.
\end{align*}\medskip

Le deuxi\`eme est la d\'emonstration de
$$\lim_{h\rightarrow 0}\ \Vert S_{j(t)} \{B(u,v) (t + h) - B(u,v)
(t)\}\Vert_{E} = 0.$$
L'argument est semblable \`a celui utilis\'e pour (\ref{eq 29}), \`a
condition de remplacer $\varepsilon (4^j \alpha)$ par
$$\widetilde\varepsilon (t,4^{j}\alpha) = \sum_{j < j(t)}\ \sum_{n\ge
-4}\ \eta_{n} \min (1,4^{n + j}\alpha),$$
et de m\^eme pour $\varepsilon (4^j\alpha + 4^j h)$. Comme on a bien,
pour tout $t > 0$
$$\lim_{\alpha\rightarrow 0}\ \widetilde \varepsilon (t, 4^j\alpha) =
0,$$
le raisonnement est inchang\'e.\bigskip

\subsection{Autres r\'esultats de r\'egularit\'e}$\phantom{x}$\medskip

Le th\'eor\`eme 9 est en fait valable sous une condition
l\'eg\`erement plus faible, qui sera utilis\'ee dans la prochaine
partie.

\begin{prop}\label{pro 15}
Soit $E$ un espace fonctionnel invariant tel que\par
\noindent a) si $f,g\in E$ avec $\supp\widehat f\subset \Gamma_{j}$ et
$\supp\widehat g\subset B_{j - 2}, j\in \SZ$, alors $fg\in E$ et
$$\Vert fg\Vert_{E}\le C\ 2^j\ \Vert f\Vert_{E}\ \Vert g \Vert_{E},$$
o\`u $C$ est une constante ne d\'ependant que de $E$,\par
\noindent b) il existe une suite $(\eta_{n})_{n\in \SZ}$ telle que
$\displaystyle \sum_{n\ge 0}\ n \eta_{n} < +\infty$ et, si $f,g\in
E$ avec $\supp\widehat f\subset \Gamma_{k}$, $\supp\widehat g\subset
\widetilde \Gamma_{k}$, alors $\Delta_{j}(fg)\in E$ pour tout $j\in
\SZ$, et
$$\Vert \Delta_{j} (fg)\Vert_{E}\le \eta_{k - j}\ 4^k\ 2^{-j}\ \Vert
f\Vert_{E}\ \Vert g \Vert_{E}.$$
Alors, il existe un espace de Banach $G\subset E$ tel que le couple
$(E,G)$ soit admissible. De plus, si $u\in \Gg$ est une solution de
(\ref{eq 1}) pour une donn\'ee initiale $u_{0}\in E$ quelconque,
alors\par
\noindent - $u\in \Cc\ ([0,\infty[ ; E)$ si $\Ss$ est dense dans
$E$,\par
\noindent - $u\in \Cc\ (]0,\infty[ ; E)$ et $\displaystyle
\lim_{t\rightarrow 0} u(t) = u_{0}$ pour la topologie  faible $\ast$
si
$\Ss$ est dense dans le pr\'edual de $E$.
\end{prop}\medskip

\preuve On consid\`ere le m\^eme espace $\Gg$ que dans le th\'eor\`eme 9.
Le seule diff\'erence avec les preuves pr\'ec\'edentes est dans le
traitement des produits du type $\Delta_{j} u(\tau)\, S_{j -
2}v(\tau)$, o\`u $u,v\in \Gg$. L'hy\-po\-th\`e\-se a) donne
$$\Vert \Delta_{j} u(\tau)\, S_{j - 2} v(\tau)\Vert_{E}\le C\ 2^j\
\Vert\Delta_{j} u(\tau)\Vert_{E}\ \Vert S_{j - 2} v(\tau)\Vert_{E}.$$
Si $j - 2\le j(\tau)$, alors $\Vert S_{j - 2} v(\tau)\Vert_{E}\le C\
\Vert v(\tau)\Vert_{\tau,G}$, tandis que si $j - 2 > j(\tau)$, alors
\begin{align*}
\Vert S_{j - 2} v(\tau)\Vert_{E}&\le \Vert S_{j(\tau)}
v(\tau)\Vert_{E} +
\sum^{j - 1}_{k = j(\tau) + 1} \Vert \Delta_{k} v(\tau)\Vert_{E}\\
&\le C\ \Vert v(\tau)\Vert_{\tau,G}.
\end{align*}
On en d\'eduit
$$\Vert \Delta_{j} u(\tau)\ S_{j - 2} v(\tau)\Vert_{E}\le C\ 2^j\ \big(1
+
2^j\sqrt \tau\big)^{-N}\ \Vert u\Vert_{\Gg} \Vert v\Vert_{\Gg},$$
puis l'in\'egalit\'e (\ref{eq 25}). On termine la preuve comme pour
le
th\'eor\`eme 9.
\qed
\medskip

Enfin, la r\'egularit\'e de la solution peut encore \^etre
am\'elior\'ee
si l'hypo\-th\`ese a) de la proposition pr\'ec\'edente est
renforc\'ee.\medskip

\begin{prop}\label{pro 15 bis}
Soit $E$ un espace v\'erifiant toutes les hypoth\`eses
pr\'e\-c\'e\-dentes, \`a l'exception de a), remplac\'ee par :\par
\noindent a') si $f,g\in E$ avec $\supp\widehat f\subset \Gamma_{j}$
et
$\supp \widehat g\subset B_{j - 2}, j\in \SZ$, alors $fg\in E$ et
$$\Vert fg\Vert_{E}\le C\ \Vert f\Vert_{E}\ \Vert g\Vert_{L^\infty},$$
o\`u $C$ est une constante ne d\'ependant que de $E$. Alors, toute
solution $u\in \Gg$ de l'\'equation (\ref{eq 1}) est de classe
$C^\infty$ sur $]0,\infty[\times \SR^3$. De plus, pour tout entier
$n$
et tout multi-indice $\alpha$, il existe une constante $C_{n}$ et une
constante $C_{\alpha}$ telles que
\begin{equation}\label{eq 29-1}
\bigg\vert \frac{\partial^n}{ \partial t^n}\ u (t,x)\bigg\vert\le
\frac{C_{n}}{t^{n + 1/2}}
\end{equation}
et
\begin{equation}\label{eq 29-2}
\vert D^\alpha_{x}\ u(t,x)\vert\le \frac{C_{\alpha}}{ t^{(\vert
\alpha\vert + 1)/2}},
\end{equation}
pour tout $t > 0, x\in \SR^3$.
\end{prop}
\medskip

\preuve
La cl\'e est dans la modification suivante de (\ref{eq 25}) :
\begin{equation}\label{eq 29-3}
R_{j}(t)\le C\ \min (1,4^j t)\ \big(1 + 2^j \sqrt t\big)^{-N-1},
\end{equation}
valable lorsque $u,v\in \Gg$. Pour le voir, on \'ecrit si $j - 2 >
j(\tau)$
$$S_{j - 2} v(\tau) = \sum_{j(\tau)\le j'\le j - 3}\ \Delta_{j'}
v(\tau) + S_{j(\tau)}
v(\tau)$$
et on utilise la Proposition 13. On obtient, puisque $N>1$,
\begin{align*}
\Vert S_{j - 2} v(\tau)\Vert_{L^\infty}&\le  C\ 
\sum_{j(\tau)\le j'\le
j - 3} 2^{j'}\, \big(2^{j'}\sqrt \tau\big)^{-N}\Vert v\Vert_{\Gg}
+ C\ 2^{j(\tau)}\Vert v\Vert_{\Gg}\\
&\le \frac{C}{ \sqrt \tau}\ \Vert v \Vert_{\Gg}.
\end{align*}
On en d\'eduit que
$$\Vert S_{j - 2} v(\tau)\Vert_{L^\infty}\le C\ \frac{2^j}{ 1 +
2^j\sqrt \tau}\
\Vert v\Vert_{\Gg}$$
pour tout $j\in \SZ$. L'in\'egalit\'e (\ref{eq 29-3}) en
d\'ecoule en utilisant a') pour $f=\Delta_ju(\tau)$ et $g=S_{j-2}v(\tau)$.

Si on note momentan\'ement $\Gg_{N}$ au lieu de $\Gg$, afin de rendre
apparent le param\`etre $N$, les in\'egalit\'es (\ref{eq 26}) et
(\ref{eq 29-3}) montrent que $B$ envoie $\Gg_{N}\times \Gg_{N}$ dans
$\Gg_{N + 1}$ :
$$\Vert B(u,v)\Vert_{\Gg_{N + 1}}\le C_N \Vert u\Vert_{\Gg_{N}} \Vert
v\Vert_{\Gg_{N}}.$$
Comme $Su_{0}\in \Gg_{N}$ pour tout $N$ quand $u_{0}\in E$, un
argument de ``bootstrap'' classique donne $u\in \Gg_{N}$ pour tout $N$
\'egalement. On en d\'eduit (\ref{eq 29-2}) par des calculs
semblables
aux pr\'ec\'edents dont le d\'etail est laiss\'e au lecteur. Cela
implique ensuite (\ref{eq 29-1}).
\qed
\bigskip

\subsection{Existence et r\'egularit\'e locales}$\phantom{x}$\medskip

Pour r\'esoudre localement l'\'equation (1) sans condition de taille
sur la donn\'ee initiale, on pose a priori $u = Su_{0} + v$, ce qui
donne l'\'equation suivante d'inconnue $v$ :
\begin{equation}\label{eq a}
v = B(Su_{0}, Su_{0}) + 2B (Su_{0},v) + B (v,v).
\end{equation}

Si $T > 0$ est, pour le moment, quelconque, on d\'esigne par
$\Ff_{T}$ l'espace des fonctions continues $v$ de $]0,T[$ \`a valeurs
dans $F$, telles que
$$\Vert v\Vert_{\Ff_{T}} = \displaystyle \sup_{t > 0}\ \Vert
v(t)\Vert_{t,F} < +\infty.$$
L'espace $F$ est l'espace $C^{N,\infty}_{E}$ utilis\'e pour prouver le
th\'eor\`eme 8. On d\'efinit de la m\^eme mani\`ere l'espace
$\Gg_{T}$.\medskip

Soit $L$ l'op\'erateur lin\'eaire d\'efini par
$$L v = 2 B(Su_{0},v).$$
sur $\Ff_{T}$ ou $\Gg_{T}$. Bien qu'on ait ainsi d\'efini plusieurs
op\'erateurs, on les d\'esignera de la m\^eme fa\c con dans la
suite.\medskip

Il r\'esulte des preuves des th\'eor\`emes 8 et 9 que, sous leurs
hypoth\`eses, $L$ est continu sur $\Ff_{T}$ et sur $\Gg_{T}$,
uniform\'ement en $T$, et pour tout $u_{0}\in E$. Les r\'esultats
locaux seront une cons\'equence de l'observation suivante.

\begin{lem}
Si $\displaystyle \lim_{j\rightarrow +\infty}\ \Vert \Delta_{j}
u_{0}\Vert_{E} = 0$, alors les normes de $L$ sur $\Ff_{T}$ (si
$\displaystyle \sum_{n\ge 0}\ \eta_{n} < +\infty$) et sur $\Gg_{T}$ (si
$\displaystyle \sum_{n\ge 0}\ n\eta_{n} < +\infty$)  tendent vers 0
avec $T$.
\end{lem}
\medskip

Admettons-le un moment. Si $T$ est assez petit, l'op\'erateur $I - L$
est inversible sur $\Ff_{T}$ ou $\Gg_{T}$, et l'\'equation (\ref{eq
a}) est \'equivalente \`a
\begin{equation}\label{eq b}
v = (I - L)^{-1}\ B(Su_{0}, Su_{0}) + (I - L)^{-1} B (v,v).
\end{equation}
Le th\'eor\`eme 3 peut donc s'appliquer. Comme $(I - L)^{-1}\ B(Su_{0},
Su_{0}) = \frac{1}{2} (I - L)^{-1}L(
Su_{0})$, le lemme 17 implique aussi
$$\displaystyle \lim_{T\rightarrow 0}\ \Vert (I - L)^{-1}\ B(Su_{0},
Su_{0})\Vert_{\Ff_{T}} = 0$$
si $\displaystyle \sum_{n\ge 0}\ \eta_{n} < +\infty$, et de m\^eme
dans $\Gg_{T}$ si $\displaystyle \sum_{n\ge 0}\ n\eta_{n} < +\infty$.

On peut donc r\'esoudre (\ref{eq b}).\medskip

Ecrivons les r\'esultats obtenus.\medskip

\begin{theo}
Soit $E$ un bon espace et $u_{0}\in E$ tel que
$$\displaystyle
\lim_{j\rightarrow +\infty} \Vert \Delta_{j} u_{0}\Vert_{E} = 0.$$
Alors il existe $T > 0$ tel que l'\'equation (1) admette une solution
$u$ dans $\Ff_{T}$. De plus, si $\displaystyle \sum_{n\ge 0}\ n\eta_{n} <
+\infty$, alors $u \in \Gg_{T}$ et donc

- $u\in \Cc ([0,T[ ; E)$ si $S$ est dense dans $E$,

- $u\in \Cc (]0,T[ ; E)$ et $\displaystyle \lim_{t\rightarrow 0} u(t) =
u_{0}$ pour la topologie faible  $\ast$ si $S$ est dense dans le
pr\'edual de $E$.
\end{theo}\medskip

Le reste de ce paragraphe est consacr\'e \`a la preuve du lemme
17.\medskip

On se place d'abord dans $\Ff_{T}$ : soit $v = (v(t))_{0<t<T}$ tel que
$$\Vert \Delta_{j} v(t)\Vert_{E}\le \big(1 + 2^j\sqrt t\big)^{-N}$$
pour tous $j\in \SZ$ et $t\in ]0,T[$. Reprenant la preuve du th\'eor\`eme
8, on se ram\`ene \`a estimer les trois termes suivants :
$$R'_{j}(t) = 2^j \int^t_{0} \big(1 + 2^j\sqrt {t -
\tau}\big)^{-p}\ \Vert \Delta_{j} Su_{0}(\tau)\, S_{j - 2}
v(\tau)\Vert_{E}\ d\tau,$$
$$R''_{j}(t) = 2^j  \int^t_{0} \big(1 + 2^j\sqrt {t -
\tau}\big)^{-p}\ \Vert S_{j - 2} Su_{0}(\tau)\,  \Delta_{j}
v(\tau)\Vert_{E}\ d\tau,
$$
$$C_{j}(t) = 2^j  \int^t_{0} \big(1 + 2^j\sqrt {t -
\tau}\big)^{-p}\ \bigg\Vert \Delta_{j} \bigg( \sum_{k\ge j - 4}\
\Delta_{k} Su_{0}(\tau)\, \widetilde\Delta_{k}
v(\tau)\bigg)\bigg\Vert_{E}\ d\tau.$$
Chacun de ces termes est estim\'e comme dans les in\'egalit\'es
(25-26-27), en suivant soigneusement la d\'ependance par rapport \`a
$u_{0}$. On utilise toujours l'in\'egalit\'e
$$\Vert \Delta_{j} Su_{0} (\tau)\Vert_{E}\le C\ \big(1 + 2^j \sqrt
\tau\big)^{-N}\ \Vert \Delta_{j} u_{0}\Vert_{E},$$
qui se d\'emontre comme le lemme 14.\medskip

On obtient ainsi pour le premier terme rectangle :
$$R'_{j}(t)\le C\ \min (1, 4^j t)\ \Vert \Delta_{j} u_{0}\Vert_{E} \big(1 +
2^j \sqrt t\big)^{-N},$$
de la m\^eme fa\c con qu'on a obtenu (25). Pour le second, on \'ecrit
\begin{align*}
\Vert S_{j - 2} Su_{0}(\tau)\ \Delta_{j} v(\tau)\Vert_{E} &
\le C\
 \sum_{j'\le j - 3} 2^{j'}\ \Vert \Delta_{j'}
Su_{0}(\tau)\Vert_{E}\ \Vert \Delta_{j} v(\tau)\Vert_{E}
\\
&
\le C\ 2^j \big(1 + 2^j \sqrt \tau\big)^{-N} \sum_{j'\le j - 3}
2^{j' - j}\ \Vert \Delta_{j'} u_{0}\Vert_{E}.
\end{align*}
On en d\'eduit la majoration
$$R''_{j}(t)\le C\ \min (1, 4^j t)\  \sum_{j'\le j - 3} 2^{j' -
j}\ \Vert \Delta_{j'} u_{0}\Vert_{E} \big(1 + 2^j \sqrt
t\big)^{-N}.$$\medskip

Enfin, on a pour le terme carr\'e
$$C_{j}(t)\le C  \int^t_{0} \big(1 + 2^j\sqrt {t -
\tau}\big)^{-p}\  \sum_{k \ge j - 4}\ \eta_{k - j}\ \Vert \Delta_{k}
u_{0}\Vert_{E} \big(1 + 2^k\sqrt \tau\big)^{-2N}\ 4^k d\tau.$$
Comme pour (26-27), on obtient
$$C_{j}(t)\le C\ \displaystyle \sum_{k\ge j - 4}\ \eta_{k - j}\ \Vert
\Delta_{k} u_{0}\Vert_{E}\ \min (1,4^k t)\ \big(1 + 2^j \sqrt
t\big)^{-2N}.$$\medskip

L'ensemble de ces estimations donne
\begin{multline}\label{eq c}
\big(1
+ 2^j\sqrt t \big)^{N} \ \Vert\Delta_{j}
B(Su_{0},v)(t)\Vert_{E}
\le\\ C\{\min
(1,4^jt) \sum_{j'\le j} 2^{j' -
j}\Vert\Delta_{j'}u_{0}\Vert_{E} + \sum _{k\ge j - 4}
\eta_{k - j}\Vert \Delta_{k}u_{0}\Vert_{E} \min (1,4^kt)\}.
\end{multline}
Il est \'el\'ementaire de v\'erifier que le terme entre accolades tend
vers 0 avec $t$, uniform\'ement par rapport \`a $j$, si $\displaystyle
\lim_{l\rightarrow +\infty} \Vert \Delta_{l}u_{0}\Vert_{E} = 0$.

Cela prouve que la norme de $L$ sur $\Ff_{T}$ tend vers 0 avec
$T$.\medskip

Lorsque $\displaystyle \sum_{n\ge 0} n\eta_{n} < +\infty$, l'extension
\`a $\Gg_{T}$ repose encore une fois sur l'in\'egalit\'e
$$\Vert S_{j(t)} B(Su_{0},v) (t)\Vert_{E}\le \displaystyle
\sum_{j<j(t)}\ \Vert\Delta_{j} B(Su_{0},v)(t)\Vert_{E}.$$
On injecte l'in\'egalit\'e (\ref{eq c}), et on conclut : les
d\'etails sont laiss\'es au lecteur.



\section{Liens avec les r\'esultats ant\'erieurs}\bigskip

\subsection{L'espace $L^3 (\SR^3)$}$\phantom{x}$\medskip

C'est sur cet espace que Kato a construit
des solutions de (\ref{eq 1}) suivant la m\'ethode de la premi\`ere
partie.\medskip

\begin{prop}\label{por 16}
$L^3 (\SR^3)$ est un bon espace.
\end{prop}

\preuve
$L^3 (\SR^3)$ est un espace fonctionnel invariant. L'in\'egalit\'e
(\ref{eq 19}) r\'esulte de la Proposition 13 et de ce que toute
fonction born\'ee est un multiplicateur de $L^3$.\medskip

Soient $f,g\in L^3$, $\supp \widehat f\subset
\Gamma_{k}$, $\supp \widehat g\subset \widetilde
\Gamma_{k}$, et $j\le k + 4$.  Puisque $fg\in
L^{3/2}$, l'in\'egalit\'e d'Young donne
$$\Vert\Delta_{j} (fg)\Vert_{L^3}\le C\ 2^j\ \Vert
fg\Vert_{L^{3/2}}\le C\ 2^j\ \Vert f\Vert_{L^3}\ \Vert g\Vert_{L^3}.$$
Ceci montre (\ref{eq 20}) avec $\eta_{n} = 4^{-n}$, et conclut la
preuve.
\qed\medskip

Kato a employ\'e dans \cite{bib K} un espace $\Ff$ diff\'erent de celui
utilis\'e dans la preuve du th\'eor\`eme 8 (ou du th\'eor\`eme 9, qui
s'applique) : les r\'esultats d'unicit\'e de
Furioli-Lemari\'e-Terraneo montrent que ce sont les m\^emes solutions
qui sont obtenues ([F,LR,Te]).\bigskip

\subsection{Espaces de Lorentz}$\phantom{x}$\medskip

Meyer a montr\'e comment faire fonctionner la m\'ethode KW dans
$L^{(3,\infty)}$ (\cite{bib M}). On a en fait la

\begin{prop}\label{pro 17}
Si $1\le q\le \infty$, l'espace $L^{(3,q)}$ est un bon espace.
\end{prop}

\preuve
On rappelle que, si $1\le p,q\le \infty$, une fonction mesurable $f$
appartient \`a $L^{(p,q)}$ lorsque
$$\Vert f\Vert_{L^{(p,q)}} = \biggl(\frac{q}{ p}\
\int^{+\infty}_{0}\ [t^{1/p} f^\ast(t)]^q \frac{dt}{t}\biggr)^{1/q}
< +\infty,$$ et si $1\le p\le \infty$ et $q = \infty$, lorsque
$$\Vert f\Vert_{L^{(p,\infty)}} = \sup_{t > 0} t^{1/p} f^\ast(t) <
+\infty.$$
On a not\'e $f^\ast$ le r\'earrangement d\'ecroissant de $f$ :
$$f^\ast(t) = \inf\{s\ge 0 ; \vert\{\vert f \vert > s\}\vert\le
t\},\ \ t\ge 0.$$
Si $p > 1$, il existe une norme sur $L^{(p,q)}$, \'equivalente \`a
$\Vert \cdot\Vert_{L^{(p,q)}}$ (qui n'est pas une norme), et qui
rend
$L^{(p,q)}$ complet (Stein-Weiss [S,W], chapitre V). Enfin,
$L^{(p,p)}$ n'est autre que $L^p$.

Parmi tous les espaces de Lorentz, ce sont les $L^{(3,q)}$ qui sont
invariants.\eject

On a bien s\^ur $\Ss\subset L^{(3,q)}$, et l'inclusion $L^{3,q}\subset
L^{p_{1}} + L^{p_{2}}$, o\`u $1 < p_{1} < 3 < p_{2}$, montre que
$L^{(3,q)}\subset \Ss'$ et que $\displaystyle \lim_{j\rightarrow
-\infty}\ S_{j} f = 0$ si $f\in L^{(3,q)}$. Les $L^{(3,q)}$ sont donc
des espaces fonctionnels invariants. Pour montrer qu'ils sont
compatibles, on note d'abord que les fonctions born\'ees sont des
multiplicateurs de ces espaces, ce qui donne (19) (avec la
Proposition 13). D'autre part, on a
\begin{equation}\label{eq 30}
\Vert fg\Vert_{L^{(\frac{3}{ 2},q)}}\le C\ \Vert
f\Vert_{L^{(3,q)}}\
\Vert g\Vert_{L^{(3,q)}}
\end{equation}
et
 \begin{equation}\label{eq 31}
\Vert h\ast F\Vert_{L^{(3,q)}}\le C\ \Vert h\Vert_{L^{\frac{3}{
2}}}\ \Vert F\Vert_{L^{(\frac{3}{ 2},q)}}.
\end{equation}
Si on admet ces deux in\'egalit\'es, alors on montre que (\ref{eq 20})
est v\'erifi\'ee avec $\eta_{n} = 4^{-n}$, exactement comme dans le
cas de l'espace $L^3$.\medskip

Il suffit de d\'emontrer (\ref{eq 30}) dans le cas o\`u $f = g$.
L'identit\'e $(f^2)^\ast (t) = f^\ast(t)^2$ implique 
l'in\'egalit\'e
$$\Vert f^2\Vert_{L^{(\frac{3}{ 2},q)}}\le C\ \Vert
f\Vert_{L^{(3,\infty)}}\
\Vert f\Vert_{L^{(3,q)}},$$
et on conclut avec l'inclusion $L^{(3,q)}\subset L^{(3,\infty)}$.

Pour obtenir (\ref{eq 31}), on part de
$$\Vert h\ast F\Vert_{L^r}\le \Vert h\Vert_{L^{\frac{3}{ 2}}}\
\Vert f\Vert_{L^p}$$
chaque fois que $1\le p\le 3$ et $\frac{1}{ r} = \frac{1}{ p} -
\frac{1}{3}$. Le th\'eor\`eme d'interpolation de Hunt ([S,W],
p.197) fournit
$$\Vert h\ast F\Vert_{L^{(r,q)}}\le C\ \Vert h\Vert_{L^{\frac{3}{
2}}}\ \Vert f\Vert_{L^{(p,q)}}$$
pour tous $p\in ]1,3[$ et $q\in [1,\infty]$. On en d\'eduit (\ref{eq
31}).
\qed\bigskip

\subsection{Espaces de Besov et de Triebel-Lizorkin}$\phantom{x}$\medskip

Les espaces $\dot B^{\frac{3}{ p} - 1,\infty}_{p}$ ont \'et\'e
\'etudi\'es par Cannone dans sa th\`ese, notamment pour $3 < p \le 6$,
et lui ont permis de construire des solutions auto-similaires des
\'equations de Navier-Stokes ([C]).

\begin{prop}\label{pro 18}
Les espaces $\dot B^{\frac{3}{ p} - 1,q}_{p},\ 1\le q\le
\infty$ et $1\le p < \infty$, sont de bons espaces. Il en est de
m\^eme des espaces $\dot F^{\frac{3}{ p} - 1,q}_{p},\ 1\le p,q <
\infty$.
\end{prop}\medskip

\preuve
Le param\`etre de r\'egularit\'e $\frac{3}{ p} - 1$ est ajust\'e
pour que ces espaces soient invariants. Ce sont des espaces
fonctionnels : voir [Tr].

Ici, $E$ d\'esigne indiff\'eremment
l'un des $\dot B^{\frac{3}{ p} - 1,q}_{p}$ ou des
$\dot F^{\frac{3}{ p} - 1,q}_{p},\ 1\le p <
\infty$ fix\'e, $1\le q\le
\infty$ (Besov) ou $1\le q < \infty$ (Triebel-Lizorkin). On
utilisera la remarque importante que 
pour tout   $k \in \SZ$ et
$f
\in E$ avec 
$\supp\widehat f\subset
\Gamma_{k}$ on a  $\|f\|_E \sim 2^{k(\frac{3}{ p} - 1)}\Vert
f\Vert_{L^p}$ uniform\'ement pour $k \in \SZ$. Noter que
l'indice $q$ n'appara\^{\i}t pas. Cela r\'esulte de
l'\'equivalence des normes sur $E$ obtenues en partant de deux
analyses de Littlewood-Paley diff\'erentes.

Soient $f,g\in E$, $\supp\widehat f\subset
\Gamma_{k}$,\ $\supp\widehat g\subset \Gamma_{l}$, $
k,l\in \SZ$ et
$j\in \SZ$.

Si $\vert k - l\vert\ge 3$, et par exemple $l\le k - 3$, alors $fg =
\widetilde\Delta_{k} (fg)$, donc
\begin{align*}
\Vert fg\Vert_{E}&\le C\ 2^{k(\frac{3}{ p} - 1)}\ \Vert
fg\Vert_{L^p}\\ &\le \ C\ 2^{k(\frac{3}{ p} - 1)}\ \Vert
f\Vert_{L^p}\ \Vert g\Vert_{L^\infty}\\
&\le \ C\ 2^l\ \Vert f\Vert_{E}\ \Vert g\Vert_{E},
\end{align*}
d'apr\`es la Proprosition 13 une fois encore.

Si $\vert k - l\vert \le 2$, on utilise \`a nouveau l'in\'egalit\'e
d'Young. Si $p\ge 2$, on a
\begin{align*}
\Vert \Delta_{j} (fg)\Vert_{E}&\le C\ 2^{j(\frac{3}{ p} - 1)}\
\Vert
\Delta_{j} (fg)\Vert_{L^p}\\
&\le C\ 2^{j(\frac{6}{ p} - 1)}\ \Vert fg\Vert_{L^{p/2}}\\
&\le \ C\ 2^{j(\frac{6}{ p} - 1)}\ \Vert f\Vert_{L^p}\ \Vert
g\Vert_{L^p}\\
&\le \ C\ 2^{(j - k) (\frac{6}{ p} - 1)}\ 2^k\Vert f\Vert_{E}\
\Vert g\Vert_{E}.
\end{align*}
Cela prouve (\ref{eq 19}) et (\ref{eq 20}) avec $\eta_{n} =
2^{-\frac{6n}{ p}}$.\medskip

Si $1\le p < 2$, et toujours lorsque $\vert k - l\vert\le 2$, on
\'ecrit
\begin{align*}
\Vert \Delta_{j} (fg)\Vert_{E}&\le \ C\ 2^{j(\frac{3}{ p} - 1)}\
\Vert
\Delta_{j} (fg)\Vert_{L^p}\\
&\le \ C\ 2^{2j}\ \Vert fg\Vert_{L^1}\\
&\le \ C\ 2^{2j}\ \Vert f\Vert_{L^p}\ \Vert
g\Vert_{L^{p'}}\\
&\le \ C\ 2^{2j - k}\ \Vert f\Vert_{E}\ \Vert g\Vert_{{E}}
\end{align*}
ce qui prouve (19) et (20) avec $\eta_{n} = 2^{-3n}$.
\qed\medskip

\noindent {\bf Remarque :} si $p > 6,\ fg\notin E$ en g\'en\'eral, m\^eme
sous les hypoth\`eses
de (\ref{eq 20}).
\medskip

Dans ces exemples, $p = \infty$ appara\^{\i}t comme valeur
critique. On a en effet le r\'esultat suivant.\medskip

\begin{prop}\label{pro 19}
   Les espaces $\dot B^{-1,q}_{\infty},\ 1\le q\le \infty$, ne sont
    pas de bons espaces. Ils sont compatibles et v\'erifient (\ref{eq 18})
    avec $\eta_{n} = C$ pour une certaine constante $C$, et ce
    choix de la suite $\eta$ est le meilleur possible.
    \end{prop}
    \medskip

    \preuve
    La compatibilit\'e de ces espaces se d\'emontre comme dans le cas
    $p < \infty$, et on trouve $\eta_{n} = C$.

    Soit $\phi\in \Ss$ telle que $\supp\widehat
\phi\subset B_{0}$.
    Si $x_{1}$ d\'esigne la premi\`ere coordonn\'ee de $x\in \SR^3$
    et si $k\in \SN$, on pose
    $$\begin{array}{lll}
    f(x)&=&2^k\ e^{-i2^k x_{1}}\ \phi(x),\\
    g(x)&=&2^k\ e^{i2^k x_{1}}\ \phi(x).
    \end{array}$$
    Alors $f,g\in \dot B^{-1,q}_{\infty}$ pour tout $q$, de normes
    $\Vert \phi\Vert_{L^\infty}$, avec $\supp\widehat f\subset
    \Gamma_{k}$ et de m\^eme pour $g$. On a
    $$\Delta_{0} (fg) = 4^k\ \Delta_{0} (\phi^2),$$
    donc
    $$\Vert \Delta_{0} (fg)\Vert_{\dot B^{-1,q}_{\infty}} = 4^k\ \Vert
    \Delta_{0} (\phi^2)\Vert_{L^\infty}.$$
    Ceci prouve que si $\eta$ est une suite pour laquelle (\ref{eq
    18}) est vraie dans $\dot B^{-1,q}_{\infty}$, alors
    $$\eta_{n}\ge \sup_{\phi}\frac {\Vert
    \Delta_{0}(\phi^2)\Vert_{L^\infty}}{ \Vert
    \phi\Vert^2_{L^\infty}},$$
    pour tout $n\ge 0$. La proposition est d\'emontr\'ee.
    \qed
   \medskip

   Le corollaire 12 \'enonce bien entendu une propri\'et\'e beaucoup
   plus forte de l'espace $\dot B^{-1,\infty}_{\infty}$.\bigskip

\subsection{Espaces de Morrey}$\phantom{x}$\medskip

Si $1 < q \le p < \infty$, une fonction $f\in L^q_{loc}$ est dans
l'espace de Morrey $M^p_{q}$ lorsque
$$\sup\ R^{\frac{3}{ p}} \bigg(\Mean_{B(x_{0},R)}\ \vert
f\vert^q\bigg)^{1/q} < +\infty,$$
o\`u la borne sup\'erieure est prise sur tous les $x_{0}\in \SR^3$
et les
$R > 0$, et o\`u $\mean_{B}$ d\'esigne $\frac{1}{ \vert B
\vert}\int_{B}$.

Les espaces $M^p_{1}$ sont des espaces de mesures : on dit qu'une
mesure de Radon $\nu$ appartient \`a $M^p_{1}$ lorsque
$$\sup_{x_{0},R}\  R^{\frac{3}{ p}}\,
\frac{|\nu|(B(x_{0},R))}{|B(x_{0},R)|}
 < +\infty.
$$

Les espaces de Morrey invariants sont les $M^3_{q}$, $1\le q
\le 3$.\medskip

\begin{prop}\label{pro 20}
Les espaces $M^3_{q},\ 1\le q \le 3$, sont de bons espaces.
\end{prop}\medskip

\preuve
Que ces espaces soient des espaces fonctionnels (d\'efinition 5) n'est pas
\'evident,
mais cependant d\'ej\`a connu : voir par exemple [M].
\medskip

L'in\'egalit\'e (\ref{eq 19}) s'obtient avec l'argument utilis\'e \`a
plusieurs reprises, selon lequel toute fonction born\'ee est un
multiplicateur des espaces de Morrey.\medskip

Soient maintenant $f,g\in M^3_{q},\ 1\le q \le 3$, qu'on suppose de
normes 1 pour simplifier, telles que $\supp\widehat f\subset
\Gamma_{k}$, $ \supp\widehat g\subset \widetilde \Gamma_{k}$, $
k\in
\SZ$, et soit $j\in
\SZ,\ j\le k + 4$. Il s'agit d'estimer
$$I_j= R\bigg(\Mean_{B(x_{0},R)}\ \vert \Delta_{j}
(fg)\vert^q\bigg)^{1/q}$$
uniform\'ement par rapport \`a $x_{0}, R$. \medskip

On commence par d\'emontrer l'existence d'une constante $C$ telle que
\begin{equation}\label{eq 32}
\Vert \Delta_{j} (fg)\Vert_{L^\infty}\le
\begin{cases}%
C\ 4^j,& \mbox{si}\
q\ge 2,\\
 C\ 4^{k(1 - q/2)}\ 4^{jq/2}, & \mbox{si}\ q <
2.
\end{cases}%
\end{equation}
Pour le voir, on \'ecrit, en choisissant $N$ assez grand :
\begin{align*}
\vert \Delta_{j} (fg) (x)\vert&\le C\ 2^{3j}\ \int (1 + 2^j\vert
x - y\vert)^{-N}\ \vert f(y) g(y)\vert dy\\
&\le C\ 2^{3j}\ \int_{B(x,2^{-j})}\ \vert f(y) g(y)\vert dy\\
&\quad+ C\ \displaystyle \sum_{n\ge 0}\ 2^{3j}\ 2^{-nN}\
\int_{B(x,2^{n + 1 - j})\backslash B(x,2^{n - j})}\ \vert f(y)
g(y)\vert dy\\
&\le C\  \sum_{n\ge 0}\ 2^{-n(N-3)}\
\Mean_{B(x,2^{n - j})}\ \vert f(y)
g(y)\vert dy .
\end{align*}

\noindent Si $q\ge 2$, on applique la d\'efinition de $M^3_{q}$
directement :
\begin{align*}
\vert \Delta_{j} (fg) (x)\vert &\le  C\  \sum_{n\ge 0}\ 2^{-n(N
- 3)}\
\bigg(\Mean_{B(x,2^{n  - j})}\ \vert
fg\vert^{q/2}\bigg)^{2/q}\\ &\le   C\ \sum_{n\ge 0}\
2^{-n(N - 3)}\ 4^{j - n}\\  &\le C\ 4^j.
\end{align*}
pourvu que $N > 1$.\medskip

\noindent Si $q < 2$, on utilise la Proposition 13, qui implique ici
$\Vert f\Vert_{L^\infty}\le C\ 2^k$ et de m\^eme pour $g$, en
\'ecrivant
\begin{equation*}
\vert \Delta_{j} (fg) (x)\vert \le 
 C\ \Vert fg\Vert^{1 - q/2}_{L^\infty}\  \sum_{n\ge 0}\
2^{-n(N - 3)}\ \Mean_{B(x,2^{n  - j)}}\ \vert fg\vert^{q/2}.
\end{equation*}
Des calculs analogues aux pr\'ec\'edents donnent le
r\'esultat d\`es que $N>3-q/2$.\medskip

Revenons \`a $I_j$. Dans le cas o\`u $R\le 2^{-j}$, (\ref{eq 32})
donne
\begin{equation*}
I_j \le\Vert
\Delta_{j} (fg)\Vert_{L^\infty} 2^{-j} \le 
\begin{cases}
 C\ 4^{j - k}\ 4^k\ 2^{-j}& \mbox{si}\ q\ge 2,\\
 C\ 4^{(j - k)q/2}\ 4^k\ 2^{-j}& \mbox{si}\ q < 2.
\end{cases}
\end{equation*}
Dans le cas o\`u $R > 2^{-j}$ et $q\ge 2$, on commence par estimer
\begin{align*}
\vert\Delta_{j} (fg)(x)\vert^{q/2}
&\le  \bigg\{C\ \displaystyle \int 2^{3j}(1 + 2^j \vert x -
y\vert)^{-N}\ \vert
f(y) g(y)\vert dy\bigg\}^{q/2}\\
&\le C\ 
\displaystyle \int 2^{3j} (1 + 2^j \vert
x - y\vert)^{-N}\ \vert  f(y) g(y)\vert^{q/2} dy \\
&\le C\  \displaystyle
\int_{\vert x - y\vert\le
R} 2^{3j} (1 + 2^j \vert  x - y\vert)^{-N}\ \vert  f(y) g(y)\vert^{q/2}
dy \\
&\quad + C\displaystyle \sum_{n\ge 0}\  \int_{2^nR\le\vert x -
y\vert\le 2^{n + 1}R} 2^{3j} (2^{j + n}R)^{-N}\vert f(y)
g(y)\vert^{q/2} dy 
\end{align*}
d'o\`u l'on d\'eduit, si $N>3$, que
\begin{align*}
\Mean_{B(x_{0},R)}\vert\Delta_{j} (fg)\vert^{q/2} &\le
  C\ \displaystyle \sum_{n\ge 0} (2^{j + n}R)^{-(N - 3)}\
\Mean_{B(0, 2^{n + 2}R)}\ \vert fg\vert^{q/2}\\
&\le C\ (1 + 2^jR)^{-N + 3}\ R^{-q}\\
&\le C\ R^{-q}.
\end{align*}
On interpole ensuite cette in\'egalit\'e avec (\ref{eq 32}), pour
obtenir
\begin{align*}
I_j&\le \ R\ \Vert
\Delta_{j} (fg)\Vert_{L^\infty}^{1/2}\
\bigg(\displaystyle\Mean_{B(x_{0},R)}
\vert \Delta_{j}
(fg)\vert^{q/2}\bigg)^{1/q}\\
&\le \ C\ 2^j \\
&= \  C\ 4^{j - k}\ 4^k\ 2^{-j}.
\end{align*}

Dans le cas o\`u $R > 2^{-j}$ et $q < 2$, on commence par observer
comme ci-dessus que 
\begin{align*}
\displaystyle\Mean_{B(x_{0},R)}\ \vert \Delta_{j}
(fg)\vert
\le C\ \displaystyle \sum_{n\ge 0}\ (2^{j + n} R)^{-N + 3}\
\Mean_{B(x_{0},2^{n + 2}R)} \vert fg \vert.
\end{align*}
Ecrivant ensuite
$$\vert f(y) g(y) \vert\le \Vert fg \Vert_{L^\infty}^{1 - q/2}\ \vert
f(y) g(y) \vert^{q/2},$$
on obtient
$$\Mean_{B(x_{0},R)}\ \vert \Delta_{j} (fg) \vert\le C\ 4^{k(1 - q/2)}
R^{-q}.$$
Enfin, on interpole avec (\ref{eq 32}) :
\begin{align*}
I_j &\le \ R\ \Vert
\Delta_{j} (fg)\Vert_{L^\infty}^{1 - 1/q}\
\bigg(\Mean_{B(x_{0},R)}\
\vert
\Delta_{j} (fg)\vert \bigg)^{1/q}\\
&\le\ C\ 4^{k(1 - q/2)}\ 2^{j(q - 1)} \\
& =\ C\ 2^{(j - k)q}\ 4^k\
2^{-j}.
\end{align*}
\medskip

Au total, on a prouv\'e (\ref{eq 20}), avec $\eta_{n} = 4^{-n}$ si
$q\ge 2$ et $\eta_{n} = 2^{-qn}$ si $q < 2$, ce qui cl\^ot la
d\'emonstration.
\qed
\medskip

Le m\^eme type de calculs sera utilis\'e dans la partie suivante,
o\`u l'espace de Morrey $M^3_{2}$ joue implicitement un
r\^ole.\bigskip

\subsection{Espaces de Konozo et Yamazaki}$\phantom{x}$\medskip

Il s'agit d'espaces de Besov au-dessus d'espaces de Morrey, 
not\'es ici $\dot B^{s,r}_{M^p_q}$, $s \in \SR$, $1 \le q \le p
<\infty$ et
$1\le r \le \infty$ (${\mathcal N}^s_{p,q,r}$ dans [Ko,Y]). La
norme se calcule par
$$
\big\|(2^{js}\|\widetilde\Delta_jf\|_{M^p_q}\big)_j\|_{\ell^r}.
$$

Ceux  d'entre eux qui sont invariants sont
caract\'eris\'es par $s=\frac{3}{p} -1$.

\begin{prop}\label{pro KonoetYama} Pour $1 \le q \le p <\infty$
 et $1\le r \le \infty$, 
    les espaces $\dot B^{\frac{3}{p} -1,r}_{M^p_q}$ sont de bons
espaces.
\end{prop}

La preuve est une adaptation des calculs des deux
sections pr\'ec\'edentes en utilisant que si $E=\dot
B^{\frac{3}{p} -1,r}_{M^p_q}$, $f \in E$ et $k\in \SZ$,
$$\|\Delta_kf \|_E = 2^{k(\frac{3}{p} -1)}
\|\Delta_kf\|_{M^p_q}.$$ 
On trouve
$$\eta_n= \begin{cases} 4^{-\frac{3n}{p}},  &\quad   q\ge 2,\\
4^{-\frac{3nq}{2p}},  &\quad   q\le 2.
\end{cases}
$$
Les d\'etails sont laiss\'es au lecteur.

\subsection{Espaces d\'efinis par des conditions sur la transform\'ee
de Fourier}$\phantom{x}$\medskip

Le Jan et Sznitman ont r\'esolu l'\'equation (\ref{eq 1}) avec
donn\'ee initiale dans l'espace des distributions temp\'er\'ees $f$
telles que
$\vert \xi
\vert^2\ \widehat f (\xi)$ soit born\'e sur $\SR^3$ ([Le,Sz]).
Cet espace est caract\'eris\'e par la condition
$$\sup_{j\in \SZ} 4^j\ \Vert \widehat{\Delta_{j} f}\Vert_{L^\infty} <
+\infty.$$\medskip

L'in\'egalit\'e
$$\Vert \widehat{\Delta_{j} f} \Vert_{L^q}\le C\ 2^{3j/q}\ \Vert
\widehat{\Delta_{j} f} \Vert_{L^\infty}$$
sugg\`ere une g\'en\'eralisation. Si $1\le q\le \infty$, on note
$F_{q}$ l'espace des distributions temp\'er\'ees $f$ telles que
$$\Vert f\Vert_{F_{q}} = \sup_{j\in \SZ}\ 2^{j(2 - 3/q)}\ \Vert
\widehat{\Delta_{j} f} \Vert_{L^q} < +\infty.$$
Remarquer que $F_{2} = \dot B^{1/2,\infty}_{2}$, et que $F_{1}$ est
l'espace de Besov $\dot B^{-1,\infty}_{\Aa}$ construit au-dessus de
l'alg\`ebre de Wiener $\Aa$.\medskip

\begin{prop}\label{pro 21}
    Les espaces $F_{q}$ sont de bons espaces pour $1 < q\le \infty$.
\end{prop}

\preuve
On laisse au lecteur le soin de v\'erifier que ce sont des espaces
fonctionnels invariants.

Soient $f,g\in F_{q}$, $\supp\widehat f\subset \Gamma_{k}$,
$\supp\widehat g\subset \Gamma_{l}$, $ k,l\in \SZ$, et soit $j\in
\SZ$.

Si $\vert k - l \vert\ge 3$, et par exemple $l\le k - 3$, on \'ecrit,
puisque $fg = \widetilde \Delta_{k} (fg)$,
\begin{align*}
\Vert fg \Vert_{F_{q}} &\le  C\ 2^{k(2 - 3/q)}\ \Vert \widehat
f\ast
\widehat g\Vert_{L^q}\\
&\le C\ 2^{k(2 - 3/q)}\ \Vert \widehat f\,\Vert_{L^q}\ \Vert
\widehat g
\Vert_{L^1}\\
&\le C\ \Vert f \Vert_{F_{q}}\ 2^{3l/q'}\ \Vert \widehat g
\Vert_{L^q}\\ 
&\le C\ 2^l\ \Vert f\Vert_{F_{q}}\ \Vert g
\Vert_{F_{q}},
\end{align*}
ce qui d\'emontre (\ref{eq 19}).\medskip

Si $\vert k - l \vert\le 2$, on a
\begin{align*}
\Vert \Delta_{j} (fg) \Vert_{F_{q}}&\le C\ 2^{j(2 - 3/q)}\ \Vert
\widehat f\ast \widehat g \Vert_{L^q}\\
 &\le C\ 2^{j(2 - 3/q)}\ \Vert \widehat f\,\Vert_{L^1}\ \Vert
\widehat g
\Vert_{L^q}\\
&\le C\ 2^{k(2 - 3/q)}\ 2^{3k/q'}\Vert \widehat f\, \Vert_{L^q}\
\Vert \widehat g
\Vert_{L^q}\\
&\le C\ 2^{(j - k) 3/q'}\ 4^k\ 2^{-j}\ \Vert f\Vert_{F_{q}}\
\Vert g
\Vert_{F_{q}}.
\end{align*}
Cela signifie que (\ref{eq 20}) est vraie, avec $\eta_{n} =
2^{-3n/q'}$, et termine la preuve.
\qed

\noindent {\bf Remarque :} on voit que (\ref{eq 20}) est vraie dans
$F_{1}$, avec $\eta_{n} = 1$. Le contre-exemple de la Proposition
19 permet \'egalement de montrer qu'il n'y a pas de meilleur
choix de $\eta$ possible, et par cons\'equent que $F_{1}$ n'est pas
un bon espace.\bigskip

\subsection{Discussion de quelques conditions abstraites}$\phantom{x}$\medskip

Plusieurs auteurs ont donn\'e des conditions abstraites que doit ou
peut poss\'eder un espace $E$ pour qu'il soit possible de r\'esoudre
(\ref{eq 1}) par la m\'ethode KW, sp\'ecialis\'ee par un choix
particulier de l'espace $F$ (g\'en\'eralement $E$ lui-m\^eme ou un
sous-espace de $E$).\medskip

Si par exemple on cherche \`a prouver la continuit\'e de
l'application $B$ sur $\Cc([0,\infty[ ; E)$, alors il est n\'ecessaire
d'avoir
\begin{equation}\label{eq 33}
\Vert (-\Delta)^{-1/2} fg \Vert_{E}\le C\ \Vert f \Vert_{E}\ \Vert g
\Vert_{E}
\end{equation}
pour tous $f,g\in E$ (voir \cite{bib M}). Mais cela n'est pas
suffisant, comme l'a montr\'e Oru (\cite{bib O}) avec $E = L^3$.
Noter que (\ref{eq 33}) est plus forte qu'une condition propos\'ee
par Cannone \`a la fin de sa th\`ese (\cite{bib C}), s'\'ecrivant
\begin{equation}\label{eq 34}
\Vert \Delta_{j} (fg) \Vert_{E}\le C\ 2^j\ \Vert f \Vert_{E}\ \Vert g
\Vert_{E}.
\end{equation}
En la modifiant convenablement, Furioli, Lemari\'e et Terraneo ont
toutefois obtenu la continuit\'e de $B$ sur $\Cc([0,\infty[ ; \dot
B^{0,\infty}_{E})$ ([F,LR,Te]). Leur condition est la suivante
\begin{equation}\label{eq 35}
\Vert \Delta_{j} (fg) \Vert_{E}\le C\ 2^j\ \displaystyle \sup_{j'}\ \Vert
\Delta_{j'}f \Vert_{E}\ \sup_{j'}\ \Vert \Delta_{j'}g \Vert_{E}.
\end{equation}
C'est (\ref{eq 34}) dans l'espace $\dot B^{0,\infty}_{E}$, et aussi
(\ref{eq 33}).\medskip

Meyer et Muschietti d'une part, Furioli, Lemari\'e et Terraneo
d'autre part, ont donn\'e d'autres conditions permettant de
r\'esoudre (\ref{eq 1}) par la m\'ethode de Kato :
\begin{equation}\label{eq 36}
\Vert fg \Vert_{E}\le C\ (\Vert f \Vert_{L^\infty} + \Vert
(-\Delta)^{1/2} f \Vert_{E})\ \Vert g \Vert_{E}
\end{equation}
est \'etudi\'ee dans [M],
\begin{equation}\label{eq 37}
\Vert \Delta_{j} (fg)\Vert_{E}\le C\ (\Vert f \Vert_{E}\ \Vert g
\Vert_{L^\infty} + \Vert f \Vert_{L^\infty}\ \Vert g\Vert_{E})
\end{equation}
a \'et\'e propos\'ee dans [F,LR,Te], ainsi que la forme
l\'eg\`erement plus pr\'ecise
\begin{equation}\label{eq 38}
\Vert \Delta_{j} (fg) \Vert_{E}\le C\  \Vert f \Vert_{E}\ \Vert
g\Vert_{L^\infty}
+ C\  2^{j(1 - \alpha)} \Vert f \Vert^{1 - \alpha}_{E}\
\Vert f\Vert^\alpha_{L^\infty}\ \Vert g\Vert_{E}\ ,
\end{equation}
o\`u $\alpha\in ]0,1[$.\medskip

\begin{prop}\label{pr 22}
Tout espace fonctionnel invariant v\'erifiant l'une des
propri\'et\'es (\ref{eq 33}) \`a (\ref{eq 38}) est un bon espace, au
sens de la variante de la Proposition 15.
\end{prop}
\medskip

\preuve
Soient d'abord $f,g\in E$, $\supp\widehat f\subset \Gamma_{k}$ et
$\supp\widehat g\subset B_{k - 2}$, $ k\in \SZ$. On a alors $fg =
\widetilde \Delta_{k} (fg)$, et d'apr\`es (\ref{eq 21}),
\begin{align*}
\ \ \ \ \ \ \ \ fg &=2^k\ \widetilde \Delta_{k, -1}\
((-\Delta)^{-1/2} fg),\\ (-\Delta)^{1/2} f&=  2^k\ \widetilde
\Delta_{k,-1} (f).
\end{align*}
Ces identit\'es inject\'ees dans (\ref{eq 33} - \ref{eq 38}) 
donnent
$$\Vert fg\Vert_{E}\le C\ 2^k\ \Vert f \Vert_{E}\ \Vert g\Vert_{E},
$$
\'eventuellement en utilisant aussi la Proposition 13
(c'est-\`a-dire $\|f\|_{L^\infty} \le C\, 2^k\, \|f\|_E$).

Si ensuite l'hypoth\`ese sur $g$ est chang\'ee en $\supp\widehat
g\subset \widetilde \Gamma_{k}$, et si $j\in \SZ$, alors on obtient
(\ref{eq 20}) par des consid\'erations analogues, avec $\eta_{n} =
4^{-n}$ sous les conditions (\ref{eq 33} - \ref{eq 34} - \ref{eq 35}),
et $\eta_{n} = 2^{-n}$ sous (\ref{eq 36} - \ref{eq 37} - \ref{eq 38}).
Les d\'etails sont laiss\'es au lecteur.
\qed
\medskip

\noindent {\bf Remarque :} tout espace fonctionnel invariant
v\'erifiant (\ref{eq 36}) est \'ega\-le\-ment un bon espace au sens
de la D\'efinition 7, puisque (\ref{eq 19}) est vraie dans ce cas.



\section{Les espaces $M(\eta)$ et les contre-exemples}\bigskip

Le but principal de cette partie est de d\'emontrer le th\'eor\`eme 11.
Les espaces $E$ qui vont fournir les contre-exemples d\'esir\'es
appartiennent \`a une m\^eme classe, celle des espaces
$M(\eta)$.\bigskip

\subsection{Les espaces $M(\eta)$}$\phantom{x}$\medskip

\begin{defin}\label{defin 23}
Soit $\eta = (\eta_{n})_{n\in \SZ}$ une suite positive, d\'ecroissante
et r\'eguli\`ere. L'espace $M(\eta)$ est l'ensemble des distributions
temp\'er\'ees $f$ telles que $\displaystyle \lim_{j\rightarrow -\infty}
S_{j}f = 0$ dans $\Ss'$, et pour lesquelles il existe une constante
$C$ v\'erifiant pour tous $x_{0}\in \SR^3, j, j'\in \SZ$ avec $j\ge
j'$
\begin{equation}\label{eq 39}
\displaystyle\Mean_{B(x_{0},2^{-j'})}\ \vert \Delta_{j}
f\vert^2\le C^2\
\eta^2_{j - j'}\ 4^j,
\end{equation}
o\`u $\mean_{B}$ d\'esigne $\displaystyle \frac{1}{
\vert B
\vert}\int_{B}$.
La norme dans $M(\eta)$ est la plus petite constante $C$ possible.
\end{defin}
\medskip

La r\'egularit\'e de la suite $\eta$ implique que $M(\eta)$ ne
d\'epend pas d'un choix particulier des op\'erateurs $\Delta_{j}$.
Des in\'egalit\'es (\ref{eq 39}) avec $j' = j$ et de (\ref{eq 15}),
on d\'eduit que $M(\eta)$ s'injecte dans $\dot B^{-1,\infty}_{\infty}$.
Il s'ensuit que les valeurs de $\eta$ pour les entiers $n\le 0$ n'ont
pas d'importance : on suppose dor\'enavant $\eta_{n} = \eta_{0}$ si
$n\le 0$, de sorte que (\ref{eq 39}) est vraie quels que soient $j$
et $j'$.\medskip

Si $\displaystyle\liminf 2^{3n}\,
\eta^2_{n} = 0$, alors $M(\eta)$ est r\'eduit \`a $\{0\}$, de sorte
qu'on suppose \'egalement $\displaystyle\liminf 2^{3n}\,
\eta^2_{n} > 0$. 

Lorsque $\eta_{n} = 2^{-n\alpha}$, $ n\ge 0$,
pour un certain
$\alpha\in [0,3/2]$, $ M(\eta)$ est l'espace de Besov homog\`ene $\dot
B^{\alpha - 1,\infty}_{M_{2}^{3/\alpha}}$ (voir section 3.5). En particulier,
il co\"{\i}ncide avec $\dot B^{1/2,\infty}_{2}$ pour $\alpha =
3/2$, et avec $\dot B^{-1,\infty}_{\infty}$ pour $\alpha =
0$. 

 Dans le cas o\`u $\eta$ est quelconque, on
a toujours
$$\dot B^{1/2,\infty}_{2}\subset M(\eta)\subset \dot
B^{-1,\infty}_{\infty},$$
avec injections continues. Si la suite $\eta$ est telle que
$\displaystyle \sum_{n\ge 0}\ \eta^2_{n} < +\infty$, alors $M(\eta)$
s'injecte contin\^ument dans $\dot F^{-1,\infty}_{2}$, l'espace des
d\'eriv\'ees des fonctions de BMO.\medskip

\begin{prop}\label{prop 24}
Les espaces $M(\eta)$ sont invariants, et v\'erifient 
(\ref{eq 18})
avec la suite $C\eta = (C\eta_{n})_{n\in \SZ}$, pour une certaine
constante absolue $C$.
\end{prop}\medskip

\preuve
L'invariance des $M(\eta)$ est ais\'ee \`a montrer et laiss\'ee au
lecteur. V\'erifions (\ref{eq 18}).

Soient $f,g\in M(\eta)$, $\supp\widehat f\subset \Gamma_{k}$,
 $
\supp\widehat g\subset \Gamma_{l}$, $ k,l\in \SZ$.

Si $\vert k - l\vert \ge 3$, et par exemple $l\le k - 3$, alors $fg =
\widetilde \Delta_{k} (fg)$. Soient $j'\le k$ et $x_{0}\in \SR^3$. On a
pour tout $n\ge 0$
\begin{equation}\label{eq 40}
\displaystyle\Mean_{B(x_{0},2^{-j' + n})}\ \vert fg \vert^2\le
C\
\eta^2_{j + n -
j'}\ \Vert f\Vert^2_{M(\eta)}\ \Vert g\Vert^2_{L^\infty}.
\end{equation}
Si $x\in \SR^3$ est momentan\'ement fix\'e, et en choisissant $p >
3$, on a
$$\vert \Delta_{k} (fg) (x)\vert\le C\ 2^{3k}\ \displaystyle\int (1 +
2^k\vert x -
y\vert)^{-p}\ \vert f(y) g(y)\vert\ dy,$$
d'o\`u
\begin{align*}
\vert \Delta_{k} (fg) (x)\vert^2&\le C\ 2^{3k}\ 
\displaystyle\int (1 + 2^k\vert x -
y\vert)^{-p}\ \vert f(y) g(y)\vert^2 dy\\
&\le C\ 2^{3k}\ \displaystyle\int _{B(x,2^{-j'})} (1 + 2^k \vert
x - y\vert)^{-p}\
\vert f(y) g(y)\vert^2 dy\\
&\quad + C\ \displaystyle \sum_{n\ge 0} 2^{(3 - p) (k - j' + n)}\
\Mean_{B(x,2^{-j' + n + 1})}\
\vert fg\vert^2.
\end{align*}
On utilise (\ref{eq 40}) et la d\'ecroissance de $\eta$ :
\begin{align*}
\vert \Delta_{k} (fg) (x)\vert^2 &\le C\ 2^{3k}\
\displaystyle\int_{B(x,2^{-j'})}\ (1 +
2^k\vert x - y\vert)^{-p}\ \vert f(y) g(y)\vert^2 dy\\
&\quad+ C\ \eta^2_{k - j'}\ \Vert f\Vert^2_{M(\eta)}\ \Vert
g\Vert^2_{L^\infty}.
\end{align*}
On applique cette in\'egalit\'e \`a tout $x\in
B(x_{0},2^{-j'})$ et on utilise \`a nouveau (\ref{eq 40}) :
\begin{align*}
\displaystyle\Mean_{B(x_{0},2^{-j'})}\ \vert \Delta_{k} (fg)
(x)\vert^2 dx \le 
 C\ \eta_{k - j'}^2\ \Vert f \Vert_{M(\eta)}^2\ \Vert g
\Vert^2_{L^\infty}.
\end{align*}

On obtient de m\^eme l'in\'egalit\'e analogue pour $\Delta_{k - 1}
(fg)$ et $\Delta_{k + 1} (fg)$. Compte tenu de la Proposition 13,
ceci d\'emontre (\ref{eq 19}).\medskip

On suppose maintenant $\vert k - l\vert\le 2$, et on se donne $j\in
\SZ,\ j\le k + 4$.

Si $x\in \SR^3$, on a
\begin{align*}
\vert \Delta_{j} (fg) (x)\vert \le C\ 2^{3j}\
\displaystyle\int_{B(x,2^{-j})}\ \vert fg\vert
+C\ \displaystyle\sum_{n\ge 0}\ 2^{(3 - p)n}\
\displaystyle\Mean_{B(x,2^{-j + n +
1})}\vert fg\vert,
\end{align*}
o\`u on a choisi $p > 3$ comme plus haut. La d\'efinition de $M(\eta)$
et la d\'ecroissance de $\eta$ impliquent
\begin{equation}\label{eq 41}
\vert \Delta_{j} (fg) (x)\vert\le C\ \eta^2_{k - j}\ 4^k\ \Vert
f\Vert_{M(\eta)}\ \Vert g\Vert_{M(\eta)}.
\end{equation}
Par ailleurs, si $x_{0}\in \SR^3$ et $j'\le j$ sont fix\'es
\begin{align*}
 \vert\Delta_{j} (fg)
(x)\vert  \le C\ 
\int 2^{3j}\ (1 + 2^j
\vert x - y\vert)^{-p}\ \vert f(y) g(y)\vert dy dx
\end{align*}
d'o\`u
\begin{align*}
\Mean_{B(x_{0},2^{-j'})}\ \vert\Delta_{j} (fg)
(x)\vert dx 
\le C\ 2^{3j'}\ \displaystyle\int (1 + 2^{j'} \vert x_{0} -
y\vert)^{-p}\ \vert f(y) g(y)\vert dy.
\end{align*}
La m\^eme technique de d\'ecoupage du domaine
d'int\'egration, en utilisant cette fois les boules
$B(x_{0},2^{-j' + n}),\ n\ge 0$, donne
$$\displaystyle\Mean_{B(x_{0},2^{-j'})}\ \vert\Delta_{j}
(fg)\vert\le C\
\eta^2_{k -
j'}\ 4^k\ \Vert f\Vert_{M(\eta)}\ \Vert g \Vert_{M(\eta)}.$$
Avec (\ref{eq 41}), ceci implique
$$\displaystyle\Mean_{B(x_{0},2^{-j'})}\ \vert\Delta_{j}
(fg)\vert^2\le C\
\eta^2_{k -
j}\ \eta^2_{k - j'}\ 4^{2k}\ \Vert f\Vert^2_{M(\eta)}\ \Vert g
\Vert^2_{M(\eta)}.$$
Puisque $\eta_{k - j'}\le \eta_{j - j'}$, on obtient
$\Delta_{j}(fg)\in M(\eta)$, et
$$\Vert\Delta_{j} (fg)\Vert_{M(\eta)}\le C\ \eta_{k - j}\ 4^k\ 2^{-j}\ \Vert
f\Vert_{M(\eta)}\ \Vert g \Vert_{M(\eta)}.$$
La preuve est termin\'ee.
\qed\medskip

\noindent{\bf Remarque :} on pose, pour tout $n\ge -4$
$$\varepsilon_{n} = \sup_{m\ge 0}\ \frac{\eta_{n + m}}{\eta
_{m}}.$$
 Les calculs pr\'ec\'edents montrent qu'on a
l'in\'egalit\'e plus forte
$$\Vert\Delta_{j} (fg)\Vert_{M(\eta)}\le C\ \eta_{k - j}\
\varepsilon_{k - j}\ 4^k\ 2^{-j}\ \Vert
f\Vert_{M(\eta)}\ \Vert g \Vert_{M(\eta)}.$$
L'am\'elioration est sensible dans certains cas, par exemple pour
$\eta_{n} = 2^{-n\alpha}$. Mais pour les suites $\eta\notin l^2(\SN)$
qui vont fournir les contre-exemples, l'in\'egalit\'e
$\varepsilon_{n}\le C$ est la seule raisonnable.\bigskip

\subsection{Construction des contre-exemples}$\phantom{x}$\medskip

Soit $\eta$ une suite positive d\'ecroissante et r\'eguli\`ere telle que
$$\displaystyle \sum_{n\ge 0}\ 2^{-3n}\displaystyle\inf_{k\ge n}\ 2^{3k}\
\eta_{k}^2
= +\infty.$$

La preuve du th\'eor\`eme 11 se fait par l'absurde en consid\'erant l'espace
fonctionnel compatible
$M_{0}(\eta)$, fermeture de $\Ss$ dans
$M(\eta)$, et en supposant qu'il peut
\^etre assorti d'un espace $F$ rendant le couple
$(M_{0}(\eta),F)$ admissible.\medskip

La contradiction r\'esultera des deux lemmes suivants.\eject

\begin{lem}\label{lem 25}
On peut trouver $m_{0}\in \Ss,\ \psi\in \Ss,\ x\in \SR^3$ et
$\Lambda\subset \SZ^3 - \{0\}$ tels que, si
$$u(1,\cdot) = \sum_{l\in \Lambda}\ d_{l}\ m(\cdot - l),$$
o\`u $m = e^\Delta m_{0}$ et $(d_{l})_{l\in \Lambda}$ est une famille
quelconque de coefficients presque tous nuls, si
$$u(t,\cdot) = \frac{1}{ \sqrt t}\ u\bigg(1,\frac{\cdot}{\sqrt
t}\bigg),\ t > 0,$$
 alors il existe pour tout $l$ un coefficient
$\varepsilon_{l}\in
\{-1,+1\}$ tel que, en posant
$$v(1,\cdot) = \displaystyle\sum_{l\in \Lambda}\ \overline
d_{l}\
\varepsilon_{l}\ m(\cdot - l),$$
$$v(t,\cdot) = \displaystyle\frac{1}{ \sqrt t}\
v\bigg(1,\frac{\cdot}{\sqrt t}\bigg),\ t > 0,$$
on ait
\begin{equation}\label{eq 42}
\displaystyle \sum_{l\in \Lambda}\ \frac{\vert d_{l}\vert^2}{
\vert l\vert^3}\le C\ \vert \psi\ast B(u,v) (1) (x)\vert + C\
\Vert d_{l}\Vert^2_{l^\infty},
\end{equation}
pour une certaine constante $C$ ind\'ependante de la famille
$(d_{l})_{l\in \Lambda}$. De plus, on peut choisir $\vert x\vert$ aussi
grand qu'on veut.
\end{lem}
\medskip

\begin{lem}\label{lem 26}
Avec les notations pr\'ec\'edentes, et si $\vert x\vert$ est assez
grand, alors d\`es que $\displaystyle \sum_{n\ge 0}\ 2^{-3n}\,
\displaystyle\inf_{k\ge n}\
2^{3k}\, \eta_{k}^2 = +\infty$, on peut trouver une famille $(c_{l})_{l\in
\Lambda}$ dans
$l^\infty (\Lambda)$ telle que
\begin{equation}\label{eq 43}
\displaystyle \sum_{l\in \Lambda}\ \frac{\vert c_{l}\vert^2}{
\vert l\vert^3} = +\infty,
\end{equation}
\begin{equation}\label{eq 44}
\mbox{pour toute famille}\ (\varepsilon_{l})_{l\in \Lambda}\in
\{-1,1\}^\Lambda,\ \mbox{on a}\ \displaystyle\sum_{l\in \Lambda}\
\varepsilon_{l} c_{l}\ m_{0}
(\cdot - l)
\end{equation}
$\in M(\eta)$, {\it o\`u} $m_{0}$ {\it est la fonction du lemme \ref{lem 25},
et ce uniform\'ement par rapport \`a} $(\varepsilon_{l})_{l\in \Lambda}$.

\end{lem}
\medskip

Admettant pour le moment ces deux r\'esultats, on revient au
th\'eor\`eme 11 o\`u l'on a suppos\'e  le couple
$(M_{0}(\eta),F)$ admissible.
Soit $(c_{l})_{l\in \Lambda}$ la famille donn\'ee par le lemme \ref{lem 26},
o\`u $\Lambda$ est l'ensemble donn\'e par le lemme \ref{lem 25}. Pour tout
$N$ entier positif, on pose $d_{l} = c_{l}$ si
$l\in \Lambda$ et $\vert l\vert\le N,\ d_{l} = 0$ sinon. Ensuite, on d\'efinit
les fonctions $u_{0}$ et $v_{0}$ (d\'ependant de $N$) par
$$u_{0} (y) = \displaystyle\sum_{l}\ d_{l}\ m_{0} (y - l)\ ,\ \
v_{0} (y) = \displaystyle \sum_{l} \overline d_{l}\ \varepsilon_{l}\
m_{0} (y - l).$$\medskip

Il r\'esulte du lemme \ref{lem 26} que $u_{0},\ v_{0}\in M_{0}(\eta)$
uniform\'ement par rapport \`a $N$.\medskip

Soient alors $u$ et $v$ d\'efinies comme au lemme \ref{lem 25}.
Puisque $u(1,\cdot) = e^\Delta u_{0}$, et de m\^eme pour
$v(1,\cdot)$, on a
$u(1,\cdot),\ v(1,\cdot)\in F$, ce qui implique $u,v \in
\Ff$. La continuit\'e de $B$ sur $\Ff$ entra\^{\i}ne en
particulier que
$B(u,v) (1)\in \Ss'$. Plus pr\'ecis\'ement, il existe une
constante $C$ ne d\'ependant pas de $N$ telle que
\begin{align*}
\vert \psi\ast B(u,v) (1) (x) \vert&\le C\ \Vert B (u,v)
(1)\Vert_{F}\\ &\le C\ \Vert B\Vert\ \Vert u\Vert_{\Ff}\ \Vert
v\Vert_{\Ff}\\ &\le C\ \Vert B\Vert\ \Vert u(1,\cdot)\Vert_{F}\
\Vert v(1,\cdot)\Vert_{F}\\ &\le C\ \Vert B\Vert\ \ \Vert
e^\Delta
\Vert_{F,M(\eta)}^2\ \Vert u_{0}\Vert_{M(\eta)}\ \Vert
v_{0}\Vert_{M(\eta)},
\end{align*}
o\`u $x$ est le point donn\'e par le lemme \ref{lem 25}.\medskip

On d\'eduit alors de (\ref{eq 42}) et de ce qui pr\'ec\`ede
l'existence d'une constante $M$ ind\'ependante de $N$ telle que
$$\sum_{l\in \Lambda,  \vert
l\vert\le N}\
\displaystyle\frac{\vert c_{l}\vert^2}{ \vert l\vert^3}\le
M.$$ Ceci est contradictoire avec (\ref{eq 43}), et ram\`ene
la preuve du th\'eor\`eme 11 \`a celle des lemmes \ref{lem
25} et \ref{lem 26}.\medskip

\noindent {\bf Remarque :}
    Puisque la suite $(c_{l})_{l\in \Lambda}$ est born\'ee, on a
    toujours $u_{0}, v_{0}\in \dot B^{-1,\infty}_{\infty}$,
    uniform\'ement par rapport \`a $N$. Ceci d\'emontre le corollaire 12.

    Par ailleurs, il r\'esulte du lemme \ref{lem
25} qu'on a \'egalement $u_{0},
    v_{0}\in \dot F^{-1,2}_{\infty}$, uniform\'ement par rapport \`a
    $N$ : on ne peut donc pas faire fonctionner la m\'ethode KW
    \`a partir de l'espace $\dot F^{-1,2}_{\infty}$ des d\'eriv\'ees
    des fonctions de $BMO$\footnote{Alors que cet article \'etait en
fin de r\'edaction, nous avons appris que  dans un
travail intitul\'e ``Well posedness for the Navier-Stokes equations''
(preprint, Northwestern University) Koch et Tataru ont d\'emontr\'e
l'existence de solutions, pour petite donn\'ee initiale dans $\dot
F^{-1,2}_{\infty}$,  dans un espace  de
fonctions  localement de carr\'e int\'egrable en temps et espace pour lequel
ils obtiennent la bicontinuit\'e de $B$. Cet espace ne se compare pas au
n\^otre et il n'y a donc pas de contradiction.}
.\bigskip

\subsection{Preuve du lemme \ref{lem
25}}$\phantom{x}$\medskip

Quitte \`a changer de rep\`ere, on suppose que le symbole $P(\xi)$ de
l'op\'erateur $P(D)$ n'est jamais nul dans un c\^one ouvert de sommet
0 et d'axe le demi-axe $0z$ des cotes positives. L'ensemble $\Lambda$
sera contenu dans un c\^one, image du pr\'ec\'edent par une rotation
d'angle $\pi/2$.\medskip

On commence par construire $m_{0}$.\medskip

\begin{lem}\label{lem 27}
Il existe $m_{0}\in \Ss$, paire et \`a valeurs r\'eelles, dont la
transform\'ee de Fourier est support\'ee par une couronne compacte ne
contenant pas 0, telle que les fonctions $m(\cdot - l)$, $l\in \SZ^3$,
forment une famille orthonorm\'ee de $L^2(\SR^3)$, o\`u $m = e^\Delta
m_{0}$.
\end{lem}\medskip

\preuve
On part de $h\in \Ss$, paire et \`a valeurs r\'eelles, telle que
$\widehat h$ soit support\'ee dans $\{2\pi\le \vert \xi \vert\le
8\pi\}$ par exemple, avec
$$\displaystyle \sum_{l\in \SZ^3}\ \vert \widehat h (\xi + 2\pi
l)\vert^2\ge \delta > 0.$$
On pose
$$\sigma(\xi) = \sum_{l\in \SZ^3}\ e^{-2\vert \xi + 2\pi l \vert^2}\
\vert \widehat h (\xi + 2\pi l) \vert^2,$$
et $\widehat m_{0} (\xi) = \sigma (\xi)^{-1/2}\ \widehat h (\xi)$. Par
construction, on a
$$\sum_{l\in \SZ^3}\ \vert \widehat m (\xi + 2\pi l) \vert^2 = 1,$$
d'o\`u on d\'eduit le lemme par la formule de Poisson.
\qed\medskip

On pose a priori
$$f(x) = u(1,x) = \sum_{l}\ d_{l}\ m(x - l),$$
$$g(x) = v(1,x) = \sum_{l}\ e_{l}\ m(x - l),$$
o\`u $(d_{l})$ et $(e_{l})$ sont deux suites \`a support fini. Si $t >
0$, on pose ensuite
$$u(t,x) = \frac{1}{ \sqrt t}\ u\left(1,\frac{x}{ \sqrt
t}\right),$$
$$v(t,x) = \frac{1}{ \sqrt t}\ v\left(1,\frac{x}{ \sqrt
t}\right).$$ On note enfin
$$\psi\ast B (u,v) (1) = L (fg),$$
o\`u $\psi\in \Ss$ est pour le moment quelconque, et $L(x,y)$ le
``noyau''  de $L$, au sens o\`u
$$L(fg) (x) = \displaystyle \int L(x,y)\ f(y)\ g(y)\ dy.$$
\medskip

\begin{lem}\label{lem 28}
$$L(x,y) = \displaystyle \int^1_{0}\ \theta_{1 - \tau} (x - \sqrt \tau
y)\ \sqrt \tau\ d\tau,$$
o\`u les fonctions $\theta_{1 - \tau}$ sont d\'efinies par
\begin{equation}\label{eq 45}
\widehat \theta_{1 - \tau}\ (\xi) = \widehat \psi (\xi)\ P(\xi)\ e^{-(1 -
\tau)\ \vert \xi \vert^2}.
\end{equation}
\end{lem}
\medskip

\preuve
On part de $B(u,v) (1) = \displaystyle \int^1_{0} b(\tau) \, d\tau$ o\`u
$$b(\tau)= e^{(1 - \tau)\Delta}\ P(D)\
\frac{1}{ \tau}\ f\left(\frac{\cdot}{\sqrt \tau}\right)\
g\left(\frac{\cdot}{\sqrt \tau}\right).$$
Soit $h_{0}$ d\'efinie par $\widehat h_{0}(\xi) = P(\xi)\ e^{-\vert
\xi\vert^2}$. On a
\begin{align*}
b(\tau)(z)
&=\displaystyle \int\frac{1}{ (1 - \tau)^2}\
h_{0}\left(\frac{z - y}{
\sqrt{1 - \tau}}\right) \frac{1}{ \tau} f\left(\frac{y}{ \sqrt
\tau}\right) g\left(\frac{y}{ \sqrt \tau}\right) dy\\
&=\displaystyle \int\frac{\sqrt \tau}{ (1 - \tau)^2}\
h_{0}\left(\frac{z -
\sqrt \tau y}{ \sqrt{1 - \tau}}\right) f(y)\ g(y)\ dy.
\end{align*}
On en d\'eduit
\begin{align*}
L(x,y)&=\displaystyle \int^1_{0} \int \psi (x - z)\ 
\frac{\sqrt \tau}{
(1 - \tau)^2}\ h_{0}\left(\frac{z - \sqrt \tau y}{ \sqrt{1 -
\tau}}\right) dz\ d\tau\\
&=\displaystyle \int^1_{0}\ \ \theta_{1 - \tau} (x - \sqrt \tau y)\
\sqrt \tau\ d\tau,
\end{align*}
o\`u
$$\theta_{1 - \tau} (x) = \displaystyle \int\ \psi (x - z)\ 
\frac{1}{ (1 -
\tau)^2}\ h_{0}\left(\frac{z - \sqrt \tau y}{ \sqrt{1 -
\tau}}\right)\ dz.$$
On a bien (\ref{eq 45}), et le lemme est d\'emontr\'e.
\qed
\medskip

On d\'efinit maintenant $\theta$ par
$$\widehat \theta (\xi) = \widehat \psi (\xi)\ P(\xi)\ e^{-\vert \xi
\vert^2},$$
et le noyau $L_{0}$ 
par
$$L_{0}(x,y) = \displaystyle \int^1_{0}\ \theta(x - \sqrt \tau y)\
\sqrt \tau\ d\tau.$$\medskip

\begin{lem}\label{lem 29}
$$L(fg) = L_{0} (fg) + w_{0},$$
o\`u $w_{0}$ est une fonction born\'ee v\'erifiant
$$\Vert w_{0}\Vert_{L^\infty}\le C\ \Vert d_{l}\Vert_{l^\infty}\
\Vert e_{l} \Vert_{l^\infty}.$$
\end{lem}
\medskip

\preuve
On d\'efinit les fonctions $\sigma_{1 - \tau}$ par la formule
$$\theta_{1 - \tau} = \theta + \tau\ \sigma_{1 - \tau},$$
autrement dit
\begin{align*}
\widehat \sigma_{1 - \tau} (\xi)&=\widehat \psi (\xi)\ 
P(\xi)\ \frac{e^{-(1 -
\tau)\vert \xi \vert^2} - e^{-\vert \xi \vert^2}}{
\tau}\\ &= \widehat \psi (\xi)\ P(\xi)\ \vert \xi \vert^2\
\frac{1}{ \tau}\
\displaystyle \int^1_{1 - \tau}\ e^{-s \vert \xi \vert^2}\ ds.
\end{align*}
Si $s\in [0,1]$, alors la fonction dont la transform\'ee de Fourier
est $\widehat \psi (\xi)\ P(\xi)\ \vert \xi \vert^2\ e^{-s \vert \xi
\vert^2}$ est dans $L^1 (\SR^3)$, uniform\'ement par rapport \`a $s$.
Il en r\'esulte
\begin{equation}\label{eq 46}
\Vert \sigma_{1 - \tau}\Vert_{L^1}\le C
\end{equation}
pour tout $\tau$.

La fonction $w_{0}$ s'\'ecrit
$$w_{0}(x) = \displaystyle \int\ L'_{0} (x,y)\ f(y)\ g(y)\ dy,$$
o\`u
$$L'_{0} (x,y) = \displaystyle \int^1_{0}\ \sigma_{1 - \tau}\ (x -
\sqrt \tau y)\ \tau^{3/2}\ d\tau.$$
On d\'eduit de (\ref{eq 46}) que
$$\sup_{x}\ \displaystyle \int\ \vert L'_{0} (x,y) \vert\ dy\le C,$$
ce qui implique
$$\Vert w_{0} \Vert_{L^\infty}\le C\ \Vert fg \Vert_{L^\infty}\le C\
\Vert d_{l}\Vert_{l^\infty}\ \Vert e_{l} \Vert_{l^\infty}.$$
\qed
\medskip

\begin{lem}\label{lem 30}
$$L_{0}(fg) (x) = 2\sum_{l\in \SZ^3} d_{l} e_{l}\ \displaystyle
\int^1_{0} \theta(x - \tau l) \tau^2 d\tau + w_{1},$$
o\`u $w_{1}$ est born\'ee, avec $\Vert w_{1}\Vert_{L^\infty}\le C\
\Vert d_{l}\Vert _{l^\infty}\ \Vert e_{l} \Vert_{l^\infty}$.
\end{lem}
\medskip

\preuve
Partant de
\begin{align*}
L_{0} (fg) (x)&=\displaystyle \int^1_{0}\!\! \int\ \theta(x - \sqrt
\tau y)\ f(y)\ g(y)\ \sqrt \tau\ dy\ d\tau\\
&=2\ \displaystyle \int^1_{0}\!\! \int\ \theta(x - \tau y)\ f(y)\
g(y)\
\tau^2\ dy\ d\tau,
\end{align*}
on d\'eveloppe le produit $f(y)\ g(y)$. Ceci permet de calculer
$w_{1}$ sous la forme $w_{1} = w_{2} + w_{3}$, o\`u
$$w_{2} = 2\sum_{l\in \SR^3} d_{l} e_{l} \displaystyle \int^1_{0}\!\!
\int \tau^2 \{\theta(x -\tau y) - \theta(x - \tau l)\} m^2(y - l) dy
d\tau,$$
et
$$w_{3} = 2\sum_{\substack{\scriptstyle l\in \SR^3\\ \scriptstyle l\not= l'}}
d_{l} e_{l'} \displaystyle \int^1_{0}\!\!\int \theta(x -\tau y) m(y - l)
m(y - l') \tau^2 dy d\tau.$$

On majore $|w_{2}(x)|$ par
$$2 \Vert d_{l}\Vert_{\infty} \Vert
e_{l}\Vert_{\infty}\displaystyle \int^1_{0}\sum_{l}\int \tau^2
\vert\theta(x - \tau y) - \theta(x - \tau l)\vert \vert m(y -
l)\vert^2 dy d\tau.$$
Puisque $m,\ \theta$ et $\nabla \theta$ sont \`a d\'ecroissance
rapide (en choisissant convenablement $\psi$, par exemple de sorte
que $0\notin \supp \widehat \psi$), on a pour tout $p$ assez
grand et
$x,y$ fix\'es
\begin{align*}
&\displaystyle\sum_{l}\ \vert \theta (x - \tau y) - \theta (x - \tau
l)\vert\ \vert
m (y - l)\vert^2\\
&\qquad\le C \tau \displaystyle\sum_{l, \vert y -
l\vert\le\frac{1}{
\tau}}
\vert y - l\vert\
\vert m (y - l)\vert^2\ (1 + \vert x - \tau y\vert)^{-p}\\
&\qquad\quad
+ C\displaystyle\sum_{l, \vert y - l\vert > \frac{1}{
\tau}}
\vert m (y - l)\vert^2\ \{(1 + \vert x - \tau y\vert)^{-p} + (1 +
\vert x - \tau l\vert)^{-p}\}\\
&\qquad\le  C  \tau (1 + \vert x - \tau y\vert)^{-p}\\
&\qquad\quad+ C \displaystyle\sum_{l, \vert y - l\vert > \frac{1}{
\tau}}
\tau^4\ (1 +
\vert y -
l\vert)^{-p}\ (1 + \vert x - \tau l\vert)^{-p}\\
&\qquad\le C \tau (1 + \vert x - \tau y\vert)^{-p}.
\end{align*}
Il en r\'esulte que $\Vert w_{2}\Vert_{L^\infty}\le\ C\ \Vert
d_{l}\Vert_{L^\infty}\ \Vert e_{l}\Vert_{L^\infty}.$\medskip

Pour traiter $w_{3}$, on utilise l'orthogonalit\'e des fonctions
$m(. - l)$ et $m(\cdot - l')$ lorsque $l\not= l'$ (lemme \ref{lem 27}) pour
\'ecrire
$$w_{3}(x) = 2\sum_{\substack{\scriptstyle l,l'\in \SZ^3\\\scriptstyle
l\not= l'}}d_{l} e_{l'}\int^1_{0}\!\!\int \{\theta (x - \tau y) - \theta (x -
\tau l)\}\tau^2 m(y - l) m(y - l') dy d\tau.$$
Comme $ \sum_{l'}\ \vert m (\cdot - l')\vert$ est une
fonction born\'ee, on a
$$\vert w_{3}(x)\vert\le C \Vert d_{l}\Vert_{l^\infty} \Vert
e_{l}\Vert_{l^\infty} \sum_{l}\!\!
\int^1_{0}\!\!\int \vert \theta (x - \tau y) - \theta (x - \tau
l)\vert \tau^2 \vert m(y - l)\vert dy d\tau,$$
et on conclut comme pour $w_{2}$. Le lemme est d\'emontr\'e.
\qed\medskip

Les lemmes \ref{lem 29} et \ref{lem 30} montrent qu'il existe une
constante $C$ telle que, pour tout $x$
$$\bigg\vert \sum_{l} d_{l} e_{l} \displaystyle \int^1_{0} \theta(x -
\tau l)\ \tau^2 d\tau\bigg\vert \le \frac{1}{ 2}\vert \psi\ast
B(u,v) (1) (x)\vert + C \Vert d_{l}\Vert_{l^\infty}\ \Vert
e_{l}\Vert_{l^\infty}.$$
On pose, si $n$ est fix\'e
\begin{equation}\label{eq 47}
\varepsilon_{l}(x) = \textrm{sgn}\ \displaystyle \int^1_{0}\ \theta(x - \tau
l)\
\tau^2\ d\tau.
\end{equation}
Alors, si $e_{l} = \overline d_{l}\ \varepsilon_{l} (x)$, on obtient
l'in\'egalit\'e
\begin{equation}\label{eq 48}
\sum_{l} \vert d_{l}\vert^2\ \bigg\vert \displaystyle \int^1_{0} \theta(x -
\tau l)\ \tau^2 d\tau\bigg\vert\le \frac{1}{2}\vert \psi\ast B(u,v)
(1) (x)\vert + C \Vert d_{l}\Vert^2_{l^\infty}.
\end{equation}
Il s'agit maintenant d'en d\'eduire (\ref{eq 42}). Pour cela, on
choisit $\theta$ de mani\`ere appropri\'ee.\medskip

On rappelle que $\widehat \theta (\xi) = \widehat \psi (\xi)\ P(\xi)\
e^{-\vert \xi\vert^2}$, o\`u $\psi\in \Ss$ et $0\notin
\supp\widehat \psi$, et que $P(\xi)\not= 0$ pour tout $\xi$
tel que
$\sqrt{\xi^2_{1} + \xi^2_{2}}\le \delta\ \xi_{3}$, pour un certain
$\delta > 0$. Si $x\in \SR^3,\ x = (x_{1},x_{2},x_{3})$, on note
$x' = (x_{2},x_{3})$ et $x = (x_{1},x')$. On peut alors choisir
$\psi$ de sorte que
$$\theta(x) = \varphi(x_{1})\ \psi'(x'),$$
o\`u $\varphi\in \Ss(\SR)$ est \`a valeurs r\'eelles, $\varphi(0) >
0$ et $\displaystyle \int_{\SR}\varphi = 1$, et o\`u $\psi'\in
\Ss(\SR^2),\ \psi'(0)$ r\'eel $> 0$.\medskip

\begin{lem}\label{lem 31}
Il existe $A\ge 0,\ \alpha,\beta > 0$ tels que si\ $x = (x_{1},0)$ et
$l = (l_{1},l')$ v\'erifient les conditions $x_{1}\ge A,\ l_{1}\ge
4x_{1}$ et $d(x,[0,l])\le \alpha$, alors
$$\bigg\vert\displaystyle \int^1_{0}\ \theta(x - \tau l)\ \tau^2\
d\tau\bigg\vert\ge \beta\ \frac{x^2_{1}}{ l^3_{1}}.$$
\end{lem}
\medskip

\preuve
Soient $x = (x_{1},0),\ x_{1} > 0$, et $l\in \SZ^3$ tels que
$d(x,[0,l])\le \alpha$, o\`u $\alpha > 0$ est \`a choisir. Soit
$\tau_{0} l,\ 0\le \tau_{0}\le 1$, le point de $[0,l]$ en lequel
$d(x,[0,l])$ est atteinte, et $\gamma$ l'angle entre $(0,x)$ et
$(0,l)$. On a
$$(x_{1} - \tau_{0} l_{1})^2 + \vert \tau_{0} l'\vert^2\le \alpha^2,$$
et $\sin\gamma\le \frac{\alpha}{ x_{1}}$. Si $x_{1}\ge A\ge
2\alpha$, alors $0\le \gamma\le \frac{\pi}{ 6}$, d'o\`u
$\tau_{0}\frac{l_{1}}{ x_{1}} =\cos^2 \gamma\ge \frac{3}{4}$,
c'est-\`a-dire
$$x_{1} - 2\tau_{0} l_{1}\le -\frac{x_{1}}{ 2}.$$
Enfin, on a $\tau_{0}\le \frac{1}{ 4}$ d\`es que $l_{1}\ge
4x_{1}$.

On part alors de la d\'ecomposition
$$\displaystyle \int^1_{0}\ \theta(x - \tau l)\ \tau^2\ d\tau =
\int^{2\tau_{0}}_{0} \ldots + \int^1_{2\tau_{0}}\ldots .$$
On a d'une part
\begin{multline*}
\bigg\vert\displaystyle \int^{2\tau_{0}}_{0}\theta(x - \tau l)\tau^2
d\tau\bigg\vert\ge\vert \psi'(0)\vert
\bigg\vert\int^{2\tau_{0}}_{0}\varphi(x_{1} -
\tau l_{1})\tau^2
d\tau\bigg\vert\\
-2\alpha \Vert \nabla \psi'\Vert_{L^\infty} \displaystyle
\int^{2\tau_{0}}_{0}\vert\varphi(x_{1} - \tau l_{1})\vert\tau^2
d\tau,
\end{multline*}
et d'autre part
$$\bigg\vert\displaystyle \int^1_{2\tau_{0}}\theta(x - \tau l)\tau^2
d\tau\bigg\vert\le \Vert \psi'\Vert_{L^\infty}
\int^1_{2\tau_{0}}\vert\varphi(x_{1} - \tau l_{1})\vert \tau^2 d\tau.$$
On calcule les trois int\'egrales qui apparaissent dans ces
in\'egalit\'es :
$$\int^{2\tau_{0}}_{0}\varphi(x_{1} - \tau l_{1})\tau^2 d\tau =
\frac{1}{ l^3_{1}} \int^{x_{1}}_{x_{1} - 2\tau_{0}l_{1}}
(x_{1} - v)^2\
\varphi (v)\ dv,$$
de m\^eme si $\varphi$ est remplac\'e par $\vert \varphi \vert$, et
enfin
$$\int^1_{2\tau_{0}}\vert\varphi(x_{1} - \tau l_{1})\vert\tau^2 d\tau =
\frac{1}{ l^3_{1}} \int^{x_{1} - 2\tau_{0}l_{1}}_{x_{1} -
l_{1}} (x_{1} - v)^2\
\vert\varphi (v)\vert\ dv.$$
Puisque $x_{1} - 2\tau_{0}l_{1}\le -\frac{x_{1}}{ 2}$, on a
$$\displaystyle \int^{x_{1}}_{x_{1} - 2\tau_{0}l_{1}}\ \varphi (v)\
dv\ge \frac{1}{ 2}$$
d\`es que $x_{1}$ est assez grand. On en d\'eduit
$$\bigg\vert\displaystyle\int^{2\tau_{0}}_{0}\ \varphi(x_{1} - \tau l_{1})\
\tau^2\ d\tau\bigg\vert\ge \frac{1}{ 2}\ \frac{x_{1}^2}{ l^3_{1}} -
C\ \frac{1 + x_{1}}{l^3_{1}}$$
pour une certaine constante $C$ ne d\'ependant que de $\varphi$. On a
\'egalement
$$\displaystyle\int^{2\tau_{0}}_{0}\ \vert\varphi(x_{1} - \tau
l_{1})\vert\ \tau^2\ d\tau\le  C\ \frac{1 + x_{1} + x_{1}^2}{
l^3_{1}},$$ 
et
$$\displaystyle\int^1_{2\tau_{0}}\ \vert\varphi(x_{1} - \tau
l_{1})\vert\ \tau^2\ d\tau\le \frac{C}{l^3_{1}}\
\int^{-x_{1}/2}_{-\infty}\ v^2\ \vert\varphi(v)\vert\ dv.$$
Rassemblant les diverses in\'egalit\'es obtenues, il vient
\begin{multline*}
\left\vert\displaystyle \int^1_{0} \theta(x - \tau l)\tau^2
d\tau\right\vert\ge\displaystyle\frac{1}{ 2} \vert
\psi'(0)\vert
\displaystyle\frac{x^2_{1}}{ l^3_{1}}\\
-\displaystyle\frac{C}{ l^3_{1}}\left(1 + x_{1} + \alpha
x^2_{1} +
\displaystyle\int^{-x_{1}/2}_{-\infty} v^2\vert\varphi (v)\vert
dv\right),
\end{multline*}
si $x_{1}$ est assez grand. Le lemme en r\'esulte.
\qed\medskip

L'ensemble des lemmes pr\'ec\'edents, avec les in\'egalit\'es
(\ref{eq 47} - \ref{eq 48}), d\'emontre le lemme \ref{lem 25}. On
rel\`eve que n'importe quel $x = (x_{1},0)$ assez grand convient (\`a
un \'eventuel changement de rep\`ere pr\`es), et que $\Lambda$ est
l'ensemble des $l = (l_{1},l')$ tels que $l_{1}\ge 4x_{1}$ et
$d(x,[0,l])\le \alpha$, c'est-\`a-dire des points de $\SZ^3$ se
trouvant \`a l'int\'erieur d'un certain tronc de c\^one. Noter que
$\vert l\vert\le C\ l_{1}$, pour tout $l\in \Lambda$.\bigskip

\subsection{Preuve du lemme \ref{lem
26}}$\phantom{x}$\medskip

On fixe maintenant un tel $x$ et l'ensemble $\Lambda$ associ\'e.
Sans perdre de g\'en\'eralit\'e, on suppose $\alpha\le 1$, et
$x_{1} > 4(\alpha + 1)$. Enfin, on se place dans
l'hypoth\`ese o\`u $\displaystyle \sum_{n\ge 0}\
2^{-3n}\displaystyle\inf_{k\ge n}\ 2^{3k}\,
\eta_{k}^2 = +\infty$.\medskip

Si $j\ge 4$ et $l = (l_{1},l')\in \SZ^3$ est tel que $\vert
l'\vert\le \frac{\alpha}{ \sqrt 2}\ 2^j$ et $l_{1}\in 2^j\
x_{1} +
\left[-\frac{\alpha}{ \sqrt 2}\ 2^j, \frac{\alpha}{ \sqrt 2}\
2^j\right]$, on pose a priori $c_{l} = \delta_{j}$, pour un
certain $\delta_{j} > 0$. On note $\Lambda_{j}$ l'ensemble de
ces points $l$ : les $\Lambda_{j}$ sont deux \`a deux
disjoints et inclus dans $\Lambda$. Si $l\in \Lambda$ et
$l\notin \displaystyle \bigcup_{j\ge 4}\Lambda_{j}$, on pose
$c_{l} = 0$. Dans ces conditions, les s\'eries
$\displaystyle\sum_{l\in \Lambda}\ \frac{\vert c_{l}\vert^2}{
\vert l\vert^3}$ et $\displaystyle\sum_{j\ge 4}\
\delta^2_{j}$ sont de m\^eme nature.\medskip

Si $n\ge 0$, soit $\sigma_{n}=\displaystyle \inf_{k\ge n}\ 2^{3k}\
\eta_{k}^2$. On
choisit $\delta_{4}^2 = 2^{-12}\ \sigma_{4}$ et, si $j\ge 5,\
\delta^2_{j} = 2^{-3j}\ (\sigma_{j} - \sigma_{j - 1})$. On a alors par
construction les deux propri\'et\'es suivantes :
\begin{equation}\label{eq 49}
\displaystyle\sum_{j\ge 4}\ \delta^2_{j} = +\infty,
\end{equation}

\begin{equation}\label{eq 50}
\displaystyle\sum_{j\ge n}\ 2^{3(j - n)}\ \delta^2_{j}\le \eta^2_{n}\ \ \ \ \
\mbox{pour tout}\ n\ge 4.
\end{equation}
\medskip

Comme on a d\'ej\`a remarqu\'e que $\displaystyle\sum_{j\ge
4}\ \delta^2_{j}$ et $\displaystyle\sum_{l\in \Lambda}
\frac{\vert c_{l}\vert^2}{ \vert l\vert^3}$ sont de m\^eme
nature, (\ref{eq 49}) montre que (\ref{eq 43}) est
vrai.\medskip

Il reste \`a prouver (\ref{eq 44}) : on pose
$$a(y) = \sum_{l\in \Lambda}\ \varepsilon_{l}\ c_{l}\ m_{0} (y - l),$$
o\`u $\varepsilon_{l}\in \{-1,1\}$ pour tout $l$.\medskip

D'apr\`es le lemme \ref{lem 27}, la transform\'ee de Fourier de $a$
est support\'ee par une couronne compacte ne contenant pas 0.
L'appartenance de $a$ \`a $M(\eta)$ est donc \'equivalente \`a
l'existence d'une constante $C$ telle que pour tous $y_{0}\in \SR^3$
et $k\ge 0$
\begin{equation}\label{eq 51}
\displaystyle{\Mean_{B(y_{0},2^k)}}\ \vert a\vert^2\le C\
\eta^2_{k}.
\end{equation}
\medskip

Soient $y_{0}\in \SR^3$ et $k\ge 0$. Trois cas sont \`a consid\'erer,
le premier \'etant celui o\`u il existe $j\ge k + 1$ et $l_{0}\in
\Lambda_{j}$ tels que $l_{0}\in B(y_{0},2^{k + 1})$.\medskip

On v\'erifie d'abord, dans ce cas, que
$$\Lambda_{j'}\cap B(y_{0},2^{k + 1})= \emptyset$$
si $j'\not= j$. En effet, on sait que $\vert l_{0} - y_{0}\vert\le
2^{k + 1}$ et $\vert l_{0} - 2^j x\vert\le \alpha\ 2^j$, d'o\`u $\vert
y_{0} - 2^j x\vert\le (\alpha + 1) 2^j$. S'il existait $j'\not= j$ et
$l\in \Lambda_{j'}\cap B(y_{0},2^{k + 1})$, on aurait de m\^eme
$$\vert y_{0} - 2^{j'} x\vert\le (\alpha + 1)\ 2^{\max(j',k + 1)},$$
d'o\`u
$$\vert 2^j - 2^{j'}\vert\ \vert x\vert\le (\alpha + 1)\ (2^j +
2^{\max(j',k + 1)}),$$
puis $\vert x\vert\le 4 (\alpha + 1)$, ce qui est faux.\medskip

On \'ecrit alors
$$a(y) = \sum_{l\in \Lambda_{j}}\ \varepsilon_{l}\ c_{l}\ m_{0}(y -
l) + \sum_{l\notin \Lambda_{j}}\ \varepsilon_{l}\ c_{l}\ m_{0}(y - l).$$
La premi\`ere somme donne
\begin{align*}
\displaystyle\Mean_{B(y_{0},2^k)}\bigg\vert\sum_{l\in \Lambda_{j}}
\varepsilon_{l} c_{l} m_{0}(y - l)\bigg\vert^2
dy&\le\delta^2_{j}\displaystyle
\Mean_{B(y_{0},2^k)}\bigg(\sum_{l\in \Lambda_{j}}\vert m_{0}(y
- l)\vert\bigg)^2 dy\\
&\le C\ \delta^2_{j}\\
&\le C\ \eta^2_{k}
\end{align*}
d'apr\`es la d\'efinition de $\delta_{j}$ et la d\'ecroissance de $\eta$.
Pour traiter la
deuxi\`eme somme, on utilise la d\'ecroissance rapide de $m_{0}$. Si
$j'\not= j$, on a
$$\vert y_{0} - 2^{j'}x\vert\ge \vert (\alpha + 1) 2^j - (2^j -
2^{j'})x\vert\ge C \max (2^j,2^{j'}).$$
On en d\'eduit que $\vert y - l\vert\ge C 2^j$ si $y\in B(y_{0},2^k)$
et $l\in \displaystyle\bigcup_{j'\not= j}\ \Lambda_{j'}$, puis que
$$\bigg\vert\sum_{l\notin \Lambda_{j}} \varepsilon_{l} c_{l}\, m_{0}(y -
l)\bigg\vert\le C \sum_{l\notin \Lambda_{j}}\ \vert m_{0}(y - l)\vert\le
C_{N} 2^{-Nj}$$
pour tout $N\ge 0$.\medskip

La r\'egularit\'e de la suite $\eta$ implique l'existence de $N\ge 0$
tel que $\eta_{0} 2^{-Nj}\le \eta_{j}$ pour tout $j$. On a donc finalement
\begin{equation}\label{eq 52}
\displaystyle\Mean_{B(y_{0},2^k)}\bigg\vert\sum_{l\notin
\Lambda_{j}}
\varepsilon_{l} c_{l}\, m_{0}(y - l)\bigg\vert^2 dy\le \le C \eta^2_{j} \le C
\eta^2_{k},
\end{equation}
d'o\`u l'in\'egalit\'e (\ref{eq 51}) dans le premier cas.\medskip

On suppose maintenant qu'il existe $j_{0}\le k$ et $l_{0}\in
\Lambda_{j_0}$ tels que $l_{0}\in B(y_{0},2^{k + 1})$. On d\'ecompose
alors $a$ en
$$a(y) = \sum_{\substack{\scriptstyle l\in
\Lambda_{j}\\\scriptstyle j\le k}}
 \varepsilon_{l} c_{l}\, m_{0}(y - l) + \sum_{\substack{\scriptstyle
l\in
\Lambda_{j}\\\scriptstyle j\ge k+1}}\ \varepsilon_{l} c_{l}\, m_{0}(y
- l).$$
Pour la premi\`ere somme on \'ecrit
\begin{align*}
\displaystyle\Mean_{B(y_{0},2^k)}\bigg\vert\sum_{\substack{\scriptstyle
l\in
\Lambda_{j}\\\scriptstyle j\le k}}\varepsilon_{l} c_{l} m_{0}(y -
l)\bigg\vert^2 dy&\le C\ 2^{-3k}\displaystyle
\int_{\SR^3}\bigg\vert\sum_{\substack{\scriptstyle l\in
\Lambda_{j}\\\scriptstyle j\le k}}\varepsilon_{l} c_{l} m_{0}(y -
l)\bigg\vert^2 dy\\
&\le C\ 2^{-3k}\ \displaystyle\sum_{\scriptstyle l\in
\Lambda_{j},\atop\scriptstyle j\le k}\ \vert c_{l}\vert^2\\
&\le C\ 2^{-3k} \displaystyle\sum_{j\le k} 2^{3j}\,\delta^2_{j}\\
&\le C\ \eta^2_{k}
\end{align*}
d'apr\`es (\ref{eq 50}) et le lemme \ref{lem 27}. Pour la deuxi\`eme
somme, on proc\`ede comme au premier cas : si $y\in B(y_{0},2^k)$ et
$l\in \Lambda_{j},\ j\ge k + 1$, alors $\vert y - l\vert\ge C\ 2^k$,
et on conclut comme pour (\ref{eq 52}). Au total, on a obtenu (\ref{eq
51}).\medskip

Le troisi\`eme et dernier cas est celui o\`u aucun $\Lambda_{j}$
n'intersecte $B(y_{0},2^{k + 1})$. On a donc $\vert y - l\vert\ge 2^k$
pour tous $l\in \displaystyle\bigcup_{j\ge 4}\ \Lambda_{j}$ et $y\in
B(y_{0},2^k)$. On raisonne une fois encore comme pour (\ref{eq
52}).\medskip

On a ainsi prouv\'e (\ref{eq 51}) dans tous les cas, ce qui implique
(\ref{eq 44}) : le lemme \ref{lem 26} est enti\`erement prouv\'e, et la
d\'emonstration du th\'eor\`eme 11 est compl\`ete.



\section{Un cas non invariant :\ espaces 2-microlocaux}\bigskip

Si l'invariance de l'espace $E$, dans les th\'eor\`emes 8 et 9, est
une hypoth\`ese naturelle, elle n'est pas pour autant n\'ecessaire.
Il est par exemple possible de r\'esoudre (\ref{eq 1}) par la
m\'ethode KW lorsque $u_{0}$ est une distribution de norme assez
petite dans un espace 2-microlocal $\dot C^{-1,s'} (x_{0}),\ x_{0}$
fix\'e. Caract\'eris\'es par la condition
$$\vert \Delta_{j} f(x)\vert\le C\ 2^j (1 + 2^j \vert x -
x_{0}\vert)^{-s'},$$
pour tous $x\in \SR^3$ et $j\in \SZ$, ces espaces ne sont invariants ni
par translation ni par homoth\'etie de centre autre que $x_{0}$. Ce
r\'esultat, et sa g\'en\'eralisation \`a un cadre g\'eom\'etrique
plus vaste, font l'objet de cette section.\bigskip

\subsection{Bons espaces non invariants}$\phantom{x}$\medskip

Il faut d'abord \'etendre convenablement les th\'eor\`emes 8 et 9. On
observe pour cela que l'hypoth\`ese d'invariance de $E$ est
utilis\'ee, tout au long des preuves de ces th\'eor\`emes, de deux
fa\c cons. En premier lieu, la Proposition 13 est fr\'equemment
invoqu\'ee : il suffit, si $E$ n'est pas invariant, de supposer qu'il
est inclus dans $\dot B^{-1,\infty}_{\infty}$. En second lieu,
l'invariance par translation de $E$ entra\^{\i}ne que tout op\'erateur
de convolution par une fonction int\'egrable agit contin\^ument sur
$E$, et cela est appliqu\'e par exemple aux op\'erateurs $\Delta_{j}$
et $S_{j}$. Un examen attentif des preuves montre cependant que
l'hypoth\`ese suivante est suffisante :\par
\begin{equation}\label{eq 53}
\mbox{pour toute fonction}\ \phi\in \Ss(\SR^3),\ \mbox{il existe une
constante}\
C\ \mbox{telle}
\end{equation}
que, pour tous $t > 0$ et $f\in E,\ f\ast \phi_{t}\in E$, o\`u 
$\phi_{t} =
t^{-3}
\ \phi \left({\cdot\over t}\right)$, avec
$$\Vert f\ast \phi_{t}\Vert_{E}\le C\ \Vert f\Vert_{E}.$$

On est donc amen\'e \`a \'etendre la d\'efinition des bons
espaces.\medskip

\begin{defin}\label{defin 32}
Si $E$ est un espace fonctionnel non invariant, on dit que c'est un
bon espace non invariant lorsque
\begin{enumerate}
\item[a)] $E$ v\'erifie (\ref{eq 53}),
\item[b)] $E$ s'injecte contin\^ument dans $\dot
B^{-1,\infty}_{\infty}$,
\item[c)] $E$ v\'erifie (\ref{eq 18}) avec une suite $\eta\in l^1
(\SZ)$.
\end{enumerate}
\end{defin}\medskip

Il faut \'egalement g\'en\'eraliser la notion de couple admissible,
c'est-\`a-dire \'etendre au cadre non invariant les propri\'et\'es
(P1), (P2) et (P3) du paragraphe 1.2.\medskip

Si $E$ est un espace fonctionnel non invariant ayant les
propri\'et\'es a) et b) de la d\'efinition pr\'ec\'edente, on dira que
le couple $(E,F)$, o\`u $F$ est un espace de Banach, est admissible
lorsque les propri\'et\'es (P2') et (P3') suivantes sont
satisfaites.\medskip

Propri\'et\'e (P2')
\begin{enumerate}
\item[a)] $F$ s'injecte contin\^ument dans $\Ss'$,
\item[b)] $F$ v\'erifie (\ref{eq 53}),\par
\item[c)] pour tout $t > 0$, il existe une norme sur $F$, not\'ee
$\Vert \cdot
\Vert_{t,F}$, \'equivalente \`a la norme $\Vert \cdot\Vert_{F}$,
uniform\'ement sur tout compact de $]0,\infty[$, et il existe une
constante $C > 0$ telle que $$\Vert e^{t\Delta} u_{0}\Vert_{t,F}\le
C\Vert u_{0}\Vert_{E}$$ pour tout $u_{0}\in E$.
\end{enumerate}\medskip

Propri\'et\'e (P3')\par
Si $\Ff$ d\'esigne l'ensemble des fonctions $u$, d\'efinies sur
$]0,\infty[$ et continues \`a valeurs dans $F$, telles que $$\Vert
u\Vert_{\Ff} = \sup_{t > 0}\ \Vert u(t)\Vert_{t,F} < \infty,$$
alors
\begin{enumerate}
\item[a)] $B$ est continue de $\Ff\times \Ff$ dans $\Ff$,
\item[b)] $\displaystyle \lim_{t\rightarrow 0}\ B(u,v) (t) = 0$ dans
$S'$, pour tous $u,v\in \Ff$.

\end{enumerate}

Le r\'esultat suivant est l'extension \`a ce cadre des th\'eor\`emes 8
et 9, et se d\'emontre de la m\^eme fa\c con. On se contente de
l'\'enoncer.\medskip

\begin{theo}\label{theo 33}
Si $E$ est un bon espace non invariant, il existe un espace de Banach
$F$ formant avec lui un couple admissible. Si de plus
$\displaystyle\sum_{n\ge 0}\ n\eta_{n} < +\infty$, on peut choisir
$F\subset E$, et toute solution $u\in \Ff$ de (\ref{eq 1}) est
r\'eguli\`ere, au sens du th\'eor\`eme 9.
\end{theo}\bigskip

\subsection{Espaces 2-microlocaux g\'en\'eralis\'es}$\phantom{x}$\medskip

Afin d'autoriser des distributions aussi singuli\`eres que possible,
on g\'en\'eralise la d\'efinition des $\dot C^{-1,s'} (x_{0})$ en
rempla\c cant le singleton $\{x_{0}\}$ par un ensemble $S$, pour le
moment un ferm\'e quelconque de $\SR^3$.\medskip

\begin{defin}\label{defin 34}
Une distribution temp\'er\'ee $f$ appartient \`a $\dot C^{-1,s'}_{S},\
s' > 0$, si $\displaystyle\lim_{j\rightarrow -\infty} S_{j} f = 0$
dans $\Ss'$, et s'il existe une constante $C$ telle que, pour tous
$x\in \SR^3$ et $j\in \SZ$, on ait
$$\vert\Delta_{j} f(x)\vert\le C\, 2^j\, (1 + 2^j  d_{S}(x))^{-s'},$$
o\`u $d_{S}(x)$ d\'esigne la distance de $x$ \`a l'ensemble $S$.
\end{defin}
\medskip

Il n'est pas imm\'ediat que cette d\'efinition ne d\'epende pas du
choix des op\'erateurs $\Delta_{j}$. Cela provient du lemme suivant,
qu'on utilisera \`a plusieurs reprises.\medskip

\begin{lem}\label{lem 35}
Pour tout $N > 3 + s'$, il existe une constante $C$ telle que, pour
tous $x\in \SR^3$ et $j\in \SZ$
$$\displaystyle\int\ 2^{3j}(1 + 2^j \vert x - y\vert)^{-N}\ (1 + 2^j\
d_{S}(y))^{-s'} dy\le C\ (1 + 2^j  d_{S}(x))^{-s'}.$$
\end{lem}
\medskip

\preuve
Puisque $N > 3$, on a
\begin{align*}
&\displaystyle\int_{d_{S}(y)\ge {1\over2} d_{S}(x)}2^{3j} (1 + 2^j
\vert x - y\vert)^{-N}\ (1 + 2^j d_{S}(y))^{-s'} dy\\
&\qquad\le C\ (1 + 2^j d_{S}(x))^{-s'} \displaystyle\int_{d_{S}(y)\ge
{1\over2} d_{S}(x)}2^{3j} (1 + 2^j \vert x - y\vert)^{-N} dy\\
&\qquad\le C\ (1 + 2^j d_{S}(x))^{-s'}.
\end{align*}
D'autre part, si $d_{S}(y)\le {1\over 2}\ d_{S}(x)$, alors $\vert x -
y\vert\ge {1\over 2} d_{S}(x)$. Cette observation d\'ecoule de
l'in\'egalit\'e \'el\'ementaire
\begin{equation}\label{eq 54}
\forall\ x,y\in \SR^3\ \ \ \ \ d_{S}(x)\le \vert x - y\vert + d_{S}(y).
\end{equation}
On a donc
\begin{align*}
&\displaystyle\int_{d_{S}(y)\le {1\over2} d_{S}(x)}2^{3j} (1 + 2^j
\vert x - y\vert)^{-N}\ (1 + 2^j d_{S}(y))^{-s'} dy\\
&\qquad\le C\ (1 + 2^j d_{S}(x))^{-s'} \displaystyle\int_{d_{S}(y)\le
{1\over2} d_{S}(x)}2^{3j} (1 + 2^j \vert x - y\vert)^{-N + s'} dy\\
&\qquad\le C\ (1 + 2^j d_{S}(x))^{-s'},
\end{align*}
puisque $N - s' > 3$. Le lemme est d\'emontr\'e.
\qed\medskip

\begin{prop}\label{prop 36}
Pour tout $s' > 0,\ \dot C^{-1,s'}_{S}$ est un espace fonctionnel qui
s'injecte contin\^ument dans $\dot B^{-1,\infty}_{\infty}$ et qui
v\'erifie (\ref{eq 53}).
\end{prop}
\medskip

\preuve
L'injection continue de $\dot C^{-1,s'}_{S}$ dans $\Ss'$ ou dans $\dot
B^{-1,\infty}_{\infty}$ est \'evidente.

L'injection continue de $\Ss$ dans $\dot C^{-1,s'}_{S}$ est une
cons\'equence du lemme suivant.
\medskip

\begin{lem}\label{lem 37}
Pour tous $x_{0}\in \SR^3,\ N,p\ge 0$ et $f\in \Ss$, il existe une
constante $C$ telle que
\begin{align*}
\vert \Delta_{j} f(x)\vert&\le C\ 2^{3j}\, (1 + 2^j \vert x -
x_{0}\vert)^{-N},\quad j < 0,\\
\vert \Delta_{j} f(x)\vert &\le C\ 2^{-jp}\, (1 + 2^j \vert x -
x_{0}\vert)^{-N},\quad j\ge 0.
\end{align*}
\end{lem}\medskip

Ce lemme se d\'emontre directement dans le cas $j < 0$, et en
utilisant (\ref{eq 21}) dans le cas $j\ge 0$. Choisissant $x_{0}\in
S$, il implique $\Ss\subset \dot C^{-1,s'}_{S}$, quel que soit $s' >
0$, puisque $\vert x - x_{0}\vert\ge d_{S}(x)$ pour tout $x$. La
continuit\'e de l'injection est laiss\'ee au lecteur.

Soit $\Ee$ l'espace des distributions temp\'er\'ees $f$ telles
que $\displaystyle \lim_{j\rightarrow -\infty}\ S_{j} f = 0$ dans
$\Ss'$ et
$$\Vert f\Vert_{\Ee} = \sum_{j\in \SZ}\ 2^j \int
\vert\Delta_{j} f(x)\vert\, (1 + 2^j  d_{S}(x))^{-s'}\ dx < +\infty.$$
On v\'erifie que l'espace $\Ee$ est complet, et que $\Ss$ est
dense dans $\Ee$. Soit $\lambda$ une forme lin\'eaire continue
sur $\Ee$ et $\psi\in \Ss$ telle que $\supp \widehat\psi\subset
\Gamma_{0}$. On a pour tout $x$
\begin{align*}
\vert(\lambda,2^{3j}\psi (2^jx - 2^j\cdot))\vert&\le C\ \Vert
2^{3j}\psi (2^jx - 2^j.)\Vert_{\Ee}\\
&\le C\ 2^{4j} \displaystyle\int \vert\psi(2^jx - 2^jy)\vert (1 +
2^j d_{S}(y))^{-s'}\ dy\\
&\le C\ 2^j (1 + 2^j d_{S}(x))^{-s'},
\end{align*}
d'apr\`es le lemme \ref{lem 35}. On en d\'eduit que le dual de
$\Ee$ est $\dot C^{-1,s'}_{S}$, et que ce dernier est un
espace fonctionnel.\medskip

Il reste \`a prouver (\ref{eq 53}). Si $\phi\in \Ss$ et $t > 0$, on a
$\vert\Delta_{j}\ \phi(x)\vert\le C_{N}\ 2^{3j} (1 + 2^j\vert x
\vert)^{-N}$ pour tout $N\ge 0$, d'apr\`es le lemme \ref{lem 37}, et
par cons\'equent
$$\vert\Delta_{j} \phi_{t}(x)\vert\le C_{N}\ 2^{3j} (1 + 2^j \vert
x\vert)^{-N}$$
uniform\'ement par rapport \`a $t > 0$.\medskip

Si $f\in \dot C^{-1,s'}_{S}$ et $j\in \SZ$, il r\'esulte de (\ref{eq
15}) que
$$\Delta_{j} (\phi_{t}\ast f) = \widetilde\Delta_{j} \phi_{t}\ast
\Delta_{j} f,$$
puis que
\begin{align*}
&\vert \Delta_{j} (\phi_{t}\ast f) (x)\vert\\
&
\qquad\le
C\ \Vert f\Vert_{\dot
C^{-1,s'}_{S}}\!\! \displaystyle\int 2^{3j} (1 + 2^j \vert x - y
\vert)^{-N}\, 2^j (1 + 2^j d_{S}(y))^{-s'}dy\\
&\qquad\le C\ \Vert f\Vert_{\dot C^{-1,s'}_{S}}\ 2^j (1 + 2^j
d_{S}(x))^{-s'},
\end{align*}
en utilisant une nouvelle fois le lemme \ref{lem 35}. Ceci prouve que
$$\Vert \phi_{t}\ast f\Vert_{\dot C^{-1,s'}_{S}}\le C\ \Vert f\Vert_{\dot
C^{-1,s'}_{S}},$$
uniform\'ement par rapport \`a $t$. La Proposition \ref{prop 36} est
d\'emontr\'ee.\qed \medskip

\noindent\underbar{{\bf Exemples}} : notant $x_{1},x_{2},x_{3}$ les
coordonn\'ees de $x\in \SR^3$, les distributions v.p. ${1\over x_{1}}$ et
${1\over (x_{1}^2 +
x_{2}^2)^{1/2}}$ appartiennent respectivement
\`a $\dot C^{-1,s'}_{P},\ P = \{x\in \SR^3\ ;\ x_{1} = 0\}$ et $\dot
C^{-1,s'}_{D},\ D = \{x\in \SR^3\ ;\ x_{1} = x_{2} = 0\}$, pour tout
$s' > 0$.\bigskip

\subsection{Fonction de densit\'e et bons espaces
2-microlocaux}$\phantom{x}$\medskip

Pour \'etudier la validit\'e de (\ref{eq 18}) dans les espaces $\dot
C^{-1,s'}_{S}$, on introduit ce qu'on appelle la fonction de
densit\'e de $S$.\medskip

\begin{defin}\label{defin 38}
Si $S$ est un ferm\'e de $\SR^3$, sa fonction de densit\'e est la
fonction $\varepsilon_{S}$, d\'efinie sur $[0,1]$ par
$$\varepsilon_{S}(\delta) = \sup{1\over \vert
B(x,r)\vert}\left\vert\left\{y\in B(x,r)\ ;\ d_{S}(y)\le \delta
r\right\}\right\vert,$$
o\`u le supremum est pris sur tous les $x\in \SR^3$ et $r > 0$.
\end{defin}
\noindent Noter que $\varepsilon_{S}$ est croissante.
\medskip

\begin{prop}\label{prop 39}
Pour tout $s' > 0$ et tout ferm\'e $S$ de $\SR^3$, l'espace $\dot
C^{-1,s'}_{S}$ v\'erifie (\ref{eq 18}) avec, si $n\ge 0$,
$$\eta_{n} = C\ \sum_{m = 0}^n\ 2^{-2 s'(n - m)}\
\varepsilon_{S}(2^{-m})$$
pour une certaine constante $C$.
\end{prop}\medskip

\preuve
Soient $f,g\in \dot C^{-1,s'}_{S}$, $ \supp\widehat f\subset
\Gamma_{k}$, $\supp\widehat g\subset \Gamma_{l}$, $k,l\in \SZ$.

Si $\vert k - l\vert\ge 3$ et par exemple $l\le k - 3$, on a
$$\vert f(y)\ g(y)\vert\le 2^l\ 2^k\ (1 + 2^k d_{S}(y))^{-s'}\ \Vert
f\Vert_{\dot C^{-1,s'}_{S}}\ \Vert g \Vert_{\dot C^{-1,s'}_{S}}$$
pour tout $y$. Le lemme \ref{lem 35} montre alors que, si $\vert k -
j\vert\le 2$,
$$\vert \Delta_{j} (fg) (x)\vert\le C\ 2^l\, 2^j\, (1 + 2^j
d_{S}(y))^{-s'}\
\Vert
f\Vert_{\dot C^{-1,s'}_{S}}\ \Vert g \Vert_{\dot C^{-1,s'}_{S}}.$$
Ceci implique (\ref{eq 19}).

Si $\vert k - l\vert\le 2$ et si on se donne $j\le k + 4$, alors
$$\vert f(y)\ g(y)\vert\le 4^k\ (1 + 2^k d_{S}(y))^{-2s'}\ \Vert
f\Vert_{\dot C^{-1,s'}_{S}}\ \Vert g \Vert_{\dot C^{-1,s'}_{S}},$$
d'o\`u
$$\vert \Delta_{j} (fg) (x)\vert\le C\ I_{j,k}(x)\ 4^k\ \Vert
f\Vert_{\dot C^{-1,s'}_{S}}\ \Vert g \Vert_{\dot C^{-1,s'}_{S}},$$
avec
$$I_{j,k}(x) = \displaystyle\int\ 2^{3j}\ (1 + 2^j\vert x -
y\vert)^{-N}\ (1 + 2^k d_{S}(y))^{-2s'}\ dy,$$
$N$ d\'esignant un param\`etre \`a choisir. La Proposition \ref{prop
39} r\'esulte directement de l'in\'egalit\'e
\begin{equation}\label{eq 55}
I_{j,k}(x)\le C\ \displaystyle\sum^{k - j}_{m = 0}\ 2^{-2s'(k - j -
m)}\ \varepsilon_{S}(2^{-m})\ (1 + 2^j d_{S}(x))^{-2s'},
\end{equation}
valable lorsque $j\le k$ et sous la condition $N > 3 + 2s'$ (si $k +
1\le j\le k + 4$, le lemme \ref{lem 35} s'applique).
\medskip

On aura besoin du r\'esultat suivant.\medskip

\begin{lem}\label{lem 40}
Pour tout $N > 3$ il existe une constante $C$ telle que, pour tous
$x\in \SR^3,\ j\in \SZ$ et $\delta\in [0,1]$, on ait
$$\displaystyle\int_{d_{S}(y)\le \delta 2^{-j}}\ 2^{3j}\ (1 + 2^j
\vert x - y\vert)^{-N}\ dy\le C\ \varepsilon_{S} (\delta).$$
\end{lem}\medskip

\preuve
Il suffit de d\'ecouper l'int\'egrale de mani\`ere appropri\'ee :
\begin{align*}
&\displaystyle\int_{d_{S}(y)\le \delta 2^{-j}}\ 2^{3j}\ (1 + 2^j
\vert x - y\vert)^{-N}\ dy\\
&\quad\le\ \displaystyle\int_{d_{S}(y)\le \delta 2^{-j}, \vert x -
y\vert\le 2^{-j}} \ldots +
\displaystyle\sum_{m\ge 1} \int_{d_{S}(y)\le \delta 2^{-j}, 2^{m -
j - 1}\le \vert x - y\vert\le 2^{m - j}} \ldots\\
&\quad\le\ C \displaystyle\sum_{m\ge 0} 2^{3j} 2^{-mN}\vert\{y\in
B(x,2^{m - j}) ; d_{S}(y)\le \delta 2^{-j}\}\vert\\
&\quad\le\ C \displaystyle\sum_{m\ge 0}\ 2^{-m(N - 3)}\ \varepsilon_{S}
(2^{-m}\ \delta) \\
&\quad\le\ C\, \varepsilon_{S} (\delta).
\end{align*}
\qed\medskip

Revenant \`a la preuve de (\ref{eq 55}), on distingue les cas
$d_{S}(x)\ge 2^{-j + 1}$ et $d_{S}(x) < 2^{-j + 1}$.\medskip

Si $d_{S}(x) < 2^{-j + 1}$, on d\'ecoupe $I_{j,k}(x)$ selon
$$I_{j,k}(x) = \displaystyle\int_{d_{S}(y)\ge 2^{-j}} +
\int_{d_{S}(y)\le 2^{-j}}.$$
D\`es que $N > 3$ la premi\`ere int\'egrale est major\'ee par $C\
2^{-2s'(k - j)}$. Pour la seconde, on utilise le lemme \ref{lem 40} :
\begin{align*}
&\displaystyle\int_{d_{S}(y)\le 2^{-j}} 2^{3j}(1 +
2^j\vert x - y\vert)^{-N} (1 + 2^k d_{S}(y))^{-2s'}\ dy\\
&\quad = \displaystyle\sum_{m = 0}^{k - j - 1}\int_{2^{-m-j-1}\le
d_{S}(y)\le 2^{-m-j}}\ldots + \int_{d_{S}(y)\le 2^{-k}}\ldots\\
&\quad\le\ C \displaystyle\sum^{k-j-1}_{m=0} 2^{-2s'(k-j-m)}
\int_{d_{S}(y)\le 2^{-m-j}} 2^{3j}(1 + 2^j\vert x - y\vert)^{-N}\ dy\\
&\quad \quad+ \displaystyle\int_{d_{S}(y)\le 2^{-k}} 2^{3j}(1 + 2^j\vert
x - y\vert)^{-N}\ dy\\
&\quad\le\ C \displaystyle\sum^{k-j-1}_{m=0} 2^{-2s'(k-j-m)}
\varepsilon_{S}(2^{-m}) + C \varepsilon_{S}(2^{j - k}).
\end{align*}
Au total, on obtient bien (\ref{eq 55}) dans ce cas.\medskip

Lorsque $d_{S}(x)\ge 2^{-j+1}$, on part du d\'ecoupage
\begin{align*}
I_{j,k}(x) &= \displaystyle\int_{d_{S}(y)\ge {1\over 2}d_{S}(x)}
+
\int_{2^{-j}\le d_{S}(y)\le {1\over 2}d_{S}(x)}  +
\int_{d_{S}(y)\le 2^{-j}}\\
& = I + II + III.
\end{align*}
Quand $d_{S}(y)\ge {1\over 2}d_{S}(x)$, on a
\begin{align*}
(1 + 2^k d_{S}(y))^{-2s'}&\le C (1 + 2^k d_{S}(x))^{-2s'}\\ &\le C
2^{-2s' (k-j)} (1 + 2^j d_{S}(x))^{-2s'},
\end{align*}
puisque $k\ge j$. Cela donne (d\`es que $N > 3$)
\begin{equation}\label{eq 56}
I\le C\ 2^{-2s' (k-j)}\ (1 + 2^j d_{S}(x))^{-2s'}.
\end{equation}
Quand $d_{S}(y)\le{1\over 2} d_{S}(x)$, on a $\vert x - y\vert\ge
{1\over 2} d_{S} (x)$, en vertu de (\ref{eq 54}) d'o\`u
$$(1+2^j\vert x - y\vert)^{-N} \le (1+2^jd_{S}(x))^{-2s'} (1+2^j\vert x
- y\vert)^{-N+2s'}.
$$
 Si de plus
$d_{S}(y)\ge 2^{-j}$, alors on peut majorer $(1 + 2^k
d_{S}(y))^{-2s'}$ par $2^{-2s' (k - j)}$, ce qui permet d'obtenir
\begin{equation}\label{eq 57}
II\le C\ 2^{-2s'(k-j)}\ (1 + 2^j  d_{S}(x))^{-2s'}.
\end{equation}
Enfin, on majore $III$  par 
$$ C\ (1 + 2^j d_{S}(x))^{-2s'}
\displaystyle\int_{d_{S}(y)\le 2^{-j}} 2^{3j}(1 +
2^j \vert x - y\vert)^{-N+2s'} (1 + 2^k d_{S}(y))^{-2s'} dy.$$
On estime alors l'int\'egrale du membre de droite comme au premier
cas, ce qui donne
$$III\le C\ \displaystyle\sum^{k-j}_{m = 0}\ 2^{-2s' (k-j-m)}\
\varepsilon_{S}(2^{-m})\ (1 + 2^j d_{S}(x))^{-2s'}.$$
Avec (\ref{eq 56} - \ref{eq 57}), cela prouve (\ref{eq 55}) dans le
deuxi\`eme cas \'egalement, et ach\`eve la d\'emonstration de la
Proposition \ref{prop 39}.
\qed\medskip

\noindent Des Propositions \ref{prop 36} et \ref{prop 39} d\'ecoule le
r\'esultat suivant.\medskip

\begin{theo}\label{theo 41}
Si $S$ est un ferm\'e de $\SR^3$ dont la fonction de densit\'e
$\varepsilon_{S}$ v\'erifie la condition de Dini
$\displaystyle\int^1_{0}\varepsilon_{S} (\delta)\ \frac{d\delta}{
\delta} < +\infty$, alors $\dot C^{-1,s'}_{S}$ est un bon espace non
invariant, pour tout $s' > 0$.
\end{theo}\medskip

\preuve
Il suffit de v\'erifier que la suite $\eta_{n}$ donn\'ee par la
Proposition \ref{prop 39} est sommable.
\qed\bigskip

\subsection{R\'egularit\'e des solutions}$\phantom{x}$\medskip

Appliquant le th\'eor\`eme \ref{theo 33}, on peut r\'esoudre (\ref{eq
1}) pour toute donn\'ee $u_{0}$ assez petite dans $\dot C^{-1,s'}_{S},\
s' > 0$, d\`es que $\varepsilon_{S}$ v\'erifie la condition de Dini,
ce qu'on suppose dans ce paragraphe.\medskip

On fixe $s' > 0$ et, notant $E = \dot C^{-1,s'}_{S}$, on choisit $N$
de sorte que le couple $(E,F)$ soit admissible, avec $F =
C^{N,\infty}_{E}$ (voir la d\'emonstration du th\'eor\`eme 8). Soit
$u\in \Ff$ une solution de (\ref{eq 1}), associ\'ee \`a $u_{0}\in E$.
Le but de ce paragraphe est l'\'etude de la r\'egularit\'e de
$u$.\medskip

L'appartenance de $u$ \`a $\Ff$ signifie exactement que
\begin{equation}\label{eq 58}
\vert \Delta_{j} u(t,x)\vert\le C\ 2^j\ (1 + 2^j \sqrt t)^{-N}\ (1 +
2^j d_{S}(x))^{-s'},
\end{equation}
pour tout $t > 0,\ x\in \SR^3$ et $j\in \SZ$. On voit que, comme pour
l'\'equation de la chaleur, la solution est instantan\'ement
r\'egularis\'ee. En particulier, $u(t)$ est, pour chaque $t > 0$, une
fonction born\'ee.\medskip

\begin{prop}\label{prop 42}
On a
\begin{enumerate}
\item[a)] $\vert u(t,x)\vert\le \frac{C}{\sqrt t + d_{S}(x)}$\ \ \
si\
\ $s' > 1$,
\item[b)] $\vert u(t,x)\vert\le \frac{C}{\sqrt t + d_{S}(x)}\ln
\left(10 + {d_{S}(x)\over \sqrt t}\right)$\ \ \ si\ \ $s' = 1$,
\item[c)] $\vert u(t,x)\vert\le C\ \frac{1}{ (\sqrt t)^{1-s'}}\
\frac{C}{\sqrt t + d_{S}(x)}$\ \ \ si\ \ $0 < s' < 1$.
\end{enumerate}
\end{prop}\medskip

\preuve
Il suffit d'injecter (\ref{eq 58}) dans
$$u(t,x) = \sum_{j\in \SZ}\ \Delta_{j} u(t,x).$$
Les calculs sont laiss\'es au lecteur.
\qed\medskip

Dans le cas le plus r\'egulier, c'est-\`a-dire quand $s' > 1$, on a
aussi $\vert u_{0}(x)\vert\le \frac{C}{ d_{S}(x)}$. Ainsi, \`a chaque
instant $t > 0$, la solution $u$ est born\'ee par $\frac{C}{ \sqrt t}$
dans un voisinage de $S$ d'\'epaisseur $\sqrt t$, et retrouve le
comportement en ${1\over d_{S}}$ de la donn\'ee initiale en-dehors de
ce voisinage.\medskip

On peut aller un peu plus loin en remarquant qu'on a
$$\Vert fg\Vert_{\dot C^{-1,s'}_{S}}\le \Vert f\Vert_{\dot C^{-1,s'}_{S}}\
\Vert g\Vert_{L^\infty}$$
si $f,g\in \dot C^{-1,s'}_{S}$, $\supp \widehat f\subset \Gamma_{k}$, $
\supp \widehat g\subset \Gamma_{l}$, et $l\le k - 3$ (reprendre le
d\'ebut de la preuve de la Proposition \ref{prop 39}). Si, par une
d\'emarche semblable \`a celle de la Proposition 16, on note
momentan\'ement $\Ff_{N}$ au lieu de $\Ff$ l'espace des $u$
v\'erifiant (\ref{eq 58}), on montre alors que l'application $B$ est
continue de $\Ff_{N}\times \Ff_{N}$ dans $\Ff_{N + 1}$, pour tout $N$
assez grand. Autrement dit, la solution $u$ v\'erifie (\ref{eq 58})
pour tout $N$. Elle est donc de classe $C^{\infty}$ sur
$]0,\infty[\times \SR^3$, et ses d\'eriv\'ees satisfont des
estimations ponctuelles analogues \`a celle de la Proposition
\ref{prop 42}. Par exemple, on a dans le cas $s' > 1$ :
$$\left\vert{\partial^n\over \partial t^n}\ u(t,x)\right\vert + \vert
D^\alpha_{x}\ u(t,x)\vert\le C\ (\sqrt t + d_{S}(x))^{-1-2n}$$
si $\vert \alpha\vert = 2n < s' - 1$, et
$$\left\vert{\partial^n\over \partial t^n}\ u(t,x)\right\vert + \vert
D^\alpha_{x}\ u(t,x)\vert\le C\ (\sqrt t + d_{S}(x))^{-s'}\ (\sqrt
t)^{-1 - 2n +s'}$$
si $\vert \alpha\vert = 2n > s' - 1$.\medskip

\noindent{{\bf Remarque :}}  c'est parce que l'espace $E =
\dot C^{-1,s'}_{S}$ est lui-m\^eme de type Besov $l^\infty$,
c'est-\`a-dire coinc\"{\i}de avec $\dot B^{0,\infty}_{E}$, que la
condition $\displaystyle\sum_{n\ge 0} n\eta_{n} < +\infty$
n'appara\^{\i}t pas, bien que la solution $u$ soit r\'eguli\`ere. Dans
un tel cas, les espaces $F$ et $G$ utilis\'es respectivement aux
th\'eor\`emes 8 et 9 coinc\"{\i}dent, et il n'y a donc pas lieu de
renforcer la condition de sommabilit\'e de la suite $\eta$.\medskip

\noindent \underbar{{\bf Exemples}} : les solutions de (\ref{eq 1})
obtenues \`a partir de $u_{0}(x) = \varepsilon\ v.p.\frac {1}{x_{1}}$,
ou de $v_{0}(x) = \varepsilon\ \frac{1}{ \sqrt{x^2_{1} + x^2_{2}}}$, pour
$\varepsilon$ assez petit, sont de classe $C^\infty$ sur
$]0,\infty[\times \SR^3$ et v\'erifient toutes les estimations
pr\'ec\'edentes. Ces solutions sont auto-similaires, en vertu de
l'homog\'en\'e\"{i}t\'e de $u_{0}, v_{0}$. Voir aussi ([Ko,Y]), o\`u ce
type de donn\'ees initiales est consid\'er\'e.\bigskip

\subsection{Comparaison du terme lin\'eaire et du terme
bilin\'eaire}$\phantom{x}$\medskip

On reste sous les hypoth\`eses de la section 5.4, et on
consid\`ere $u$ solution de (\ref{eq 1}) associ\'ee \`a $u_{0}\in \dot
C^{-1,s'}_{S}$.\medskip

Le terme lin\'eaire $S u_{0}$ satisfait aux estimations (\ref{eq 58})
pour tout $N$, mais sans qu'il soit possible a priori d'am\'eliorer
l'exposant $s'$. De ce point de vue, le terme bilin\'eaire $B(u,u)$
est plus r\'egulier.\medskip

\begin{lem}\label{lem 43}
En posant $w = B(u,u)$, on a
\begin{equation}\label{eq 59}
\vert \Delta_{j} w (t,x)\vert\le C_{N}\ 2^j (1 + 2^j\sqrt t)^{-N}\ (1 +
2^j d_{S}(x))^{-\sigma}\ L(2^j d_{S}(x))
\end{equation}
pour tout $N$, o\`u $\sigma = min (2s', s' + 1)$, et $L(r) = 1$ si
$s'\not= 1,\ L(r) = \ln (10 + r)$ si $s' = 1$.
\end{lem}\medskip

\preuve
On note $E^\sharp$ l'espace $\dot C^{-1,\sigma}_{S}$ si $s'\not= 1$, et si
$s' =
1$, l'espace des distributions temp\'er\'ees $f$ telles que
$\displaystyle \lim_{j\rightarrow -\infty}\ S_{j} f = 0$ dans $\Ss'$
et v\'erifient
$$\vert \Delta_{j} f(x)\vert\le C\ 2^j (1 + 2^j d_{S}(x))^{-2}\ln
(10 + 2^j  d_{S}(x))$$
pour tous $j\in \SZ,\ x\in \SR^3$. La norme de $f$ dans $E^\sharp$ est la
meilleure constante $C$ possible.\medskip

Soit, pour tout $N,\ F_{N}^\sharp = C^{N,\infty}_{E^\sharp}$ l'espace
construit
au-dessus de $E^\sharp$, et $\Ff^\sharp_{N}$ l'espace de distributions
d\'efinies sur $]0,\infty[\times \SR^3$ qui lui est associ\'e. On
sait que $u\in \Ff_{N}$ pour tout $N$ (section 5.4), de sorte que
le lemme \ref{lem 43} d\'ecoule de l'in\'egalit\'e
\begin{equation}\label{eq 60}
\Vert B(u,v)\Vert_{\Ff^\sharp_{N}}\le C_{N}\ \Vert u\Vert_{\Ff_{N}}\ \Vert
v\Vert_{\Ff_{N}}.
\end{equation}
\medskip

Celle-ci se d\'emontre en revenant \`a la preuve du th\'eor\`eme 8.
Si $u,v\in \Ff_{N}$, on estime les deux termes
$$R_{j}(t)^\sharp = 2^j \displaystyle \int^t_{0} \big(1 + 2^j \sqrt{t -
\tau}\big)^{-p} \Vert \Delta_{j} u(\tau)\ S_{j - 2}\
v(\tau)\Vert_{E^\sharp} d\tau,$$
$$C_{j}(t)^\sharp = 2^j \displaystyle \int^t_{0} \big(1 + 2^j \sqrt{t -
\tau}\big)^{-p} \bigg\Vert \Delta_{j} \bigg( \sum_{k\ge j - 4}\
\Delta_{k} u(\tau)\
\widetilde \Delta_{k} v(\tau)\bigg)\bigg\Vert_{E^\sharp}\ d\tau,$$
pour tout $t > 0$ et pour $p$ assez grand.\medskip

Le terme rectangle se traite en partant de
$$S_{j - 2} v(\tau,y) = \sum_{j'\le j - 3}\ \Delta_{j'}v(\tau,y),$$
d'o\`u
\begin{align*}
\vert S_{j - 2}\ v(\tau,y)\vert&\le \displaystyle\sum_{j'\le j - 3}\
2^{j'} (1 +
2^{j'} d_{S}(y))^{-s'}\\
&\le C2^j (1 + 2^j d_{S}(y))^{-\nu}\ L(2^j d_{S}(y)),
\end{align*}
o\`u $\nu = \min (s',1)$. On en d\'eduit
$$\vert \Delta_{j} u(\tau,y) S_{j - 2}\ v(\tau,y)\vert\le C\ 4^j (1 +
2^{j} d_{S}(y))^{-\sigma}\, L(2^j d_{S}(y))\ \big(1 + 2^j \sqrt
\tau\big)^{-N},$$ puis
$$\Vert \Delta_{j} u(\tau) S_{j - 2}\ v(\tau)\Vert_{E^\sharp}\le C\ 2^j 
\big(1 +
2^{j} \sqrt \tau\big)^{-N}$$
en utilisant le lemme \ref{lem 35}, convenablement modifi\'e si $s' =
1$. On en d\'eduit que $R_{j}(t)^\sharp$ v\'erifie (\ref{eq 25}),
tout comme $R_{j}(t)$.\medskip

Pour le terme carr\'e, on remarque que l'in\'egalit\'e (\ref{eq 55})
implique
$$\Vert \Delta_{j} (\Delta_{k}\ u(\tau) \widetilde\Delta_{k}\
v(\tau))\Vert_{E^\sharp}\le \eta_{k - j}\ 4^k\ 2^{-j} \big(1 + 2^{k} \sqrt
\tau\big)^{-2N},$$
d'o\`u on d\'eduit que $C_{j}(t)^\sharp$ v\'erifie (\ref{eq 26} -
\ref{eq 27}) comme au th\'eor\`eme 8. L'in\'egalit\'e (\ref{eq 60}) en
r\'esulte, ce qui prouve le lemme.\qed
\eject

En r\'esum\'e, on a obtenu le
\begin{theo}\label{theo 44}
Soit $S$ un ferm\'e de $\SR^3$ dont la fonction de densit\'e ob\'eit
\`a la condition de Dini. Pour tout $s' > 0$ il existe $\alpha > 0$
tel que, si $u_{0}\in \dot C^{-1,s'}_{S}$ est de norme au plus
$\alpha$, alors l'\'equation (\ref{eq 1}) admet une solution $u$.
Celle-ci v\'erifie (\ref{eq 58}) pour tout $N$, est de classe
$C^\infty$ sur $]0,\infty[\times \SR^3$, et converge vers $u_{0}$
dans $\dot C^{-1,s'}_{S}$ quand $t\rightarrow 0$ pour la topologie
 faible $\ast$. De plus, dans la d\'ecomposition $u = Su_{0} +
B(u,u)$, le terme $B(u,u)$ v\'erifie (\ref{eq 59}) pour tout $N$.
\end{theo}

\medskip

\noindent {\bf Remerciements. } {\sl Merci \`a Mireille Berg pour sa
frappe efficace de ce long manuscrit.}




\end{document}